\documentclass[10pt,a4paper,twoside]{article}
\usepackage{amssymb,amsmath,amsthm,amsfonts,amscd}
\usepackage{graphicx,xcolor}
\usepackage[unicode,psdextra]{hyperref}
\hypersetup{pdfborder={0 0 0.06},citebordercolor=[rgb]{0.8196,0.2275,0.5098},linkbordercolor=[rgb]{0.1765,0.5490,0.8118},urlbordercolor=[rgb]{0.7059,0.5333,0.1137}}
\usepackage{tikz-cd}
\usepackage[english]{babel}
\usepackage{dsfont}
\usepackage[longnamesfirst,numbers,sort&compress]{natbib}
\usepackage{enumerate}
\usepackage{setspace}
\usepackage{changepage}
\usepackage{booktabs}
\usepackage[sort&compress,capitalise,nameinlink]{cleveref}
\crefname{section}{\textsection}{\textsection}
\crefname{subsection}{\textsection}{\textsection}
\crefname{subsubsection}{\textsection}{\textsection}
\crefname{paragraph}{\textparagraph}{\textparagraph}
\crefname{thde}{Theorem}{Theorems}
\crefname{prop}{Pr.}{Pr.}
\crefname{thm}{Theorem}{Theorems}
\usepackage{csquotes}
\usepackage[cal=euler,scr=boondoxo,bb=ams,calscaled=1.07,scrscaled=1.07]{mathalfa}
\makeatletter
\def\mathalfa@frakscaled{s*[1]}
\makeatother
\DeclareFontFamily{U}{euf}{}%
\DeclareFontShape{U}{euf}{m}{n}{<-7>\mathalfa@frakscaled eufm5
  <7-9>\mathalfa@frakscaled eufm7
  <9->\mathalfa@frakscaled eufm10}{}%
\DeclareFontShape{U}{euf}{b}{n}{<-7>\mathalfa@frakscaled eufb5
  <7-9>\mathalfa@frakscaled eufb7
  <9->\mathalfa@frakscaled eufb10}{}%
\DeclareMathAlphabet{\mathfrak}{U}{euf}{m}{n}

\DeclareMathOperator{\Tr}{Tr}

\DeclareMathOperator{\ran}{Ran}

\newtheorem{thm}{Theorem}[section]
\newtheorem{proposition}[thm]{Proposition}
\newtheorem{prop}[thm]{Proposition}
\newtheorem{thde}[thm]{Theorem-Definition}
\newtheorem*{proposition*}{Proposition}
\newtheorem{lemma}[thm]{Lemma}
\newtheorem{corollary}[thm]{Corollary}
\newtheorem{definition}[thm]{Definition}

\theoremstyle{definition}

\newtheorem{remark}[thm]{Remark}
\newtheorem*{remark*}{Remark}

\newcommand{\im}{\mathrm{Im}}
\newcommand{\re}{\mathrm{Re}}

\def\YEAR{\year}\newcount\VOL\VOL=\YEAR\advance\VOL by-1995
\def\firstpage{1}\def\lastpage{1000}
 \def\received{}\def\revised{}
 \def\communicated{}

\makeatletter
\def\magnification{\afterassignment\m@g\count@}
\def\m@g{\mag=\count@\hsize6.5truein\vsize8.9truein\dimen\footins8truein}
\makeatother

\font\eightrm=cmr8
\font\caps=cmcsc10                    
\font\Caps=cmcsc10 scaled \magstep1   
\font\scaps=cmcsc8

\def\bfseries{\normalsize\caps}

\makeatletter
\@specialpagefalse\headheight=8.5pt
\def\DocMath{{\def\th{\thinspace}\scaps Documenta Math.}}
\renewcommand{\@oddfoot}{\hfill\scaps Final version, to appear in \DocMath\hfill}
\renewcommand{\@evenfoot}{\ifnum\thepage>\lastpage\hfill\scaps
    Final version, to appear in \DocMath\hfill\else\@oddfoot\fi}%

\renewcommand{\@evenhead}{%
    \ifnum\thepage>\lastpage\rlap{\thepage}\hfill%
    \else\rlap{\thepage}\slshape\leftmark\hfill{\caps\SAuthor}\hfill\fi}%
\renewcommand{\@oddhead}{%
    \ifnum\thepage=\firstpage{$\phantom{i}$\hfill\llap{\thepage}}%
    \else{\slshape\rightmark}\hfill{\caps\STitle}\hfill\llap{\thepage}\fi}%
\makeatother

\def\TSkip{\bigskip}
\newbox\TheTitle{\obeylines\gdef\GetTitle #1
\ShortTitle  #2
\SubTitle    #3
\Author      #4
\ShortAuthor #5
\EndTitle
{\setbox\TheTitle=\vbox{\baselineskip=20pt\let\par=\cr\obeylines%
\halign{\centerline{\Caps##}\cr\noalign{\medskip}\cr#1\cr}}%
	\copy\TheTitle\TSkip\TSkip%
\def\next{#2}\ifx\next\empty\gdef\STitle{#1}\else\gdef\STitle{#2}\fi%
\def\next{#3}\ifx\next\empty%
    \else\setbox\TheTitle=\vbox{\baselineskip=20pt\let\par=\cr\obeylines%
    \halign{\centerline{\caps##} #3\cr}}\copy\TheTitle\TSkip\TSkip\fi%
\centerline{\caps #4}\TSkip\TSkip%
\def\next{#5}\ifx\next\empty\gdef\SAuthor{#4}\else\gdef\SAuthor{#5}\fi%
\ifx\received\empty\relax
    \else\centerline{\eightrm Received: \received}\fi%
\ifx\revised\empty\TSkip%
    \else\centerline{\eightrm Revised: \revised}\TSkip\fi%
\ifx\communicated\empty\relax
    \else\centerline{\eightrm Communicated by \communicated}\fi\TSkip\TSkip%
\catcode'015=5}}\def\Title{\obeylines\GetTitle}
\def\Abstract{\begingroup\narrower
    \parskip=\medskipamount\parindent=0pt{\caps Abstract. }}
\def\EndAbstract{\par\endgroup\TSkip}

\long\def\MSC#1\EndMSC{\def\arg{#1}\ifx\arg\empty\relax\else
     {\par\narrower\noindent%
     2010 Mathematics Subject Classification: #1\par}\fi}

\long\def\KEY#1\EndKEY{\def\arg{#1}\ifx\arg\empty\relax\else
	{\par\narrower\noindent Keywords and Phrases: #1\par}\fi\TSkip}

\newbox\TheAdd\def\Addresses{\vfill\copy\TheAdd\vfill
    \ifodd\number\lastpage\vfill\eject\phantom{.}\vfill\eject\fi}
{\obeylines\gdef\GetAddress #1
\Address #2 
\Address #3
\Address #4
\EndAddress
{\def\xs{4.3truecm}\parindent=0pt
\setbox0=\vtop{{\obeylines\hsize=\xs#1\par}}\def\next{#2}
\ifx\next\empty 
     \setbox\TheAdd=\hbox to\hsize{\hfill\copy0\hfill}
\else\setbox1=\vtop{{\obeylines\hsize=\xs#2\par}}\def\next{#3}
\ifx\next\empty 
     \setbox\TheAdd=\hbox to\hsize{\hfill\copy0\hfill\copy1\hfill}
\else\setbox2=\vtop{{\obeylines\hsize=\xs#3\par}}\def\next{#4}
\ifx\next\empty\ 
     \setbox\TheAdd=\vtop{\hbox to\hsize{\hfill\copy0\hfill\copy1\hfill}
                \vskip20pt\hbox to\hsize{\hfill\copy2\hfill}}
\else\setbox3=\vtop{{\obeylines\hsize=\xs#4\par}}
     \setbox\TheAdd=\vtop{\hbox to\hsize{\hfill\copy0\hfill\copy1\hfill}
	        \vskip20pt\hbox to\hsize{\hfill\copy2\hfill\copy3\hfill}}
\fi\fi\fi\catcode'015=5}}\gdef\Address{\obeylines\GetAddress}

\def\LOCAL{\jobname.files}

\makeatletter
\long\def\@makecaption#1#2{%
  \vskip\abovecaptionskip
  \sbox\@tempboxa{\scriptsize {\scshape #1:} #2}%
  \ifdim \wd\@tempboxa >\hsize
    \scriptsize {\scshape #1:} #2\par
  \else
    \global \@minipagefalse
    \hb@xt@\hsize{\hfil\box\@tempboxa\hfil}%
  \fi
  \vskip\belowcaptionskip}
\makeatother

\begin{document}

\Title
Cylindrical Wigner Measures
\ShortTitle 
\SubTitle
\Author
Marco Falconi
\ShortAuthor
M. Falconi
\EndTitle
\Abstract
In this paper we study the semiclassical behavior of quantum states acting on the C*-algebra of canonical
commutation relations, from a general perspective. The aim is to provide a unified and flexible approach
to the semiclassical analysis of bosonic systems. We also give a detailed overview of possible
applications of this approach to mathematical problems of both axiomatic relativistic quantum field
theories and nonrelativistic many body systems. If the theory has infinitely many degrees of freedom, the
set of Wigner measures, \emph{i.e.}\ the classical counterpart of the set of quantum states, coincides
with the set of all cylindrical measures acting on the algebraic dual of the space of test functions for
the field, and this reveals a very rich semiclassical structure compared to the finite-dimensional
case. We characterize the cylindrical Wigner measures and the \emph{a priori} properties they inherit from
the corresponding quantum states.
\EndAbstract
\MSC 
\emph{Primary:} 81Q20, 81S05. \emph{Secondary:} 46L99, 47L90.
\EndMSC
\KEY 
Infinite Dimensional Semiclassical Analysis, CCR algebra, Constructive Quantum Field Theory, Wigner measures.
\EndKEY
\Address 
Fachbereich Mathematik\\Universität Tübingen\\ Auf der Morgenstelle 10\\72076 Tübingen
\Address
\Address
\Address
\EndAddress

\bibliographystyle{abbrvnat}

\section{Introduction}
\label{sec:introduction-1}

The study of semiclassical and effective behaviors in quantum mechanical systems with many particles plays
an important role in mathematical, theoretical, and experimental physics. In particular, we may mention
the recent widespread mathematical interest in the rigorous derivation of bosonic and fermionic effective
theories from many-body non-relativistic Hamiltonians \citep[it would be too long to provide here an
extensive bibliography on the subject, the interested reader should refer, \emph{e.g.}, to the
reviews][and references thereof
contained]{golse2016lnamm
  ,benedikter2016sbmp
}. For relativistic and semi-relativistic (bosonic) systems, where particles may be created or destroyed,
there are fewer results \citep{eckmann1975lmp
  ,donald1981cmp
  ,ginibre2006ahp
  ,falconi2013jmp
  ,ammari2014jsp
  ,ammari2014hok
  ,leopold2016arxiv
  ,ammari2017sima
} since the situation is usually more involved (\emph{e.g.}, due to renormalization issues). Most of the
latter results, at least the most recent ones, make use of the semiclassical approach in the Fock
representation developed by \citet{ammari2008ahp
  ,ammari2009jmp
  ,ammari2011jmpa
  ,ammari2015asns
}, that is also suitable to study non-relativistic mean-field problems \citep[see][in addition to the ones
just mentioned]{liard2014jmp
  ,ammari2016cms
  ,liard2017jfa
  ,ammari2017arxiv
}. This approach overcomes the limitations of other techniques, related to the choice of (initial) quantum
many-body state of the system: it is in fact possible to study the effective behavior of a general class
of quantum states, in particular states that have no coherent semiclassical structure, or that are not
``close'' to one with such structure.

The aim of this paper is to take an even more general approach, in order to complete our mathematical
understanding of bosonic semiclassical analysis, and to collect our knowledge in a unified description,
valid for most situations of physical interest. To do so, we connect the algebraic formulation of quantum
theories to semiclassical analysis. This has the advantage of taking into account at the same time all
possible representations of the algebra of canonical commutation relations, thus allowing to study the
semiclassical states of any given bosonic theory with fixed physical parameters such as mass and spin. As
a consequence, we are able to apply our results, \emph{e.g.}, to relativistic axiomatic field theories,
for which inequivalent representations of the canonical commutation relations play a crucial role, and to
thermodynamic states in quantum statistical mechanics. This and other applications are described in detail
in~\cref{sec:phys-appl} \citep[see
also][]{falconi2017ccm
  ,correggi2017ahp
  ,correggi2017arxiv
}. This paper can thus be seen as a continuation of the program started by Ammari and Nier, where we
combine the representation-independent algebraic framework to the powerful tools of semiclassical analysis
in a way that yields new physical applications, opening the way to interesting future
directions. Mathematically, we also provide a bridge between noncommutative and commutative concepts,
highlighting the way some properties of objects in a commutative theory may descend directly from the
properties of the corresponding objects in the deformed (noncommutative) counterpart, when the deformation
parameter converges to zero.

Throughout the paper, we may append the existential quantifier to mathematical objects as a superscript,
if the range of quantification is unambiguous. That is, letting $m,n\in \mathds{N}$, instead of $(\exists x_1\in X_1)(\exists x_2\in
X_2)\dotsm (\exists x_n\in X_n)R(x_1,\dotsc,x_n;y_1,\dotsc,y_m)$ we write $R(^{\exists }x_1,\dotsc,\,\!^{\exists }x_n;y_1,\dotsc,y_m)$. In particular,
the existential quantification appears as a superscript more than once only if appended to different
symbols, and each time it means a distinct quantification; if an object is bounded by an existential
quantifier and appears more than once in the same formula, the quantifier is appended to it only once,
\emph{e.g.}, $^{\exists }f\circ g=g\circ f= \mathrm{id}$. Even if this notation is somewhat unconventional, we think
that here it helps to lighten the notation, thus improving readability. In
\cref{sec:introduction-1,sec:phys-appl} an extensive use of footnotes is made, to quickly provide
definitions of objects and properties that may be unfamiliar to the reader.

\subsection{An overview of Wigner measures}
\label{sec:conc-wign-meas}

Wigner semiclassical measures are a powerful tool in the effective or asymptotic description of quantum
systems \citep[see, \emph{e.g.},][and references thereof
contained]{colindeverdiere1985cmp
  ,helffer1987cmp
  ,gerard1991cpde
  ,lions1993rmi
  ,fermanian2002bsmf
  ,riviere2010dmj
  ,zworski2012gsm
}. The standard Wigner measures are Radon measures on the cotangent bundle $\mathrm{T}^{*}\Sigma$ (phase space)
of a \emph{finite dimensional} real vector space $\Sigma$. They are the abelian counterpart of quantum states
that are normal with respect to the irreducible representation of the corresponding Heisenberg group
$\mathbb{H}(\mathrm{T}^{*}\Sigma)$. The semiclassical measures are naturally introduced in relation to
pseudodifferential (Weyl) calculus. In fact, let $h\geq 0$ be the semiclassical parameter. Then for any
bounded family of vectors $(u_h)_{0<h\leq h_0}\subset L^2 (\Sigma)$, $h_0\in \mathds{R}^+$, there exists a subsequence
$(u_{h_j})_{j\in \mathds{N}}$, $h_j\to 0$, and a finite Radon measure $\mu$ on $\mathrm{T}^{*}\Sigma$ such that for any $a\in
C^{\infty}_0(\mathrm{T}^{*}\Sigma)$
\begin{equation}
  \label{eq:fd}
  \lim_{j\to \infty}\langle u_{h_j}  , \mathrm{Op}_{\frac{1}{2}}^{h_j}(a)u_{h_j} \rangle_2=\int_{\mathrm{T}^{*}\Sigma}^{}a(z)  \mathrm{d}\mu(z)\; .
\end{equation}
$\mathrm{Op}_{\frac{1}{2}}^h(a)$ denotes the $h$-dependent Weyl quantization of the symbol $a$. The
physical interpretation of \cref{eq:fd} is straightforward: as quantum states are non-commutative
probabilities on which it is possible to evaluate quantum observables, classical states are classical
probabilities on which it is possible to evaluate phase space functions (classical observables). In
addition, the quantum expectation of a quantized classical observable converges in the classical limit to
the classical expectation, with respect to the semiclassical measure, of the same observable.

To study the semiclassical behavior of many-body quantum systems, Heisenberg groups associated to
\emph{infinite dimensional} symplectic spaces play a crucial role, since they encode the group form of the
canonical commutation relations. Pseudodifferential calculus for infinite dimensional symplectic phase
spaces $(X,\varsigma)$ is notoriously difficult, due to the lack of a locally finite and
translation-invariant measure. In addition, the Weyl C*-algebra associated to $(X,\varsigma)$ has
uncountably many inequivalent irreducible representations. Nonetheless, as we will explain below, it is
possible to use the finite-dimensional pseudodifferential calculus (quantizing the so-called cylindrical
symbols) to obtain a semiclassical characterization of states in quantum field theory. In addition, a
formula analogous to \cref{eq:fd} holds for all smooth cylindrical symbols, and it can be adapted to
suitable non-cylindrical and possibly non-smooth symbols as well \citep[see][for additional
details]{falconi2017ccm
}. There have been attempts to construct directly a pseudodifferential calculus for symbols that are not
cylindrical, if $(X,\varsigma)$ originates from a complex separable Hilbert space and the Weyl C*-algebra
is represented on the Fock vacuum (Fock representation). Among these attempts, let us mention the Weyl
calculus in abstract Wiener spaces \citep{amour2015jfa
  ,amour2016amrex
  ,amour2017jmp
}; the Wick quantization of polynomial symbols
\citep{berezin1971ms
  ,ammari2008ahp
}, and the inductive approach adopted by
\citet{kree1978aihpa
}.

In the Fock representation, the introduction of semiclassical Wigner measures as mean field or
semiclassical counterparts of bosonic quantum many-body states is due to
\citet{ammari2008ahp
}. Analogous measures also appear in the formulation of the bosonic quantum de Finetti Theorem
\citep{lewin2015amrx
} (see \citep[Proposition 3.2]{ammari2016cms
} for a link between the two points of view). We further develop the notion of Wigner measures introducing
\emph{cylindrical Wigner measures}. The additional adjective ``cylindrical'' is due to the fact that in
general the classical counterpart of bosonic quantum states are not Radon measures on the space of
classical fields, but rather cylindrical measures, \emph{i.e.}\ finitely additive measures on the algebra
of cylinders. The linear space of cylindrical measures includes the space of Radon measures, and the
inclusion is often strict for infinite dimensional spaces. We also make precise in which sense quantum
states converge to cylindrical measures in the semiclassical limit $h\to 0$, introducing two suitable
topologies. More precisely, even if \emph{a priori} quantum states and cylindrical measures may appear as
different types of objects, they are both isomorphic to objects of the same type, and a notion of vicinity
can therefore be introduced. We dedicate the rest of \cref{sec:introduction-1} to the motivation of the
basic ideas behind cylindrical Wigner measures, and the introduction of our main results. To improve
readability, the results are stated in a simplified form, more complete statements and proofs can be found
in \crefrange{sec:class-char-stat}{sec:group-homom-push}; \cref{sec:phys-appl} is devoted to physical
applications of the main results, and it is divided in two main subsections: in
\cref{sec:nonr-quant-field} we describe applications to many-body quantum theories and non-relativistic
quantum field theories, in \cref{sec:relat-quant-field} we describe applications to relativistic quantum
field theories.

\subsection{Semiclassical quantum states}
\label{sec:regul-quant-stat}

It is well known that, from an algebraic point of view, a quantum state is a positive and norm one
functional on the algebra of observables. We relax the normalization assumption, so throughout this paper
a \emph{quantum state} is a positive functional on the C*-algebra of quantum observables. We also adopt
the terminology \emph{normalized quantum state} and \emph{complex quantum state} to denote a normalized
positive functional and a complex functional respectively.

The algebra of bosonic quantum observables $\mathfrak{B}_h$ should embed the $h$-dependent Weyl C*-algebra
$\mathbb{W}_h(X,\varsigma)$ for some real symplectic space\footnote{Throughout the paper, we write
  categories in boldface. Let $\mathbf{C}$ be a category, by a slight abuse of notation we write $C\in
  \mathbf{C}$ for an object $C$ of the category (in some cases, to avoid confusion, we explicitly write
  $C\in \mathrm{Obj}(\mathbf{C})$). Morphism are denoted by $\mathfrak{c}\in \mathrm{Morph}(\mathbf{C})$.}
$(X,\varsigma)\in \mathbf{Symp}_{\mathds{R}}$, and should thus depend on the semiclassical parameter
$h$. The Weyl C*-algebra $\mathbb{W}_h(X,\varsigma)$ is defined as the smallest C*-algebra containing the
Weyl operators $\Bigl\{W_h(x),x\in X\Bigr\}$, \emph{i.e.}\ the elements indexed by $X$ satisfying the
following three properties:
\begin{align}
  \label{eq:i}\tag{i}  &W_h(x)\neq 0 &(\forall x\in X)\\
  \label{eq:ii}\tag{ii}  &W_h(-x)= W_h(x)^{*} &(\forall x\in X)\\
  \label{eq:iii}\tag{iii}  &W_h(x)W_h(y)=e^{-ih \varsigma(x,y)}W_h(x+y) &(\forall x,y\in X)
\end{align}
Let us emphasize that the dependence on the semiclassical parameter is given by the phase $e^{-i h
  \varsigma(x,y)}$ in \cref{eq:iii}, and that it physically corresponds to the fact that the commutator
between the canonical observables is of order $h$. Hence our analysis will apply to any such situation,
and $h$ may be in turn interpreted as a scale parameter, as an analogous of Planck's constant, or as a
quantity proportional to the inverse of the number of particles; for convenience, however, we always refer
to the regime $h\to 0$ as the semiclassical regime.  We are interested in characterizing the semiclassical
behavior of families of quantum states\footnote{The range of the semiclassical parameter is chosen to be
  $(0,1)$, however $1$ could be replaced with any strictly positive number. In addition, since states are
  elements of the continuous dual of the Banach space $\mathfrak{B}_h$, we can without loss of generality suppose that
  $\mathfrak{B}_h$ is *-isomorphic to $\mathbb{W}_h(X,\varsigma)$.}  $(\omega_h)_{h\in (0,1)}\subset
(\mathfrak{B}_h)'_+$, with
$\mathbb{W}_h(X,\varsigma)\overset{\mathrm{w}_X}{\hookrightarrow}\mathfrak{B}_h$. In order to converge to
a cylindrical Wigner measure, a family of quantum states should at least satisfy two properties: being
uniformly bounded in norm with respect to $h$, and being a family of \emph{regular} quantum states. A state
$\omega_h\in (\mathfrak{B}_h)'_+$ is regular iff for any $x\in X$, the $\mathds{R}$-action
\begin{equation*}
  \mathds{R}\ni\lambda\mapsto \omega_h\bigl(W_h(\lambda x)\bigr)\in \mathds{C}
\end{equation*}
is a \emph{continuous map}. Therefore, we make the following definition.
\begin{definition}[Semiclassical quantum states]
  \label{def:1}
  A (family of) quantum state(s) $(\omega_h)_{h\in I\subseteq (0,1)}\subset (\mathfrak{B}_h)'_+$ is
  \emph{semiclassical} iff zero is adherent to $I$, and
  \begin{align*}
    \mspace{50mu}&\bullet\quad\sup_{h\in I}\omega_h\bigl(W_h(0)\bigr) <\infty &,\\
    \mspace{50mu}&\bullet\quad\omega_{h} \textrm{ is }\mathrm{regular}&(\forall h\in I)\,.
  \end{align*}
\end{definition}
\begin{remark}
  \label{rem:4}
  All the definitions and results given in \crefrange{sec:regul-quant-stat}{sec:categ-interpr} could be
  reproduced, \emph{mutatis mutandis}, substituting Weyl C*-algebras with the corresponding
  \emph{Resolvent algebras} \citep{buchholz2008jfa
  }, and the Fourier transform of a cylindrical measure with the Stieltjes transform. The resolvent
  algebra may be more convenient than the Weyl C*-algebra to study dynamical interacting theories; however
  since there is a (purity preserving) bijection from the regular states of one algebra to the regular
  states of the other \citep{buchholz2008jfa
  }, for semiclassical purposes it is sufficient to focus on one of the two.
\end{remark}
From \cref{def:1}, it follows that a very large class of quantum states admits a semiclassical description
(\emph{e.g.}, any family of \emph{normalized} regular quantum states is semiclassical by \cref{def:1}). In
fact, the starting point of our analysis is to prove the following result. Let us denote by $X^{*}_X$ the
algebraic dual $X^{*}$ of $X$ endowed with the weak $\sigma(X^{*},X)$ topology.
\begin{thm}
  \label{thm:1}
  There exists a topology $\mathfrak{P}$ on the (disjoint) union of the set of all regular quantum states
  on $\mathfrak{B}_h \,\textrm{, } h\in (0,1)$, and of the set of all cylindrical measures on $X^{*}_X$ such
  that any semiclassical quantum state $(\omega_h)_{h\in I}$ is relatively compact. In addition, every
  $\mathfrak{P}$-cluster point of $\omega_h$ as $h\to 0$ is a cylindrical measure on $X^{*}_X$, called a
  \emph{cylindrical Wigner measure} of $\omega_h$.
\end{thm}
\begin{remark}
  \label{rem:2}
  A semiclassical quantum state $\omega_h$ has at least one cluster point, but it can have (infinitely)
  many cluster points, depending on which filter (or generalized sequence) converging to zero is taken as
  $I$. Although it is possible to construct semiclassical quantum states with different cluster points,
  such examples seem rather artificial and we could not think of a physically relevant example. One can
  therefore make the simplifying assumption that $I$ is chosen to be a filter on $(0,1)$ converging to
  zero, and on which $\omega_h$ converges.
\end{remark}
\begin{remark}
  \label{rem:1}
  A family of normalized quantum states $(\omega_h)_{h\in I}$ is not semiclassical if the states are not
  regular on some filter $F\subset I$ (equivalently, on some generalized sequence) converging to 0. In
  this case it is not, in general, possible to accurately approximate such states with a cylindrical
  probability (for small $h$). This is not completely unexpected, since non-regular states are associated
  with inherently quantum physical phenomena \citep[see,
  \emph{e.g.},][]{acerbi1993jmp
  }.
\end{remark}

\subsection{Cylindrical measures}
\label{sec:cylindrical-measures}

In order to understand \cref{thm:1} and the following results, it is important to review some properties
of cylindrical measures, a concept the reader may not be completely familiar with. Let $A\in \mathbf{Set}$
be a set, and let $\mathcal{A}\subseteq \mathds{R}^A$ be a subset of all the real-valued functions on
$A$. Let us denote by $\hat{C}(A,\mathcal{A})$ the \emph{initial $\sigma$-algebra} on $A$ with respect to
$\mathcal{A}$, \emph{i.e.}, the smallest $\sigma$-algebra of $A$ that makes all functions in $\mathcal{A}$
measurable. Equivalently, if $\mathcal{A}$ is a vector space then $\hat{C}(A,\mathcal{A})$ is the
$\sigma$-algebra generated by the algebra $C(A,\mathcal{A})$ of cylinders, where a cylinder is defined as
\begin{equation*}
  C_{\alpha_1,\dotsc,\alpha_n;B_1,\dotsc,B_n}=\{a\in A, \forall 1\leq j\leq n\;,\; \alpha_j(a) \in B_j\}\; ,
\end{equation*}
where $n\in \mathds{N}_{*}$, $\alpha_j\in \mathcal{A}$ and $B_j\in \mathrm{Borel}(\mathds{R})$ for any
$1\leq j\leq n$. Let us recall that if $\mathcal{A}$ is a vector space with a basis $\Delta$ of finite
cardinality (finite-dimensional space), then $C(A,\mathcal{A})=\hat{C}(A,\mathcal{A})=\hat{C}(A,\Delta)$.

\begin{definition}[Cylindrical measure]
  \label{def:2}
  Let $A\in \mathbf{Set}$, $\mathcal{A}\subseteq \mathds{R}^A$ a vector space. A finitely additive measure
  $M$ on $C(A,\mathcal{A})$ is a \emph{cylindrical measure} iff it is finite and its restriction to every
  $\sigma$-algebra $C(A,\mathcal{F})=\hat{C}(A,\mathcal{F})$, with $\mathcal{F}\subset \mathcal{A}$ a
  finite dimensional subspace, is countably additive. Let us denote by
  $\mathcal{M}_{\mathrm{cyl}}(A,\mathcal{A})$ the set of all cylindrical measures.
\end{definition}
It is clear that any (Radon) finite measure $\mu$ on the measurable space
$\bigl(A,\hat{C}(A,\mathcal{A})\bigr)$ is a cylindrical measure, but the converse is not true in general
(since there may be cylindrical measures that fail to be countably additive on
$\hat{C}(A,\mathcal{A})$). If $A\in \mathbf{TVS}_{\mathbb{R}}$ is a real topological vector space and
$\mathcal{A}=A'$ is its continuous dual, then the cylindrical measures can be equivalently defined as
projective families of Borel measures $(\mu_{\Phi},\pi_{\Phi,\Psi})_{\Phi\supset \Psi\in F(A)}$, where
$F(A)$ is the set of $\sigma(A,\mathcal{A})$-closed subspaces of $A$ with finite codimension, $\mu_{\Phi}$
is a Borel measure on $A/\Phi$, and $\pi_{\Phi,\Psi}:A/\Phi\leftarrow A/\Psi$ is the map obtained from the
identity map passing to the quotients, that acts on measures in a projective way pushing them forward:
$\mu_{\Phi}=\pi_{\Phi\Psi}\,_{*}\,\mu_{\Psi}$. This is due to the fact that the polar of a finite
dimensional subspace of $A'$ is isomorphic to a weakly closed subspace of $A$ of finite codimension, and
all weakly closed finite codimensional subspaces are polars of some finite dimensional subspace of $A'$.

For the purpose of semiclassical analysis, cylindrical measures are identified by the vector space
$\mathcal{A}$, rather than by the space $A$ on which they act upon, and there is a good extent of
arbitrariness in choosing the latter. In other words, there are (infinitely) many classical effective
theories corresponding to a given algebra of bosonic observables
$\Bigl(\mathbb{W}_h(X,\varsigma)\overset{\mathrm{w}_X}{\hookrightarrow}\mathfrak{B}_h\Bigr)_{h\in
  (0,1)}$. The arbitrariness is given by the fact that $A$ can be any set such that $X$ is linearly
embedded in $\mathds{R}^A$. In fact, if we denote by $e_X:X\overset{1\textrm{-}1}{\longrightarrow}
\mathds{R}^A$ the linear embedding, then we can actually interpret the cluster points of \cref{thm:1} as
cylindrical measures on $A$, with respect to the $\sigma$-algebra $\hat{C}\bigl(A,e_X(X)\bigr)$. However,
all such spaces of cylindrical measures are \emph{isomorphic to one another}, and in turn isomorphic to
the space of cylindrical measures on the topological vector space\footnote{Let us remark that the
  continuous dual of $X^{*}_X$ is isomorphic to $X$, thus the set of cylindrical measures on $X^{*}_X$ is
  exactly $\mathcal{M}_{\mathrm{cyl}}(X^{*},X)$ (where $X$ is identified with $(X^{*}_X)'$).}
$X^{*}_X$. This yields the following result.
\begin{proposition}
  \label{prop:1}
  The semiclassical description of bosonic quantum systems, in terms of classical cylindrical
  probabilities and observables, always exists and it is unique up to isomorphisms.
\end{proposition}
The isomorphisms among classical theories are given by Bochner's theorem\footnote{Applied to our context,
  Bochner's theorem is reformulated as the fact that the set of functions of positive type on the phase
  space $X$, continuous when restricted to any finite dimensional subspace of $X$, is isomorphic to the
  set of cylindrical measures on $A$ with respect to $e_X(X)$ (and such set of functions clearly does not
  depend on $A$). This is why all the semiclassical descriptions are isomorphic, and why only $X$ plays a
  role in their identification.}  \citep[see,
\emph{e.g.},][]{bourbaki1969int9
  ,schwartz1973tirsm
  ,vakhania1987ma
}.
\begin{thm}[Bochner]
  \label{thm:2}
  The Fourier transform is a bijection of the set of cylindrical measures on A, with respect to
  $\mathcal{A}$, onto the set of functions $g$ from $\mathcal{A}$ to $\mathds{C}$ such that:
  \begin{align*}
    \mspace{10mu}&\bullet\quad \sum_{j,k\in \mathfrak{F}}^{}\bar{\zeta}_k\zeta_j g(\alpha_j-\alpha_k)\geq 0 &\Bigl(\forall \mathfrak{F} \textrm{ finite set, }\forall (\zeta_j)_{j\in \mathfrak{F}}\subset \mathds{C} \text{ , } \forall (\alpha_j)_{j\in \mathfrak{F}}\subset \mathcal{A} \Bigr)\\
    \mspace{10mu}&\bullet\quad g\bigr\rvert_{\mathcal{F}} \textrm{ continuous } &\Bigl(\forall \mathcal{F}\subset \mathcal{A} \textrm{ , } \mathcal{F} \textrm{ finite dimensional} \Bigr)
  \end{align*}
\end{thm}

\subsection{The topologies of semiclassical convergence}
\label{sec:topol-conv}

Now that we have introduced cylindrical Wigner measures, as cluster points of semiclassical quantum
states, it is important to investigate the topology $\mathfrak{P}$ of convergence. As we will see, such
topology has a natural physical interpretation, but it is at times not convenient to work with. We
therefore introduce another topology for testing semiclassical convergence of quantum states, denoted by
$\mathfrak{T}$, that is more convenient for the explicit characterization of Wigner measures.

The topology $\mathfrak{P}$ is a (Hausdorff) weak topology, where the convergence is tested with smooth
cylindrical functions. In other words, a generalized sequence
\begin{equation*}
  \omega_{h_{\beta}}\underset{h_{\beta}\to 0}{\overset{\mathfrak{P}}{\longrightarrow}} M=(\mu_{\Phi})_{\Phi\in F(X^{*}_X)}
\end{equation*}
if and only if all the following expectations converge:
\begin{equation}
  \label{eq:1}
  \begin{split}
      \lim_{h_{\beta}\to 0}\omega_{h_{\beta}}\Bigl(\mathrm{Op}^{h_{\beta}}_{\frac{1}{2}}(f) \Bigr)\mspace{-3mu}:= \mspace{-3mu}\lim_{h_{\beta}\to 0}\omega_{h_{\beta}}\Bigl(\int_{\Phi^{\circ}}^{}\hat{f}_{\Phi}(x)W_{h_{\beta}}(\pi x)  \mathrm{d}x \Bigr) \mspace{-3mu}\\=\mspace{-3mu} \int_{X^{*}_X/\Phi}^{}f_{\Phi}(\xi)  \mathrm{d}\mu_{\Phi}(\xi)\mspace{-3mu}=:\mspace{-3mu}\int_{X^{*}_X}^{}f(z)  \mathrm{d}M(z);
  \end{split}
\end{equation}
where $\Phi\in F(X^{*}_X)$, $\Phi^{\circ}\subset X$ is its finite dimensional polar\footnote{The polar of
  a closed subspace of $X^{*}_X$ with finite codimension is a finite-dimensional subspace of $X$ that is
  isomorphic to the dual of $X^{*}_X/\Phi$ \citep[see,
  \emph{e.g.},][]{bourbaki1981evt1-5
  }.}, $f:X^{*}_X\to \mathds{C}$ is a smooth cylindrical function\footnote{A function $f:X^{*}_X\to
  \mathbb{C}$ is a \emph{smooth cylindrical function} iff there exists a $\Phi\in F(X^{*}_X)$ (the
  cylinder base) and a smooth function $f_{\Phi}\in C_0^{\infty}(X^{*}_X/\Phi)$ such that for any $z\in
  X^{*}_X$, $f(z)=f_{\Phi}(\pi_{\Phi}z)$, where $\pi_{\Phi}:X^{*}_X\to X^{*}_X/\Phi$ is the canonical
  projection. In other words, the function $f$ is determined only by finitely many degrees of freedom.}
with base function $f_{\Phi}\in C_0^{\infty}(X^{*}_X/\Phi)$, and
$\mathrm{Op}^{h_{\beta}}_{\frac{1}{2}}(f)$ its Weyl quantization.

From the physical standpoint, the average of a cylindrical function with respect to a cylindrical Wigner
measure is the number that approximates semiclassically the expectation on the quantum state of the
quantization of the aforementioned function. This is a good starting point, in order to have a useful
effective theory with some predictive power. It is, however, not easy to characterize the cylindrical
Wigner measure explicitly using $\mathfrak{P}$-convergence. Ideally, one would like to test with the
canonical quantum observables to characterize the measure. A convenient way to do that is using Weyl
operators. The generating functional of a quantum state $\omega_h\in (\mathfrak{B}_h)'_+$ is defined as
the numerical function characterizing the expectation of Weyl operators on the state:
\begin{equation}
  \label{eq:2}
  X \ni x\mapsto \mathcal{G}_{\omega_h}(x)=\omega_h\bigl(W_h(x)\bigr)\in \mathds{C}\; .
\end{equation}
The map $\omega_h\mapsto \mathcal{G}_{\omega_h}$ is the \emph{noncommutative Fourier transform}. The name
is due to the fact that for \emph{regular states} a noncommutative version of Bochner's theorem holds
\citep{segal1961cjm
}.
\begin{thm}[Noncommutative Bochner]
  \label{thm:3}
  The noncommutative Fourier transform is a bijection of the set of regular quantum states on
  $\mathbb{W}_h(X,\varsigma)$ onto the set of functions $\mathcal{G}_h$ from $X$ to $\mathbb{C}$ such
  that:
  \begin{align*}
    &\bullet\quad \sum_{j,k\in \mathfrak{F}}^{}\bar{\zeta}_k\zeta_j \mathcal{G}_h(x_j-x_k)e^{ih\varsigma(x_j\,,\,x_k)}\geq 0 &\Bigl(\forall \mathfrak{F} \textrm{ finite set, }\forall (\zeta_j)_{j\in \mathfrak{F}}\subset \mathds{C} \textrm{ , }\\
    & & \forall (x_j)_{j\in \mathfrak{F}}\subset X \Bigr)\\
    &\bullet\quad \mathcal{G}_h\bigr\rvert_{Y} \textrm{ continuous } &\Bigl(\forall Y\subset X \textrm{ , } Y \textrm{ finite dimensional } \Bigr)
  \end{align*}
\end{thm}
Noncommutative and commutative Fourier transforms provide therefore another way to treat quantum states
and cylindrical Wigner measures on the same grounds, as complex-valued functions on $X$ that are
continuous on any finite-dimensional subspace. The (Hausdorff) topology $\mathfrak{T}$ is thus defined to
be the preimage of the topology of pointwise convergence on the complex-valued functions on $X$. In other
words,
\begin{equation}
  \label{eq:3}
  \omega_{h_{\beta}}\underset{h_{\beta}\to 0}{\overset{\mathfrak{T}}{\longrightarrow}} M=(\mu_{\Phi})_{\Phi\in F(X^{*}_X)}
\end{equation}
if and only if the corresponding generating functional $\mathcal{G}_{\omega_{h_{\beta}}}$ converges
pointwise to $\hat{M}$ as $h_{\beta}\to 0$. This convergence provides a useful way to characterize the
Wigner measure, since the generating functional can be often explicitly computed, and thus its limit as
well, and the Fourier transform characterizes a cylindrical measure uniquely.

The topologies $\mathfrak{P}$ and $\mathfrak{T}$ are not comparable, so it may happen that a semiclassical
quantum state converges to a cylindrical measure in $\mathfrak{T}$, but not in $\mathfrak{P}$, or that it
converges to two different measures in the two topologies. However, the topology $\mathfrak{P}$ is the
physically relevant one, since it guarantees the convergence of the expectation of quantum (cylindrical)
observables. Therefore, we are only interested in semiclassical quantum states that converge either only
in $\mathfrak{P}$, or in both $\mathfrak{P}$ and $\mathfrak{T}$ \emph{to the same limit}. We therefore
make use of the join topology $\mathfrak{P}\vee \mathfrak{T}$, \emph{i.e.}, the coarsest topology that is
finer than both $\mathfrak{P}$ and $\mathfrak{T}$. Not all semiclassical quantum states have limit in the
$\mathfrak{P}\vee\mathfrak{T}$ topology, but statements can be formulated that are equivalent to such
convergence; they are presented in the theorem below. Let us denote by $\psi_{\varepsilon,\Phi}\in
C_0^{\infty}(X^{*}_X/\Phi)$ an approximation of the identity function, converging to it for $\varepsilon\to 0$.
\begin{thm}
  \label{thm:4}
  Let $\omega_h$ be a semiclassical quantum state on $\mathfrak{B}_h$ such that
  \begin{equation*}
    \omega_h\underset{h\to 0}{\overset{\mathfrak{P}}{\longrightarrow}} M=(\mu_{\Phi})_{\Phi\in F(X^{*}_X)}\;.
  \end{equation*}
  Then the following statements are equivalent:
  \begin{align}
    \label{eq:4}\tag{i}
    &\omega_h\underset{h\to 0}{\overset{\mathfrak{P}\vee\mathfrak{T}}{\longrightarrow}} M&; \\
    \label{eq:5}\tag{ii}
    &\lim_{\varepsilon\to 0}\lim_{h\to 0}\omega_h\bigl(W_h(x)-\mathrm{Op}_{\frac{1}{2}}^h(^{\exists }\psi_{\varepsilon,\Phi}(\cdot )\, e^{2ix(\cdot )})\bigr)=0&\bigl(\forall \Phi\in F(X^{*}_X)\,,\,\forall x\in \Phi^{\circ}\bigr)\,;\\
    \label{eq:6}\tag{iii}
    &\lim_{\varepsilon\to 0}\lim_{h\to 0}\omega_h\bigl(1-\mathrm{Op}_{\frac{1}{2}}^h(^{\exists }\psi_{\varepsilon,\Phi})\bigr)=0&\bigl(\forall \Phi\in F(X^{*}_X)\bigr)\,;\\
    \label{eq:7}\tag{iv}
    &\lim_{h\to 0}\omega_h\bigl(W_h(0)\bigr)=M(X^{*}_X)&.
  \end{align}
\end{thm}
We call a $\mathfrak{P}$-convergent semiclassical quantum state that satisfies \crefrange{eq:4}{eq:7} a
\emph{semiclassical quantum state with no loss of mass}.

\subsection{The set of all cylindrical Wigner measures}
\label{sec:set-all-cylindrical}

It is important to characterize the set of all possible semiclassical configurations of a given system. By
\cref{thm:1,prop:1}, we know that all semiclassical configurations can be equivalently described as
cylindrical measures acting on measurable or topological vector spaces. With cylindrical measures,
however, it is only possible to integrate cylindrical functions. On the other hand, cylindrical measures
can also be seen as Borel Radon measures, but usually on a ``big'' topological space. If we take the
prototypical example of phase space in nonrelativistic quantum field theory,
$\bigl(\mathbb{F}\!\mathbb{F}L^2(\mathds{R}^d),\im \langle \,\cdot \, , \,\cdot \,
\rangle_2\bigr)$\footnote{$\mathbb{F}\!\mathbb{F}$ denotes the forgetful functor from complex to real
  vector spaces. In other words, $\mathbb{F}\!\mathbb{F}L^2$ is the vector space that has the same
  elements as $L^2$, but only multiplication by real scalars is allowed (therefore its basis is
  ``doubled''); $\mathbb{F}\!\mathbb{F}L^2$ is a real separable Hilbert space with scalar product $\re
  \langle \,\cdot \, , \,\cdot \, \rangle_2$, and a symplectic space with symplectic form $\im \langle
  \,\cdot \, , \,\cdot \, \rangle_2$.}, the cylindrical measures associated to it are Radon measures only
when acting on the weak Hausdorff completion of $L^2(\mathds{R}^d)$ (\emph{i.e.}, on the completion of
$L^2$ endowed with the $\sigma(L^2,L^2)$ topology). In addition, it is well known that there are true
cylindrical measures on $L^2$, such as the Gaussian measure, that are concentrated outside of $L^2$ when
considered as Radon measures \citep[see, \emph{e.g.}, the construction of abstract Wiener
spaces,][]{gross1967pfbsmsp
}. One may therefore hope that the set of cylindrical Wigner measures is contained strictly in the set of
cylindrical measures, and that it coincides with the space of Radon measures on some topological vector
space with nice properties. However, \emph{this is not the case}: every cylindrical measure is the Wigner
measure of at least one semiclassical quantum state.
\begin{thm}
  \label{thm:5}
  The set of all cluster points, in the $\mathfrak{P} \vee \mathfrak{T}$ topology, of semiclassical
  quantum states on $(\mathfrak{B}_h)_{h\in (0,1)}$ coincides with the cone
  $\mathcal{M}_{\mathrm{cyl}}(X^{*},X)$ of all cylindrical measures on $X^{*}_X$ (or on any other set $A$
  for which $X\overset{e_X}{\hookrightarrow}\mathds{R}^A$).
\end{thm}
This result is, in our opinion, quite interesting. It shows that cylindrical measures are indeed the
states that emerge from a bosonic field theory in its classical approximation. This perhaps suggests that
a classical field theory well-suited for quantization should either take into account such cylindrical
structure, or that it should be set in a topological space where all cylindrical Wigner measures are Radon
measures.

\subsection{Maps, convolutions, products on cylindrical Wigner measures}
\label{sec:maps-conv-prod}

Now that we know what is the classical counterpart of a regular (bosonic) state, we would like to
characterize which structures and properties of quantum states are inherited by the corresponding
classical states. First of all, let us consider the action of central group-homomorphisms between
Heisenberg groups, acting on the corresponding Weyl C*-algebras as *-homomorphisms mapping Weyl operators
into Weyl operators. A map of such type, from the Heisenberg group $\mathbb{H}(X,\varsigma)$ to
$\mathbb{H}(Y,\tau)$, is induced by a \emph{symplectic map} $s:X\to Y$,
$\tau\bigl(s(x),s(x')\bigr)=\varsigma(x,x') \; (\forall x,x'\in X)$. In turn, $s$ induces the desired
*-homomorphism between Weyl C*-algebras preserving Weyl operators:
\begin{equation}
  \label{eq:9}
  \mathfrak{s}_h : \underset{\mspace{19mu}W_h(x)\quad\longmapsto \quad W_h(s(x))}{\mathbb{W}_h(X,\varsigma)\longrightarrow\mathbb{W}_h(Y,\tau)}\; ;
\end{equation}
that can be extended to a *-homomorphism from $\mathfrak{B}_h$ to $\mathfrak{C}_h$, whenever
$\mathbb{W}_h(X,\varsigma)\overset{\mathrm{w}_X}{\hookrightarrow}\mathfrak{B}_h$ and
$\mathbb{W}_h(Y,\tau)\overset{\mathrm{w}_Y}{\hookrightarrow}\mathfrak{C}_h$. As we discuss in
\cref{sec:categ-interpr}, these are the natural morphisms among quantum bosonic theories, induced by the
(Segal) quantization functor. The map $\mathfrak{s}_h$ of quantum observables induces, by duality, a
\emph{positivity preserving continuous linear} map on quantum states, the transposed map:
\begin{equation*}
  ^{\mathrm{t}}\mathfrak{s}_h: (\mathfrak{C}_{h})'_+\longrightarrow (\mathfrak{B}_{h})'_+\; .
\end{equation*}

On the other hand, at the classical level, the map $s:X\to Y$ induces by duality a \emph{continuous
  linear} map between $Y^{*}_Y$ and $X^{*}_X$:
\begin{equation*}
  ^{\mathrm{t}}s:Y^{*}_Y \longrightarrow X^{*}_X\; ,
\end{equation*}
that acts on cylindrical measures pushing them forward\footnote{The push-forward of a cylindrical measure
  $M=(\mu_{\Phi})_{\Phi\in F(V)}$ by a weakly continuous linear function $u:V\to W$ on a topological
  vector space $V$ is defined as follows. For any $\Xi\in F(W)$, let $u^{-1}(\Xi)\in F(V)$; and let
  $u_{\Xi}:V/u^{-1}(\Xi)\to W/\Xi$ be the linear application obtained from $u$ passing to the
  quotients. Then $u\, _{*}\, M=(u_{\Xi}\, _{*}\, \mu_{u^{-1}(\Xi)})_{\Xi\in F(W)}$.}: given a cylindrical
measure $M_Y$ on $Y^{*}_Y$, then $M_X:=\,\!^{\mathrm{t}}s \,_{*}\, M_Y$ is a cylindrical measure on
$X^{*}_X$. It is therefore natural to ask whether $^{\mathrm{t}}\mathfrak{s}_h$ converges in some sense to
$^{\mathrm{t}}s$. The answer is positive when the quantum map acts on semiclassical quantum states
converging in the $\mathfrak{P}\vee\mathfrak{T}$ topology:
\begin{equation}
  \label{eq:8}
  \begin{aligned}
      &\varpi_{h_{\beta}}\underset{h_{\beta}\to 0}{\overset{\mathfrak{P}\vee\mathfrak{T}}{\longrightarrow}}M \quad \Longrightarrow \quad ^{\mathrm{t}}\mathfrak{s}_{h_{\beta}}\varpi_{h_{\beta}}\underset{h_{\beta}\to 0}{\overset{\mathfrak{P}\vee\mathfrak{T}}{\longrightarrow}} \,\!^{\mathrm{t}}s \,_{*}\, M &\bigl(\forall \varpi_{h_{\beta}}\in (\mathfrak{C}_{h_{\beta}})'_+ \textrm{ semiclassical} \\ & &\textrm{quantum state} \bigr).
  \end{aligned}
\end{equation}
As it will be outlined in \cref{sec:relat-quant-field}, such result is important to study the
semiclassical behavior of symmetry transformations of relativistic bosonic fields.

The product of (cylindrical) measures is related to statistical independence of classical systems. A
system of classical bosonic fields (or particles) composed of two subsystems can be modeled by a
topological vector space $W=V_1\times V_2$ that is product of two spaces $V_1,V_2\in
\mathbf{TVS}_{\mathds{R}}$, each one describing a subsystem. The state of the system is given in
statistical mechanics by a (cylindrical) measure $M$ on $W$. If such measure can be written as the
(tensor) product\footnote{Let $V_1$ and $V_2$ be topological vector spaces. The (tensor) product
  $M_1\otimes M_2$ of two cylindrical measures $M_1$ on $V_1$ and $M_2$ on $V_2$ is the cylindrical
  measure on $V_1\times V_2$ defined as follows.  $M_1\otimes M_2=\bigl(\mu_{V_1/\Phi_1}\otimes
  \mu_{V_2/\Phi_2}\bigr)_{\Phi_1\in F(V_1),\Phi_2\in F(V_2)}$, where $\mu_{V_1/\Phi_1}\otimes
  \mu_{V_2/\Phi_2}$ is the usual product measure on the finite dimensional space $(V_1/\Phi_1)\times
  (V_2/\Phi_2)\cong V_1\times V_2/\Phi_1\times \Phi_2$. The family $\bigl(\mu_{V_1/\Phi_1}\otimes
  \mu_{V_2/\Phi_2}\bigr)_{\Phi_1\in F(V_1),\Phi_2\in F(V_2)}$ is projective, and the set of subspaces of
  the form $\Phi_1\times \Phi_2$ is cofinal with the set of all weakly closed subspaces of finite
  codimension. Therefore $M_1\otimes M_2$ defines a cylindrical measure on $V_1\times V_2$.}
$M=M_1\otimes M_2$ of a measure on $V_1$ and one on $V_2$, then the two subsystems are independent. This
is due to the fact that in such case an event in one of the two subsystems would not affect the outcome of
events in the other. At the quantum level, a bosonic composite system is described by the tensor product
of two C*-algebras $\mathfrak{B}_h\otimes_{\alpha} \mathfrak{C}_h$, where the index $\alpha$ stands for a
suitable choice of cross norm for the product algebra, and
$\mathbb{W}_h(X,\varsigma)\overset{\mathrm{w}_{X}}{\hookrightarrow}\mathfrak{B}_h$,
$\mathbb{W}_h(Y,\tau)\overset{\mathrm{w}_{Y}}{\hookrightarrow}\mathfrak{C}_h$ (from which it follows that
$\mathbb{W}_h(X \oplus Y,\varsigma \oplus \tau)\overset{^{\exists
  }\mathrm{w}_{XY}}{\hookrightarrow}\mathfrak{B}_h\otimes_{\alpha} \mathfrak{C}_h$). If a quantum state
$\xi_h$ on $\mathfrak{B}_h\otimes_{\alpha} \mathfrak{C}_h$ is of the form $\xi_h=\omega_h\otimes
\varpi_h$, with $\omega_h\in (\mathfrak{B}_h)'_+$ and $\varpi_h\in (\mathfrak{C}_h)'_+$, then the two
subsystems are statistically independent. The lack of statistical independence in quantum systems is often
called quantum entanglement. The following proposition is a well-known physical fact about entanglement,
that we can rigorously prove in our framework.
\begin{proposition}
  \label{prop:2}
  Entanglement can at most be destroyed by the classical limit.
\end{proposition}
The proof, whose details are given in \cref{sec:semicl-einst-caus}, consists of two parts. The first is to
construct an entangled semiclassical quantum state whose limit is a product cylindrical measure, and thus
statistically independent. This can be easily done, \emph{e.g.}, using an entangled ``perturbation'' of
order $h$ of a non-entangled state. The second part consists in proving that all possible cylindrical
Wigner measures of a semiclassical quantum state $\omega_h\otimes \varpi_h$ are of the form $M_X\otimes
M_Y$, with $M_X$ a cylindrical measure on $X^{*}_X$ and $M_Y$ a cylindrical measure on $Y^{*}_Y$.

Another useful notion in classical probability is that of convolution of (cylindrical)
measures\footnote{The convolution $M_1\circledast M_2$ of two cylindrical measures $M_1$ and $M_2$ on
  $V\in \mathbf{TVS}_{\mathds{R}}$, is again a cylindrical measure on $V$, defined as follows. It is the
  pushforward of the product $M_1\otimes M_2$ by the addition map $\tilde{+}:V\times V\to V$ defined by
  $(v,w)\mapsto v+w$. In other words, $M_1\circledast M_2= \tilde{+}\;_{*}\, M_1\otimes M_2$.}. Signed or
complex cylindrical measures are an algebra, if we use the convolution as product. The Fourier transform
on cylindrical measures behaves as the usual Fourier transform with respect to convolution (and is thus an
algebra homomorphism from the algebra of signed/complex cylindrical measures to the algebra of
complex-valued functions):
\begin{equation*}
  \widehat{\bigl(M_1\circledast M_2\bigr)}( \cdot )=\hat{M}_1( \cdot )\hat{M}_2(\cdot ) \; .
\end{equation*}
This suggests an analogous definition of quantum (noncommutative) convolution, that would make the complex
quantum states an algebra (and the noncommutative Fourier transform an algebra homomorphism). Let
$\omega_h,\varpi_h\in \mathbb{W}_h(X,\varsigma)'$, and define the action $\omega_h\star\varpi_h$ on Weyl
operators by
\begin{equation*}
  (\omega_h\star\varpi_h)\bigl(W_h(x)\bigr)=\omega_h\bigl(W_h(x)\bigr)\varpi_h\bigl(W_h(x)\bigr)\; .
\end{equation*}
This action is extended by linearity to finite combinations
\begin{equation*}
  (\omega_h\star\varpi_h)\Bigl(\sum_{j\in \mathfrak{F}}^{}\zeta_jW_h(x_j)\Bigr)=\sum_{j\in \mathfrak{F}}^{}\zeta_j(\omega_h\star\varpi_h)\bigl(W_h(x_j)\bigr)
\end{equation*}
and hence to a (complex) state on $\mathbb{W}_h(X,\varsigma)$. Semiclassically, we have the following
result:
\begin{align*}
  \omega_{h_{\beta}}\underset{h_{\beta}\to 0}{\overset{\mathfrak{P}\vee\mathfrak{T}}{\longrightarrow}}M_1\;,\; \varpi_{h_{\beta}}\underset{h_{\beta}\to 0}{\overset{\mathfrak{P}\vee\mathfrak{T}}{\longrightarrow}}M_2 \quad \Longrightarrow \quad \omega_{h_{\beta}}\star\varpi_{h_{\beta}}\underset{h_{\beta}\to 0}{\overset{\mathfrak{P}\vee\mathfrak{T}}{\longrightarrow}} M_1\circledast M_2\; .
\end{align*}

\subsection{Categorical interpretation}
\label{sec:categ-interpr}

The framework of cylindrical Wigner measures admits an elegant reformulation in the language of category
theory. Such reformulation is well-adapted to algebraic formulations of relativistic quantum field theory,
such as the Locally Covariant Quantum Field Theory, and thus we present it here. Some results of
semiclassical analysis in LCQFT are given in \cref{sec:phys-appl}.

Let $\mathbf{Symp}_{\mathds{R}}$ be the category of real symplectic spaces, \emph{i.e.}, the category with
objects the real symplectic vector spaces, and morphisms the linear symplectic maps; and let
$\mathbf{C^{*}alg}$ be the category of C*-algebras, \emph{i.e.}, the category with C*-algebras as objects
and *-homomorphisms as morphisms.
\begin{definition}[Segal bosonic quantization]
  \label{def:3}
  Let $h\in (0,1)$. The functor
  \begin{equation*}
    \mathbb{W}_h: \mathbf{Symp}_{\mathds{R}}\to \mathbf{C^{*}alg}
  \end{equation*}
  that associates:
  \begin{itemize}
  \item to each symplectic space $(X,\varsigma)$ the Weyl C*-algebra $\mathbb{W}_h(X,\varsigma)$,
  
  \item to each linear symplectomorphism $s\!:(X,\varsigma)\to (Y,\tau)$ the associated
    Weyl-operator-preserving *-homomorphism $\mathbb{W}_h(s):=\mathfrak{s}_h\,$, defined by~\cref{eq:9}~,
  \end{itemize}
  is called \emph{Segal bosonic quantization}.
\end{definition}
The Segal quantization associates to each classical phase space the corresponding algebra of observables
generated by the canonical commutation relations, and to linear maps between symplectic spaces the
corresponding *-homomorphism of observables, that maps Weyl operators into Weyl operators\footnote{The
  Segal functor should not be confused with the second quantization functor
  \citep{nelson1973jfa
  }, that could be considered as a special case of the former.}. Now, let $\mathbf{BanCone}$ be the
category with objects the pointed and generating cones in real Banach spaces, and with morphisms the
linear continuous cone preserving maps between underlying Banach spaces. Hence the contravariant duality
functor $\mathbb{D}_+: \mathbf{C^{*}alg}\to \mathbf{BanCone}$ associates to each C*-algebra the
corresponding cone of quantum states (positive linear functionals), and to *-homomorphisms their
transposed maps. It follows that
\begin{equation*}
  \mathbb{S}_h=\mathbb{D}_+\circ \mathbb{W}_h: \mathbf{Symp}_{\mathds{R}}\to \mathbf{BanCone}
\end{equation*}
is the (contravariant) \emph{functor of bosonic quantum states}, associating to each classical phase space
the corresponding quantum states in the algebra of canonical commutation relations.

The classical limit of the functor of bosonic quantum states is the \emph{functor of cylindrical Wigner
  measures}, defined as follows. Let $\mathbf{CylM}$ be the category that has as objects the sets of
cylindrical measures, and as morphisms the compatible maps\footnote{A map $f:A\to B$, $A,B\in
  \mathbf{Set}$ is compatible with the sets of cylindrical measures
  $\mathcal{M}_{\mathrm{cyl}}(A,\mathcal{A})$ and $\mathcal{M}_{\mathrm{cyl}}(B,\mathcal{B})$,
  $\mathcal{A}\subseteq \mathds{R}^{A},\mathcal{B}\subseteq \mathds{R}^B$, iff for any $\beta\in
  \mathcal{B}$, $\beta\circ f\in \mathcal{A}$. In that case, it is possible to define the pushforward $f\,
  _{*}\, M \in \mathcal{M}_{\mathrm{cyl}}(B,\mathcal{B})$ for any $M\in
  \mathcal{M}_{\mathrm{cyl}}(A,\mathcal{A})$, and therefore $f$ induces a homomorphism between sets of
  cylindrical measures (again denoted by $f$).}.
\begin{definition}[Functor of cylindrical Wigner measures]
  \label{def:4}
  A contravariant functor
  \begin{equation*}
    \mathbb{S}_0: \mathbf{Symp}_{\mathds{R}}\to \mathbf{CylM}
  \end{equation*}
  is a \emph{functor of cylindrical Wigner measures} iff \footnote{With a slight abuse of notation, we use
    the quantifier $\forall $ for categories that are not small. The corresponding statements should be
    intended in the appropriate sense.}:
  \begin{align}
    \label{eq:10}\tag{i}
    & \mathbb{S}_0(X,\varsigma)=\mathcal{M}_{\mathrm{cyl}}\bigl(\,\!^{\exists }A,\,\!^{\exists }e_X(X)\bigr)\;,\;X\overset{e_X}{\hookrightarrow}\mathds{R}^A \mspace{80mu}\Bigl(\forall (X,\varsigma)\in \mathbf{Symp}_{\mathds{R}}\Bigr)\\
    \label{eq:11}\tag{ii}
    &\begin{aligned}
      &e_Y\bigl(s(x)\bigr)=e_X(x)\circ \mathbb{S}_0(s) &\mspace{80mu}\Bigl(\forall x\in X\,,\;\forall s\in \mathrm{Morph}(\mathbf{Symp}_{\mathds{R}})\,,\,\\
      & &s:(X,\varsigma)\to (Y,\tau) \Bigr)
    \end{aligned}
  \end{align}
\end{definition}
\cref{def:4} takes into account the fact that it is possible to equivalently describe Wigner measures as
cylindrical measures in different spaces, as discussed in \cref{sec:cylindrical-measures}. However, the
semiclassical description always exists and it is essentially unique, as stated in \cref{prop:1}. In the
categorical setting, such statement takes the following form.
\begin{proposition}
  \label{prop:3}
  The functor of cylindrical Wigner measures always exists. If $\,\mathbb{S}_0$ and $\mathbb{T}_0$ are two
  functors of cylindrical Wigner measures, then there exists a natural isomorphism
  $\nu:\mathbb{S}_0\overset{\mathrm{iso}}{\longrightarrow} \mathbb{T}_0$.
\end{proposition}
\begin{proof}
  The existence is guaranteed by the fact that it is possible to choose, for all $(X,\varsigma)\in
  \mathbf{Symp}_{\mathds{R}}$, $A=X^{*}_X$, identifying $X$ with $(X^{*}_X)'$. This choice satisfies
  \eqref{eq:10}. The transformed symplectic morphisms are also easily defined as
  $\mathbb{S}_0(s)=\,^{\mathrm{t}}s$, satisfying \eqref{eq:11}. Uniqueness is proved as follows. Let
  $\mathbb{S}_0$ and $\mathbb{T}_0$ be two functors of cylindrical Wigner measures. Let us show that there
  is a natural isomorphism $\nu: \mathbb{S}_0\to \mathbb{T}_0$. Let $(X,\varsigma)\in
  \mathbf{Symp}_{\mathds{R}}$. Define the component $\nu_X: \mathbb{S}_0(X,\varsigma)\to
  \mathbb{T}_0(X,\varsigma)$ in the following way. By definition of classical functor,
  $\mathbb{S}_0(X,\varsigma)=\mathcal{M}_{\mathrm{cyl}}\bigl(A,e_X(X)\bigr)$ and
  $\mathbb{T}_0(X,\varsigma)=\mathcal{M}_{\mathrm{cyl}}\bigl(B,f_X(X)\bigr)$ with
  $X\overset{e_X}{\hookrightarrow}\mathds{R}^A$ and $X\overset{f_X}{\hookrightarrow}\mathds{R}^B$. By
  Bochner's theorem, \cref{thm:2}, the Fourier transform induces an isomorphism
  \begin{equation*}
    \nu_X: \mathcal{M}_{\mathrm{cyl}}\bigl(A,e_X(X)\bigr)\to \mathcal{M}_{\mathrm{cyl}}\bigl(B,f_X(X)\bigr)\; .
  \end{equation*}
  The commutativity of the diagram
  \begin{equation*}
    \begin{tikzcd}
      \mathbb{S}_0(X,\varsigma) \arrow[d,swap,"\nu_X"]& \mathbb{S}_0(Y,\tau) \arrow[l,swap,"\mathbb{S}_0(s)"]\arrow[d,swap,"\nu_Y"]\\
      \mathbb{T}_0(X,\varsigma) &\mathbb{T}_0(Y,\tau) \arrow[l,swap,"\mathbb{T}_0(s)"]
    \end{tikzcd}
  \end{equation*}
  follows since by \eqref{eq:11} of \cref{def:4}, $e_Y\bigl(s(x)\bigr)=e_X(x)\circ\mathbb{S}_0(s)$ and
  $\,f_Y\bigl(s(x)\bigr)=f_X(x)\circ \mathbb{T}_0(s)$.
\end{proof}

All the results of \crefrange{sec:regul-quant-stat}{sec:maps-conv-prod} can conveniently be reinterpreted
as a ``semiclassical convergence'' of $\mathbb{S}_h$ to $\mathbb{S}_0$, as $h\to 0$.
\begin{thde}
  \label{th:de}
  \begin{equation*}
    \mathbb{S}_h\underset{h\to 0}{\rightarrowtail} \mathbb{S}_0\; .
  \end{equation*}
\end{thde}
\begin{proof}
  The semiclassical convergence of functors $\underset{h\to 0}{\rightarrowtail}$ holds in the following
  sense in any \emph{small} subcategory $\mathbf{S}\subset \mathbf{Symp}_{\mathds{R}}$:
  \begin{itemize}
  \item \emph{Convergence within objects} (\cref{thm:1}). There exists a topology for the element-wise
    semiclassical convergence of objects:
    \begin{equation*}
      \mathfrak{P}_{\mathbf{S}}:= \prod_{(X,\varsigma)\in \mathbf{S}}^{}\mathfrak{P}_X\; ,
    \end{equation*}
    where $\mathfrak{P}_X$ is the topology of semiclassical convergence defined in
    \cref{sec:topol-conv}. In other words, any semiclassical family of maps
    $\Bigl(\mathbf{S}\ni(X,\varsigma)\mapsto \omega_{h,X}\in \mathbb{S}_h(X,\varsigma) \Bigr)_{h\in
      I\subseteq (0,1) }$, with zero adherent to $I$ and $\omega_{h,X}$ semiclassical, is relatively
    compact in the $\mathfrak{P}_{\mathbf{S}}$ topology. Therefore the element-wise convergence of
    semiclassical objects is given by
    \begin{equation*}
      \Bigl((X,\varsigma)\mapsto \omega_{h_{\beta},X} \Bigr)\underset{h_{\beta}\to 0}{\overset{\mathfrak{P}_{\mathbf{S}}}{\longrightarrow}}\Bigl((X,\varsigma)\mapsto M_X \Bigr)\; ,
    \end{equation*}
    with $M_X\in \mathbb{S}_0(X,\varsigma)\cong\mathcal{M}_{\mathrm{cyl}}(X^{*},X)$ for any
    $(X,\varsigma)\in \mathbf{S}$.
    
  \item \emph{Convergence of morphisms} (\cref{eq:8}). There exists a topology $\mathfrak{T}_{\mathbf{S}}$
    on objects for the pointwise semiclassical convergence of morphisms:
    \begin{equation*}
      \mathfrak{T}_{\mathbf{S}}:=\prod_{(X,\varsigma)\in \mathbf{S}}^{}\mathfrak{T}_X\;,
    \end{equation*}
    where $\mathfrak{T}_X$ is the other topology of semiclassical convergence defined in
    \cref{sec:topol-conv}. In fact, let $\Bigl(\mathrm{Morph}(\mathbf{S})\ni \bigl(s\!:\!(X,\varsigma)\to
    (Y,\tau)\bigr)\,\longmapsto \,\mathbb{S}_h(s)= \,\!^{\mathrm{t}}\mathfrak{s}_h\in
    \mathcal{L}\bigl(\mathbb{S}_h(Y,\tau),\mathbb{S}_h(X,\varsigma)\bigr) \Bigr)_{h\in (0,1)}$, be a
    family of maps of morphisms. Then on any semiclassical family of
    $(\mathfrak{P}_{\mathbf{S}}\vee\mathfrak{T}_{\mathbf{S}})$-convergent maps
    \begin{equation*}
      \Bigl((X,\varsigma)\mapsto \omega_{h_{\beta},X} \Bigr)\underset{h_{\beta}\to 0}{\overset{\mathfrak{P}_{\mathbf{S}}\vee\mathfrak{T}_{\mathbf{S}}}{\longrightarrow}}\Bigl((X,\varsigma)\mapsto M_X \Bigr)\; ,
    \end{equation*}
    we have the pointwise convergence $\Bigl(s\mapsto
    \mathbb{S}_{h_{\beta}}(s)\Bigr)\underset{h_{\beta}\to 0}{\longrightarrow}\Bigl(s\mapsto
    \mathbb{S}_0(s)\Bigr)$:
    \begin{equation*}
      \Bigl((Y,\tau)\mapsto \mathbb{S}_{h_{\beta}}(s) \omega_{h_{\beta},Y} \Bigr)\underset{h_{\beta}\to 0}{\overset{\mathfrak{P}_{\mathbf{S}}\vee\mathfrak{T}_{\mathbf{S}}}{\longrightarrow}}\Bigl((Y,\tau)\mapsto \mathbb{S}_0(s)\,_{*}\, M_Y \Bigr)\; .
    \end{equation*}
    
  \item $\mathrm{Ran}\bigl(\mathbb{S}_h\bigr\rvert_{\mathbf{S}}\bigr)\underset{h\to
      0}{\rightarrowtail}\mathrm{Ran}\bigl(\mathbb{S}_0\bigr\rvert_{\mathbf{S}}\bigr)$ (\cref{thm:5}). We
    define the range of a functor $\mathbb{F}:\mathbf{A}\to \mathbf{B}$, $\mathbf{A}$ small, as the
    set\footnote{We use here the following standard set-theoretic notation \citep[see,
      \emph{e.g.},][]{jech2003smm
      }: given a set $S$, we denote by $\bigcup S$ the set
      \begin{equation*}
        \bigcup S=\Bigl\{t\,,\, t\in \,\!^{\exists }s\in S\Bigr\}\; .
      \end{equation*}
    }
    \begin{equation*}
      \mathrm{Ran} (\mathbb{F})=\bigcup \Bigl\{\mathbb{F}(a)\,,\,a\in \mathbf{A}\Bigr\} \; .
    \end{equation*}
    The notation $\mathrm{Ran}\bigl(\mathbb{S}_h\bigr\rvert_{\mathbf{S}}\bigr)\underset{h\to
      0}{\rightarrowtail}\mathrm{Ran}\bigl(\mathbb{S}_0\bigr\rvert_{\mathbf{S}}\bigr)$ should be
    interpreted as the fact that every element $M_{^{\exists }X}\in
    \mathrm{Ran}\bigl(\mathbb{S}_0\bigr\rvert_{\mathbf{S}}\bigr)$ is the limit point (in the
    $\mathfrak{P}_X\vee\mathfrak{T}_X$ topology) of at least one semiclassical state $^{\exists
    }\omega_{h,M_X}\in \mathrm{Ran}\bigl(\mathbb{S}_h\bigr\rvert_{\mathbf{S}}\bigr)$.
  \end{itemize}
\end{proof}
\begin{remark}
  \label{rem:3}
  By means of the $\mathfrak{P}$ topology, we could also formulate a result of ``weak semiclassical
  convergence'' of the Segal quantization functor $\mathbb{W}_h$ to a functor of classical observables, at
  least when restricted to cylindrical observables (and tested on semiclassical quantum states).
\end{remark}

\section{Physical applications}
\label{sec:phys-appl}

In this section we outline some applications of the results introduced in \cref{sec:introduction-1}, to
physical systems described by bosonic field theories. Some applications are new, and some are developed in
detail elsewhere; we provide more details for the new ones, and mainly refer to the literature for the
others. The section is divided in two main parts, \cref{sec:nonr-quant-field} dealing with the
semiclassical and mean-field analysis of nonrelativistic (or semirelativistic) systems, and
\cref{sec:relat-quant-field} dealing with semiclassical relativistic systems.

\subsection{Nonrelativistic and semirelativistic Quantum Field Theories}
\label{sec:nonr-quant-field}

Nonrelativistic and semirelativistic quantum field theories are the ones either describing nonrelativistic
particles (invariant with respect to Galilei symmetry transformations), or describing the interaction of
nonrelativistic particles with relativistic force-carrier bosonic fields (invariant with respect to
Lorentz symmetry transformations). In the first case, the quantum field theoretic description comes into
play when the limit of a large number of particles is considered (mean field description), and if
particles can be pumped in or absorbed by the environment. In the second case, the relativistic
force-carriers can be created and destroyed by the interaction, and are thus naturally described by
quantum fields.

\subsubsection{Fock representation}
\label{sec:fock-representation}

Since most of the nonrelativistic and semirelativistic quantum field theories are studied in the Fock
representation, let us briefly recall some well known results about it here, mostly to fix the
notation. Let $\mathfrak{H}$ be a separable complex Hilbert space. As briefly recalled
in~\cref{sec:set-all-cylindrical}, a complex Hilbert space induces a real Hilbert and symplectic space
\emph{via} the forgetful functor $\mathbb{F}\!\mathbb{F}$ from complex to real vector spaces. Therefore,
$\mathbb{F}\!\mathbb{F}\mathfrak{H}$ is a real vector space with the same elements\footnote{If $\eta\in
  \mathfrak{H}$, let us denote by $\mathbb{F}\!\mathbb{F}\eta$ the corresponding element on
  $\mathbb{F}\!\mathbb{F}\mathfrak{H}$.} as $\mathfrak{H}$, complete with respect to the inner product
$\re \langle \,\cdot \, , \,\cdot \, \rangle_{\mathfrak{H}}$ and with symplectic form $\im \langle \,\cdot
\, , \,\cdot \, \rangle_{\mathfrak{H}}$. Let us denote the associated Weyl C*-algebra by
\begin{equation*}
  \mathbb{W}_h(\mathfrak{H}):=\mathbb{W}_h(\mathbb{F}\!\mathbb{F}\mathfrak{H},\im\langle \,\cdot \, , \,\cdot \, \rangle_{\mathfrak{H}})\; .
\end{equation*}
The \emph{Fock vacuum} $\varpi_h$ on $\mathbb{W}_h(\mathfrak{H})$ is the regular quantum state with
generating functional
\begin{align*}
  \mspace{180mu}&\mathcal{G}_{\varpi_h}\bigl(\mathbb{F}\!\mathbb{F}\eta\bigr)=e^{-\frac{h}{2}\lVert \eta  \rVert_{\mathfrak{H}}^2} &(\forall \eta\in \mathfrak{H})\;.
\end{align*}
\begin{definition}[Fock representation]
  \label{def:5}
  The Fock representation of $\mathbb{W}_h(\mathfrak{H})$ is the GNS representation
  $(\mathscr{H}_{\varpi_h},\pi_{\varpi_h},\Omega_h)$ generated by the Fock vacuum. The Hilbert space
  $\mathscr{H}_{\varpi_h}=\Gamma_{\mathrm{s}}(\mathfrak{H})$ is the \emph{symmetric Fock space}, and the
  cyclic vector $\Omega_h\in \Gamma_{\mathrm{s}}(\mathfrak{H})$ is the \emph{(bosonic) Fock vacuum
    vector}.
\end{definition}
The Fock space $\Gamma_{\mathrm{s}}(\mathfrak{H})$ is defined as
\begin{equation*}
  \Gamma_{\mathrm{s}}(\mathfrak{H}):=\mathds{C} \oplus \bigoplus_{n\in \mathds{N}_{*}} \,\bigvee_{j=1}^{n} \mathfrak{H}\;,
\end{equation*}
where $\bigvee$ stands for the symmetric tensor product. The vacuum vector $\Omega_h$ is the vector
\begin{equation*}
  \Omega_h=(1,0,\dotsc,0,\dotsc)\; .
\end{equation*}
The dependence on the semiclassical parameter $h$ is given by the canonical variables of the Fock
representation, the so-called \emph{creation and annihilation operators}. For any $\eta,\xi\in
\mathfrak{H}$, let $a^{*}_h(\eta)$ and $a_h(\xi)$ be the creation and annihilation (closed) operators on
$\Gamma_{\mathrm{s}}(\mathfrak{H})$ satisfying the commutation relations
\begin{equation*}
  [a_h^{*}(\eta),a_h(\xi)]=-h\langle \xi  , \eta \rangle_{\mathfrak{H}}\; .
\end{equation*}
They are explicitly defined by the action on Fock space vectors
$\psi_h=(\psi_{0,h},\psi_{1,h},\dotsc,\psi_{n,h},\dotsc)$ with finitely many non-zero components, and each
component of the form
\begin{align}
  \label{eq:14}
  \psi_{n,h}= \bigvee_{j=1}^n\eta^{(n)}_{j,h}\; .
\end{align}
The set of such vectors is a core for every $a_h^{*}(\eta)$ and $a_h(\eta)$.
\begin{align*}
  &\bigl(a^{*}_h(\eta)\psi_h\bigr)_n=\sqrt{hn}\, \eta\vee \psi_{n-1,h}=\sqrt{hn}\, \eta\vee\bigvee_{j=1}^{n-1} \eta^{(n-1)}_{j,h}\; ;\\
  &\bigl(a_h(\eta)\psi_h\bigr)_n=\sqrt{h(n+1)}\,\langle \eta  , \psi_{n+1,h} \rangle_{\mathfrak{H}}= \sqrt{\frac{h}{n+1}}\,\sum_{k=1}^{n+1}\langle \eta  , \eta^{(n+1)}_{k,h} \rangle_{}\bigvee_{j\neq k}\eta^{(n+1)}_{j,h}\;.
\end{align*}
Intuitively, the spaces $\bigvee_{j=1}^n\mathfrak{H}$ represent the subspaces with $n$ particles (and $\mathds{C}$ the
subspace with no particles). With this interpretation in mind, it is easy to see that by definition the
creation operator $a^{*}_h(\eta)$ creates a particle in the configuration $\eta\in \mathfrak{H}$, and $a_h(\eta)$
annihilates a particle (in the same configuration). The Weyl operators $W_h\bigl(\mathbb{F}\!\mathbb{F} \eta\bigr)$, $\eta\in
\mathfrak{H}$, are Fock-represented by the following unitary operators on $\Gamma_{\mathrm{s}}(\mathfrak{H})$:
\begin{equation*}
  \pi_{\varpi_h}\Bigl(W_h\bigl(\mathbb{F}\!\mathbb{F} \eta\bigr)\Bigr)=e^{\,i\,\bigl(a^{*}_h(\eta)+a_h(\eta)\bigr)}\; .
\end{equation*}

The Fock representation is the most used, since it describes free quantum field theories, both from the
nonrelativistic and relativistic standpoint. Let us also remark that if $\mathfrak{H}\subset
\mathcal{S}'(G)$, for some locally compact abelian group $G$ (\emph{e.g.}, $G=\mathds{R}^d$), then the
Fock representation is also a representation of the Weyl C*-algebra of test functions
$\mathbb{W}_h\bigl(\mathbb{F}\!\mathbb{F} \mathcal{S}(G),\im\langle \,\cdot \, , \,\cdot \,
\rangle_{\mathfrak{H}}\bigr)$, where $\mathcal{S}(G)\subset \mathfrak{H}$ is the nuclear space of rapid
decrease functions on $G$. The non-Fock representations of $\mathbb{W}_h\bigl(\mathbb{F}\!\mathbb{F}
\mathcal{S}(G),\im\langle \,\cdot \, , \,\cdot \, \rangle_{\mathfrak{H}}\bigr)$ play an important role in
interacting relativistic field theories, as it is discussed in \cref{sec:relat-quant-field}.

\subsubsection{Thermodynamic limit of trapped ideal Bose gases}
\label{sec:therm-analys-bose}

A system of $N$ non-relativistic $d$-dimensional bosons in a harmonic trap is usually described by the
following Hamiltonian on $\bigvee_{j=1}^NL^2 (\mathds{R}^d )$:
\begin{equation}
  \label{eq:12}
  H_{N,V_N}=\sum_{j=1}^N(-\Delta_j+\omega_N^2x_j^2) + \frac{1}{N}\sum_{j<k}^{}V_N(x_j-x_k)\; ;
\end{equation}
where $\omega_N\in \mathds{R}$ is proportional to the frequency of the trap, and $V_N$ is a symmetric
two-body interaction potential, with suitable regularity properties yielding the self-adjointness of
$H_{N,V_N}$ on a suitable domain. The two body interaction $V_N$ may be, \emph{e.g.}, independent of $N$
(mean field regime), or of the form $V_N(\cdot )= N^2V(N\,\cdot \,)$ for $d=3$ (Gross-Pitaevskii
regime). The Hamiltonian $H_{N,V_N}$ agrees with the restriction to the $N$-particle sector of a
particle-preserving Hamiltonian on the Fock representation\footnote{The one-particle and two-particle
  second quantizations $\mathrm{d}\Gamma_h$ and $\mathrm{d}\Gamma^{(2)}_h$ are maps from the one or
  two-particle self-adjoint operators on $\mathfrak{H}$ and $\mathfrak{H}\vee \mathfrak{H}$ respectively
  to the self-adjoint operators on $\Gamma_{\mathrm{s}}(\mathfrak{H})$. They are defined by the action on
  core vectors $\psi_h=(\psi_{0,h},\psi_{1,h},\dotsc,\psi_{n,h},\dotsc)$ with finitely many non-zero
  components, and each component of the form \eqref{eq:14} (with each $\eta_{j,h}$ or $\eta_{j,h}\vee
  \eta_{j',h}$ in the domain of the one or two-particle self-adjoint operator):
  \begin{align*}
    &\bigl(\mathrm{d}\Gamma_h(K^{(1)})\psi\bigr)_n=h\sum_{k=1}^n (K^{(1)}\eta_k)\,\vee\, \bigvee_{j\neq k}\eta_j &\Bigl(\forall n\in \mathbb{N}^{*}\,,\, \bigl(\mathrm{d}\Gamma_h(K^{(1)})\psi\bigr)_0=0\Bigr)\\
    &\bigl(\mathrm{d}\Gamma^{(2)}_h(K^{(2)})\psi\bigr)_n=h^2\sum_{k\neq l=1}^n (K^{(2)} \eta_k \vee \eta_l)\,\vee\, \bigvee_{j\neq k,l}\eta_j &\Bigl(\forall n\geq 2\,,\, \bigl(\mathrm{d}\Gamma^{(2)}_h(K^{(2)})\psi\bigr)_0=\\ & &\bigl(\mathrm{d}\Gamma^{(2)}_h(K^{(2)})\psi\bigr)_1=0\Bigr)
  \end{align*}
} $\Gamma_{\mathrm{s}}\bigl(L^2(\mathds{R}^d)\bigr)$, provided that $h=N^{-1}$:
\begin{equation}
  \label{eq:13}
  H_{h,V_h}=\mathrm{d}\Gamma_h\bigl(k_{\omega_{h^{-1}}}\bigr)+\mathrm{d}\Gamma_h^{(2)}\bigl(\tfrac{1}{2}V_{h^{-1}}\bigr)\; ,
\end{equation}
where $k_{\lambda}(\nabla_x,x)=-\Delta_x+\lambda^2x^{2}$ is a differential self-adjoint operator on $L^2
(\mathds{R}^d )$. In this context, the semiclassical parameter $h$ is therefore interpreted as a quantity
proportional to the inverse of the expected number of particles in the system.

Let us consider now a trapped \emph{ideal} nonrelativistic Bose gas in $d$ dimensions. In the Fock
representation, the system is therefore described by the Hamiltonian
\begin{equation*}
  H_{h,0}=\mathrm{d}\Gamma_h(k_{\omega_{h^{-1}}})\;,
\end{equation*}
with $k_{\omega_{h^{-1}}}$ the one-particle harmonic oscillator of frequency $\omega_{h^{-1}}\in
\mathds{R}^+$. The corresponding grand-canonical Gibbs state $\gamma_h$ at inverse temperature $\beta_h>
0$ and chemical potential\footnote{The chemical potential $\mu_h\in \mathds{R}$ should be suitably chosen,
  \emph{i.e.}\ it should be s.t.\ $\beta_h \bigl(k_{\omega_{h^{-1}}}-\mu_h\bigr)> 0$ uniformly with respect to
  $h\in (0,1)$.}  $\mu_h\in \mathds{R}$ has the following form:
\begin{equation*}
  \gamma_h(\cdot )=\frac{\Tr_{\Gamma_{\mathrm{s}}}\Bigl( \;\, \cdot \;\; e^{-\beta_h H_{h,0}+\mu_h\mathrm{d}\Gamma_h(1)}\Bigr)}{\Tr_{\Gamma_{\mathrm{s}}}\Bigl( e^{-\beta_hH_{h,0}+\mu_h\mathrm{d}\Gamma_h(1)}\Bigr)}\; .
\end{equation*}
It corresponds to the generating functional \citep[see][Proposition
5.2.28]{bratteli1997tmp2
}
\begin{equation}
  \label{eq:15}
  \mathcal{G}_{\gamma_h}(\eta)=\exp\biggl(-\tfrac{h}{2}\Bigl\langle \,\eta\,  ,\, \biggl(e^{-\beta_h(k_{\omega_{h^{-1}}}-\mu_h)}\Bigl(1-e^{-\beta_h(k_{\omega_{h^{-1}}}-\mu_h)}\Bigr)^{-1}\biggr)\,\eta\, \Bigr\rangle_2\biggr)\; .
\end{equation}

We are interested in the behavior of the grand-canonical Gibbs state in the thermodynamic limit. The
thermodynamic limit is usually defined as the limit $N,V\to \infty$ ($N$ expected number of particles, $V$
volume), with $\frac{N}{V}=\rho$ constant. In this case, $N=h^{-1}$, and therefore $N\to \infty$ is
equivalent to $h\to 0$. The effective volume occupied by the trapped system depends on the frequency of
the trap, and more precisely it is proportional to $\omega_{h^{-1}}^{-d}$. Therefore, to achieve the
thermodynamic limit, one should relax the trap in the following way: $\omega_{h^{-1}}=\omega
h^{\frac{1}{d}}$, $\omega\in \mathds{R}^+$, as $h\to 0$ \citep[see, \emph{e.g.},][or any other book on
Bose-Einstein
condensation]{pitaevskii2003oup
}. The thermodynamic limit can therefore be reinterpreted as a semiclassical limit $h\to 0$. Now, let
$h\in I\subseteq (0,1)$, with zero adherent to $I$. Since any state normal with respect to the Fock
representation is regular\footnote{A state $\omega$ on a C*-algebra is normal with respect to a given
  representation $\mathscr{H}$ iff $\omega(\cdot )=\Tr_{\mathscr{H}}(\,\cdot \,\,^{\exists}\!\rho)$, with
  $\rho\in \mathfrak{S}^1_+(\mathscr{H})$ (positive and trace class).}, and
\begin{equation*}
  \sup_{h\in I}\mathcal{G}_{\gamma_h}(0)=1\; ,
\end{equation*}
it follows that the grand canonical Gibbs state $\gamma_h$ is a \emph{semiclassical quantum state}, and
thus it has at least one cylindrical Wigner measure associated to it. We can without loss of generality suppose that
\begin{equation*}
  \gamma_h\underset{h\to 0}{\overset{\mathfrak{P}}{\longrightarrow}}G\; \textrm{ (otherwise change }I\textrm{),}
\end{equation*}
with\footnote{To be precise, $M_{\gamma_h}\in
  \mathcal{M}_{\mathrm{cyl}}\bigl(\mathbb{F}\!\mathbb{F}L^2(\mathds{R}^d),\mathbb{F}\!\mathbb{F}L^2
  (\mathds{R}^d)\bigr)$ (with the usual implicit identification $\mathbb{F}\!\mathbb{F}L^2
  (\mathds{R}^d)\cong \mathbb{F}\!\mathbb{F}L^2 (\mathds{R}^d)'$), however
  $\mathcal{M}_{\mathrm{cyl}}\bigl(L^2(\mathds{R}^d),L^2 (\mathds{R}^d)\bigr)$ is the equivalent
  description on the original complex vector spaces.}  $G\in
\mathcal{M}_{\mathrm{cyl}}\bigl(L^2(\mathds{R}^d),L^2(\mathds{R}^d)\bigr)$.

We would like to reinterpret the well-known analysis of condensation for ideal Bose gases using the
thermodynamic (classical) state $G$. 
Let us consider the orthonormal basis $\Bigl\{\varphi^{(m)}_{\omega_{h^{-1}}}\Bigr\}_{m\in \mathds{N}^d}$
of $L^2 (\mathds{R}^d )$ given by the eigenvectors of the harmonic oscillator $k_{\omega_{h^{-1}}}$,
satisfying
\begin{equation*}
  k_{\omega_{h^{-1}}}\varphi^{(m)}_{\omega_{h^{-1}}}=\omega_{h^{-1}}(m_1+\dotsm +m_d+1)\varphi^{(m)}_{\omega_{h^{-1}}}\; ,
\end{equation*}
where we have used the notation $m=(m_1,\dotsc,m_d)$. The quantum observable
$a^{*}_{h}(\varphi^{(m)}_{\omega_{h^{-1}}})a_h(\varphi^{(m)}_{\omega_{h^{-1}}})$, $m\in \mathds{N}^d$,
counts the relative number of particles in the configuration $\varphi^{(m)}_{\omega_{h^{-1}}}$. From
\cref{eq:15}, it is not difficult to see that
\begin{equation}
  \label{eq:16}
  \begin{split}
      \gamma_h\Bigl(a^{*}_{h}(\varphi^{(m)}_{\omega_{h^{-1}}})a_h(\varphi^{(m)}_{\omega_{h^{-1}}})\Bigr)=h\frac{e^{-\beta_h\bigl(\omega_{h^{-1}}(m_1+\dotsc +m_d+1)-\mu_h\bigr)}}{1-e^{-\beta_h\bigl(\omega_{h^{-1}}(m_1+\dotsm +m_d+1)-\mu_h\bigr)}}\\=\frac{h}{e^{\beta_h\bigl(\omega h^{\frac{1}{d}}(m_1+\dotsm +m_d+1)-\mu_h\bigr)}-1} \; .
  \end{split}
\end{equation}
Let us recall that $h$ is interpreted as the inverse of the expected number of particles, \emph{i.e.}\ it
satisfies the consistency condition
\begin{equation*}
  h^{-1}=\gamma_h\Bigl(\sum_{m\in \mathds{N}^d}^{}a^{*}_1(\varphi^{(m)}_{\omega_{h^{-1}}})a_1(\varphi^{(m)}_{\omega_{h^{-1}}})\Bigr)\; .
\end{equation*}
Therefore, the quantity in \cref{eq:16} is a number between zero and one uniformly in $h$ (as expected).

Now, if $m=(0,\dotsc,0)$, \cref{eq:16} counts the relative number of particles in the ground state. If, in
the thermodynamic limit $h\to 0$, this number is bigger than zero, then a macroscopic fraction of
particles is in the ground state, and thus the system exhibits condensation. We characterize such
macroscopic fraction of particles in the condensate indirectly\footnote{The indirect approach is customary
  in physics. In addition, it is related to the fact that, in the thermodynamic limit, the measure $G$
  loses mass, and such mass corresponds to the condensed fraction \citep[see][for additional
  details]{ammari2017arxiv
  }.}: we would like to find suitable cylindrical symbols $s_m$, $m\neq (0,\dotsc,0)$, such that
\begin{equation}
  \label{eq:17}
  \lim_{h\to 0}\gamma_h\Bigl(\sum_{m\neq (0,\dotsc,0)}^{}a^{*}_{h}(\varphi^{(m)}_{\omega_{h^{-1}}})a_h(\varphi^{(m)}_{\omega_{h^{-1}}})\Bigr)=\sum_{m\neq (0,\dotsc,0)}^{}\int_{L^2(\mathds{R}^d)}^{}s_m(z)  \mathrm{d}G(z)\;.
\end{equation}
If $\varphi^{(m)}_{\omega_{h^{-1}}}$ was independent of $h$, the symbol would be
\begin{equation*}
  s_m(z)=\Bigl\lvert \langle \varphi^{(m)}_{\omega_{h^{-1}}}\,  , \, z \rangle_2  \Bigr\rvert_{}^2\; .
\end{equation*}
In the thermodynamic limit, however, we are relaxing the harmonic trap, hence
$\varphi^{(m)}_{\omega_{h^{-1}}}$ depends on $h$, through $\omega_{h^{-1}}\underset{h\to
  0}{\longrightarrow} 0$. On the other hand, each $\varphi^{(m)}_{\omega_{h^{-1}}}$ belongs to the unit
ball of $L^2 (\mathds{R}^d )$ uniformly with respect to $h\in I$, and thus by Banach-Alaoglu's theorem it has a
weak limit $\varphi^{(m)}_0\in L^2(\mathds{R}^d)$ (up to a redefinition of $I$), and therefore
\begin{equation*}
  s_m(z)=\Bigl\lvert \langle \varphi^{(m)}_0  \,,\, z \rangle_2  \Bigr\rvert^2\; .
\end{equation*}
The quantity
\begin{equation*}
  f_0(G)=1-\sum_{m\neq (0,\dotsc,0)}^{}\int_{L^2(\mathds{R}^d)}^{}s_m(z)  \mathrm{d}G(z)
\end{equation*}
$0\leq f_0(G)\leq 1$, is the \emph{fraction of condensed particles in the (classical) thermodynamic state
  $G$}. Such fraction is easily calculated using the following physical assumptions: $\beta_h=\beta
h^{\frac{d-1}{d}}$, and $\mu_h=\omega h^{\frac{1}{d}}$, $\beta\in \mathds{R}^+$. With these assumptions
the thermodynamic state $G$ only depends on the macroscopic inverse temperature $\beta$, and thus we
denote it by $G_{\beta}$. It is well-known that there is a $\beta_{*}\in \mathds{R}^{+}$ (the inverse of
the critical temperature) for which there is a phase transition: $f_0(G_{\beta})=0$ for any $\beta\leq
\beta_{*}$ (no condensation), and $f_0(G_{\beta})>0$ (condensation) for any $\beta>\beta_{*}$, with
$f_0(G_{\beta})\underset{\beta\to \infty}{\longrightarrow} 1$ (complete condensation at temperature zero).

\subsubsection{High-temperature limit and Gibbs measures}
\label{sec:high-temp-limit}

The effective behavior of grand-canonical quantum Gibbs states as classical measures has also been
recently studied, in order to provide a microscopic derivation of Gibbs measures
\citep{lewin2015jepm
  ,frohlich2017cmp
  ,lewin2017arxiv
  ,frohlich2017arxiv
}. The scaling used in these works is, however, very different from the one used in
\cref{sec:therm-analys-bose} above: they analyze a \emph{high-temperature limit}, in which there is no
condensation (and thus the classical state could be a Gibbs measure, that has no atoms, \emph{i.e.}\ no
macroscopic occupation of the ground state). Nonetheless, their results could be reinterpreted in our
general framework as the convergence (in the $\mathfrak{P}\vee\mathfrak{T}$ topology) of grand canonical
Gibbs states to their unique cylindrical Wigner measures, the corresponding classical Gibbs
measures. Without entering into the details, one interesting fact is that in dimension $d=2,3$, the limit
measures are, already in the non-interacting case, truly cylindrical in $L^2 (\Lambda_d)$ (where
$\Lambda_d$ is either $\mathds{R}^d$ or the $d$-dimensional torus), and concentrated as Radon measures in
$H^{r(d)}(\Lambda_d)\smallsetminus L^2 (\Lambda^d) $, for a suitable $r(d)<0$. Of course by \cref{thm:5}
it is known that there are semiclassical quantum states that converge to true cylindrical measures, but
these grand-canonical Gibbs states provide a physically relevant example (other such examples will be
given in \cref{sec:rough-effect-potent} below).

\subsubsection{Rough effective potentials in the quasi-classical limit of semi-relativistic theories}
\label{sec:rough-effect-potent}

Using semiclassical analysis for semirelativistic quantum field theories it is possible to study how
controllable effective external potentials acting on quantum particles can be created by making the latter
interact with semiclassical radiation fields
\citep{correggi2017ahp
  ,correggi2017arxiv
}, as it is common practice in experimental physics. In order to produce interesting potentials, such as
harmonic traps and uniform magnetic potentials, it is necessary to use semiclassical quantum states whose
Wigner measures are \emph{truly cylindrical}. This is due to the fact that the effective potentials
generated by semiclassical quantum states whose cylindrical Wigner measures are Radon measures cannot be
singular functions (they are continuous and vanishing at infinity).

Let us briefly overview the main ideas. The interaction of nonrelativistic particles with radiation is
modeled by composite systems described by self-adjoint Hamiltonian operators acting on spaces of the form
$\mathscr{H}\otimes \Gamma_{\mathrm{s}}(\mathfrak{H})$, with $\mathscr{H},\mathfrak{H}$ separable
Hilbert. The \emph{quasi-classical approximation} is the regime in which only the field behaves
semiclassically, while the quantum nature of the particles is still relevant. Mathematically, it amounts
to say that the semiclassical parameter $h$ appears only in the canonical commutation relations of the
field, \emph{i.e.},
\begin{align*}
  \mspace{150mu}&\bigl[1\otimes a_h^{*}(\eta),1\otimes a_h(\xi)\bigr]=-h\langle \eta  , \xi \rangle_{\mathfrak{H}}&(\forall \eta,\xi\in \mathfrak{H})\;.
\end{align*}
Fock-normal semiclassical states are therefore represented by density matrices $\rho_h\in
\mathfrak{S}_+^1(\mathscr{H}\otimes \Gamma_{\mathrm{s}}(\mathfrak{H}))$, with trace uniformly bounded
with respect to $h\in I\subseteq (0,1)$. The additional degrees of freedom given by $\mathscr{H}$ are reflected in
the fact that the Wigner measures associated to $\Tr\bigl(\,\cdot \, \rho_h\bigr)$ are in general not
scalar-valued, but vector-valued (with values in the Banach cone of positive states on the C*-algebra
$\mathcal{L}(\mathscr{H})$). This generalization is taken into account below in
\crefrange{sec:class-char-stat}{sec:states-c-algebra}, where we prove (generalizations of) the results
discussed in \cref{sec:introduction-1}. However, if the semiclassical quantum state is of the form
\begin{align*}
  \bigl\lvert \,\!^{\exists }\Psi\otimes\,\!^{\exists }\psi_h\bigr\rangle\bigl\langle\,\Psi\otimes \psi_h\, \bigr\rvert\; ,\; \Psi\in \mathscr{H}\,,\, \psi_h\in \Gamma_{\mathrm{s}}(\mathfrak{H})\; ,
\end{align*}
then the corresponding cylindrical measure factorizes to the product $\Tr_{\mathscr{H}}\bigl(\; \cdot \;
\lvert\Psi\rangle\langle\Psi\rvert\bigr) \times M$, where $M$ is the scalar cylindrical Wigner measure of
$\Tr_{\Gamma_{\mathrm{s}}}\bigl(\; \cdot \; \lvert \psi_h\rangle\langle\psi_h \rvert\bigr)$. Now, let
$\psi_h\in \Gamma_{\mathrm{s}}(\mathfrak{H})$ be a semiclassical quantum vector (\emph{i.e.},
$\Tr_{\Gamma_{\mathrm{s}}}\bigl(\; \cdot \; \lvert \psi_h\rangle\langle\psi_h \rvert\bigr)$ is a
semiclassical quantum state) converging to $M\in \mathcal{M}_{\mathrm{cyl}}(\mathfrak{H},\mathfrak{H})$ in
the $\mathfrak{P}\vee\mathfrak{T}$ topology. As discussed above we are interested in studying the
quasi-classical behavior of a composite system consisting of particles and a bosonic field; in particular
we would like to characterize the effects of the interaction on the subsystem of particles. A simple model
for the interaction is the following. Let $\mathscr{H}=L^2 (\mathds{R}^{Nd} )$ be the Hilbert space of $N$
nonrelativistic particles (without spin), and let
\begin{equation*}
  H_h=\bigl(-\Delta+V(x)\bigr)\otimes 1 + 1\otimes \mathrm{d}\Gamma_h(k)+a_h^{*}\bigl(\lambda(x)\bigr)+a_h\bigl(\lambda(x)\bigr)\; ,
\end{equation*}
where $V$, when considered as a multiplication operator, is a small Kato perturbation of $-\Delta$, $k\geq
0$ is self-adjoint on $\mathfrak{H}$, and $\lambda\in L^{\infty}(\mathds{R}^{Nd},\mathfrak{H})$. These
types of models were first introduced in the community of mathematical physics by
\citet{nelson1964jmp
}. If we take the partial trace of $H_h$ with respect to $\lvert \psi_h\rangle\langle\psi_h \rvert$, we
obtain an operator $\langle H_h \rangle_{\psi_h}$ acting only on $\mathscr{H}=L^2 (\mathds{R}^{Nd} )$. In
addition, we subtract from $\langle H_h \rangle_{\psi_h}$ the multiple of the identity $c_h=\langle
\mathrm{d}\Gamma_h(k) \rangle_{\psi_h}$, since it only amounts to a spectral shift. Then in the
quasi-classical limit $h\to 0$ we are able to prove in
\citep{correggi2017ahp
} that, as long as the semiclassical vector $\psi_h$ is sufficiently regular\footnote{In particular, we
  have that $M=\mu\in \mathcal{M}_{\mathrm{rad}}(\mathfrak{H})$, with $\int_{\mathfrak{H}}^{}\lVert z
  \rVert_{\mathfrak{H}}^{} \mathrm{d}\mu(z)<\infty$.}, the following convergence holds in the norm
resolvent sense:
\begin{equation*}
  \langle H_h \rangle_{\psi_h}-c_h\;\underset{h\to 0}{\overset{\textrm{norm-res}}{\longrightarrow}}\; H(\mu)=-\Delta+V(x)+2\re\int_{\mathfrak{H}}^{}\langle \lambda(x)  , z \rangle_{\mathfrak{H}}  \mathrm{d}\mu(z)\; .
\end{equation*}
Therefore, the interaction with the field generates, in the quasi-classical limit, an effective potential
$2\re\int_{\mathfrak{H}}^{}\langle \lambda(x) , z \rangle_{\mathfrak{H}} \mathrm{d}\mu(z)$ acting on the
particles. Such potential is controllable by tuning the quantum field configuration $\psi_h$. One drawback
in requiring that the Wigner measure of $\psi_h$ is a Radon measure on $\mathfrak{H}$ is that it is not
possible to generate ``physically interesting'' potentials, \emph{e.g.}, unbounded ones that could trap
the particles (such as $\omega^2x^2$, that played an important role in \cref{sec:therm-analys-bose} to
study condensation). In fact, for every $\mu\in \mathcal{M}_{\mathrm{rad}}(\mathfrak{H})$ the effective
potential is either continuous and vanishing at infinity (if $\lVert \,\cdot \, \rVert_{\mathfrak{H}}^{}$
is $\mu$-integrable), or undefined (if $\int_{\mathfrak{H}}^{} \lVert z
\rVert_{\mathfrak{H}}^{}\mathrm{d}\mu(z)$ diverges). In order to obtain interesting potentials it is
therefore necessary to consider semiclassical vectors whose cylindrical Wigner measures \emph{are not}
Radon measures on $\mathfrak{H}$. In fact, it is, \emph{e.g.}, possible to prove that given any potential
$W\in L^2_{\mathrm{loc}} (\mathds{R}^{Nd},\mathbb{R}_+)$, then there exists at least one semiclassical
vector $\psi_{h,W}\in \Gamma_{\mathrm{s}}\bigl(L^2 (\mathds{R}^d )\bigr)$, with explicit form and
converging to a ``true'' cylindrical measure, and (infinitely many) coupling functions $\lambda_W\in
L^{\infty}\bigl(\mathbb{R}^{Nd},L^2 (\mathds{R}^d )\bigr)$ such that
\begin{equation*}
  \langle H_{h,\lambda_W} \rangle_{\psi_{h,W}}-c_h\;\underset{h\to 0}{\overset{\textrm{strong-res}}{\longrightarrow}}\; H_W=-\Delta+V(x)+W(x)\; .
\end{equation*}
Therefore this procedure describes a simple way of producing any given external potential, acting on
quantum particles, exploiting their interaction with a semiclassical bosonic field (\emph{e.g.}, a phonon
field, a radiation field, \dots); and it relies on semiclassical states whose semiclassical measures are
truly cylindrical to produce strong potentials.

\subsection{Relativistic Quantum Field Theories}
\label{sec:relat-quant-field}

If the applications described in \cref{sec:nonr-quant-field} motivate the importance of cylindrical
measures as Wigner measures, the ones described in this section motivate the abstract, representation
independent, algebraic approach to semiclassical states taken in \cref{sec:introduction-1}. In fact, there
are features of relativistic quantum field theories that make such an approach necessary. One is that the
phase space is often taken to be a space of test functions, such as the nuclear space $\mathcal{S}(G)$ or
$\mathcal{D}(G)=C_0^{\infty}(G)$ with the inductive limit topology (where $G$ in both cases is some
locally compact abelian group). Such choice actually plays an important role in defining the right
representation of the canonical commutation relations, at least in the few rigorously definable
interacting theories that we know of \citep[see, \emph{e.g.},][and references thereof
contained]{glimm1987qp
}. It is therefore too restrictive, in relativistic quantum field theories, to assume the phase spaces to
be (pre)Hilbert (with the natural induced topology), and to focus solely on Fock-normal states. In
addition, since the phase space is considered, from the physical standpoint, as a space of test functions,
it becomes clear why the natural space of classical fields (the one on which the cylindrical Wigner
measures should act upon) is a space of ``distributions'', \emph{i.e.}\ a space in duality with the phase
space\footnote{In the aforementioned examples, it would naturally be respectively $\mathcal{S}'(G)$ with
  the $\sigma\bigl(\mathcal{S}'(G),\mathcal{S}(G)\bigr)$ topology, and $\mathcal{D}'(G)$ with the
  analogous weak topology.}. Another feature that emphasizes the importance of a
representation-independent analysis of the semiclassical states is the so-called \emph{Haag's theorem}
\citep[see, \emph{e.g.},][]{haag1992tmp
}, asserting that if we consider two relativistic invariant states\footnote{A typical example would be to
  consider the ground state of a non-interacting theory (the Fock vacuum for the scalar field), and the
  ground state of an interacting theory (such as $(\varphi^4)_2$).} of the same Weyl C*-algebra
(\emph{i.e.}\ two states that are invariant with respect to the action of the Poincaré group on the C*-algebra),
then they are either equal or disjoint\footnote{Two states $\omega_1,\omega_2$ on a C*-algebra are
  disjoint iff $\omega_1$ is not normal with respect to the GNS representation of $\omega_2$, and vice-versa. It
  follows that the GNS representations given by $\omega_1$ and $\omega_2$ are \emph{inequivalent},
  \emph{i.e.}\ they are not related by a unitary isomorphism.}. Since the ground states of a free and an
interacting theory cannot be the same, this means that a relativistically covariant interacting theory
should be in a representation that is inequivalent to the free (Fock) representation. Since the results of
\cref{sec:introduction-1} are all \emph{representation-independent}, they are well suited for application
to relativistic quantum field theories.

\subsubsection{Semiclassical Axiomatic Quantum Field Theory}
\label{sec:axiom-quant-field}

There are various axiomatic formulations of relativistic quantum field theories. The most used sets of
axioms are the ``Hilbert space'' axioms known as Gårding-Wightman
\citep{wightman1965afksv
  ,streater2000plp
}, and the ``algebraic'' axioms known as Haag-Kastler
\citep{haag1964jmp
  ,haag1992tmp
}. Due to the algebraic formulation of semiclassical analysis given in \cref{sec:introduction-1}, it is
natural for us to choose algebraic axioms as a starting point. In order to describe also quantum field
theories in non-Minkowski spacetimes, we use the so-called \emph{locally covariant} algebraic axioms
\citep[see,
\emph{e.g.},][]{brunetti2009lnp
  ,brunetti2014rmp
  ,fredenhagen2016jmp
}.

The starting point of the locally covariant approach to algebraic bosonic quantum field theory is a
\emph{local bosonic quantization functor}. Let $\mathbf{A}$ be a \emph{small} category, to be interpreted
as the category that has as objects local spacetimes, \emph{e.g.}, local (bounded) regions of a given
spacetime, and as morphisms suitably regular mappings between local spacetimes. The local quantization
functor is then the composition of a functor from $\mathbf{A}$ to (local) phase spaces (more accurately,
to spaces of local test functions) and the Segal bosonic quantization functor.

\begin{definition}[Local bosonic quantization functor]
  \label{def:6}
  A functor $\mathbb{L}\!\mathbb{W}_h:\mathbf{A}\to \mathbf{C^{*}alg}$ is a local bosonic quantization
  functor iff
  \begin{align*}
    \mspace{210mu}&\mathbb{L}\!\mathbb{W}_h= \mathbb{W}_h\circ \,\!^{\exists }\mathbb{L} &\bigl(\mathbb{L}:\mathbf{A}\to \mathbf{Symp}_{\mathds{R}}\bigr)\; .
  \end{align*}
\end{definition}
The local bosonic quantization functor associates to each local spacetime a C*-algebra of canonical
commutation relations, the local bosonic observables\footnote{The algebras of local bosonic observables
  are in general taken to be larger, containing the algebras of canonical commutation relations as
  subalgebras. However, since we are interested in semiclassical regular quantum states, we can without loss of generality
  restrict the observables to the ones in the Weyl C*-algebra.}. By another composition with the duality
functor $\mathbb{D}_+: \mathbf{C^{*}alg}\to \mathbf{BanCone}$, defined in \cref{sec:categ-interpr}, we
obtain the \emph{functor of local bosonic quantum states}
\begin{equation*}
  \mathbb{L}\mspace{-1mu}\mathbb{S}_h:=\mathbb{D}_+\circ \mathbb{L}\!\mathbb{W}_h=\mathbb{D}_+\circ \mathbb{W}_h\circ \mathbb{L}\; .
\end{equation*}
Clearly, the results formulated in \cref{sec:categ-interpr} for $\mathbb{S}_h$ hold, \emph{mutatis
  mutandis}, for $\mathbb{L}\mspace{-1mu}\mathbb{S}_h$. To this extent, let us define a \emph{functor of local
  cylindrical Wigner measures} to be any functor of the form
\begin{equation*}
  \mathbb{L}\mspace{-1mu}\mathbb{S}_0:=\mathbb{S}_0\circ \mathbb{L}\; ,
\end{equation*}
where $\mathbb{S}_0$ is any functor of cylindrical Wigner measures, see \cref{def:4}. By \cref{prop:3}, it
follows that for any $\mathbb{L}$ a functor of local cylindrical Wigner measures exists and it is unique
up to natural isomorphisms. In addition, \cref{th:de} yields
\begin{equation}
  \label{eq:18}
  \mathbb{L}\mspace{-1mu}\mathbb{S}_h\underset{h\to 0}{\rightarrowtail} \mathbb{L}\mspace{-1mu}\mathbb{S}_0\; .
\end{equation}

It is now time to introduce the four axioms of locally covariant algebraic quantum field theory:
\emph{isotony}, \emph{covariance}, \emph{Einstein causality}, and \emph{time-slice}. The \emph{spectrum
  condition}, that together with the existence of a preferred \emph{ground state} (or vacuum) is very
important in Minkowski spacetime, is locally implemented in suitable non-Minkowski spacetimes using
\emph{Hadamard states}. We discuss the semiclassical behavior of ground states in
\cref{sec:semicl-kms-ground}.

In a relativistic theory, the quantum fields should be locally covariant. In the Gårding-Wightman
formulation (in Minkowski global spacetime), this corresponds to requesting the existence of Poincaré
invariant operator-valued field distributions \citep[see \emph{e.g.}][IX.8 Property
6]{reed1975II
}. In the locally covariant algebraic formulation, the covariance axiom has a very simple form.
\begin{equation}
  \label{eq:20}\tag{Cov}
  \mathbb{L} \textrm{ is a covariant functor.}
\end{equation}
An immediate consequence of \eqref{eq:20} is that $\mathbb{L}\!\mathbb{W}_h$ is a covariant functor, and
that $\mathbb{L}\mspace{-1mu}\mathbb{S}_h$ is a contravariant functor (since $\mathbb{D}_+$ is contravariant).

The isotony axiom formalizes the fact that observables of a local spacetime region are also observables of
any region that includes the former. In mathematical terms, the functor of local quantum observables
preserves injective morphisms $\hookrightarrow$ (embeddings):
\begin{align}
  \label{eq:19}\tag{Iso}
  \mspace{53mu}&\mathfrak{a}: A\hookrightarrow B\; \Rightarrow\; \mathbb{L}(\mathfrak{a}): \mathbb{L}(A)\hookrightarrow \mathbb{L}(B)&\bigl(\forall \mathfrak{a} \in \mathrm{Morph}(\mathbf{A})\bigr)\; .
\end{align}
\begin{lemma}
  \label{lemma:1}
  If the \eqref{eq:19} axiom holds, then $\mathbb{L}\!\mathbb{W}_h(\mathfrak{a}):
  \mathbb{L}\!\mathbb{W}_h(A) \hookrightarrow \mathbb{L}\!\mathbb{W}_h(B)$ for any embedding
  $\mathfrak{a}: A\hookrightarrow B$ of $\mathbf{A}$.
\end{lemma}
\begin{proof}
  $\mathbb{L}\!\mathbb{W}_h(\mathfrak{a})$ is the Weyl-operator-preserving *-homomorphism of Weyl
  C*-algebras induced by an injective homomorphism of Heisenberg groups. More precisely, it is induced by
  the group homomorphism between $\mathbb{H}\bigl(\mathbb{L}(A)\bigr)$ and
  $\mathbb{H}\bigl(\mathbb{L}(B)\bigr)$, in turn induced by the injective symplectomorphism
  $\mathbb{L}(\mathfrak{a})$. Therefore, $\mathbb{L}\!\mathbb{W}_h(\mathfrak{a})$ is injective as well.
\end{proof}

In relativistic theories, quantum fields should also respect causality. The best-known formulation of the
causality axioms is that local observables belonging to space-like separated regions of spacetime should
always commute. As showed by \citet{roos1970cmp
}, given two commuting subalgebras of observables statistical independence is equivalent to a
non-vanishing condition on the product of elements from each different subalgebra, that is always
satisfied by tensor products (provided the embedding algebra also embeds the tensor product of the two
subalgebras). In fact, the sensible notion of Einstein causality in the locally covariant formulation
turns out to be the fact that the local functor $\mathbb{L}$ \emph{preserves tensor structures}. A
category is \emph{monoidal} (or tensor) if there exists a monoidal structure on it\footnote{More
  precisely, a category $\mathbf{C}$ is \emph{monoidal} (usually denoted by $\mathbf{C}^{\otimes }$) if
  there exists a monoidal structure:
  \begin{itemize}
  \item $\exists \, \otimes : \mathbf{C}\times \mathbf{C} \to \mathbf{C}$ (tensor or monoidal product);
  
  \item $\exists \, I\in \mathbf{C}$ (identity);
  
  \item $\exists\, \nu^{\mathrm{ass}},\nu^{\mathrm{l}},\nu^{\mathrm{r}}$ natural isomorphisms with
    components:
    \begin{itemize}
    \item $\nu^{\mathrm{ass}}_{A,B,C}:(A\otimes B)\otimes C \overset{\mathrm{iso}}{\longrightarrow}
      A\otimes (B\otimes C)$,
    
    \item $\nu^{\mathrm{l}}_A:I\otimes A \overset{\mathrm{iso}}{\longrightarrow} A$,
    
    \item $\nu^{\mathrm{r}}_A:A\otimes I \overset{\mathrm{iso}}{\longrightarrow} A$.
    \end{itemize}
  \end{itemize}}. A functor that preserves monoidal structures is called \emph{homomorphic}, or
\emph{monoidal}\footnote{A functor $\mathbb{F}^{\otimes } :\mathbf{C}^{\otimes }\to \mathbf{D}^{\otimes }$
  is homomorphic iff:
  \begin{itemize}
  \item $\mathbb{F}^{\otimes }(A\otimes B) \cong \mathbb{F}^{\otimes }(A)\otimes \mathbb{F}^{\otimes
    }(B)$;
  
  \item $\mathbb{F}^{\otimes }\bigl(\mathfrak{c}_1\otimes \mathfrak{c}_2:A\otimes B\to C\otimes D\bigr)
    \cong \mathbb{F}^{\otimes }(\mathfrak{c}_1)\otimes \mathbb{F}^{\otimes }(\mathfrak{c}_2);$
  
  \item $\mathbb{F}^{\otimes }(I_{\mathbf{C}}) \cong I_{\mathbf{D}}$.
  \end{itemize}}. On local spacetimes there is a natural monoidal structure, therefore we suppose that
$\mathbf{A}^{\otimes }$ is monoidal. On symplectic spaces, the monoidal structure is given by direct
products: $(X,\varsigma)\otimes (Y,\tau)=(X\oplus Y,\varsigma \oplus \tau)$, for any
$(X,\varsigma),(Y,\tau)\in \mathbf{Symp}_{\mathds{R}}$ (with the space $\{0\}$ as identity).
\begin{equation}
  \label{eq:21}\tag{Ein}
  \mathbb{L}^{\otimes }: \mathbf{A}^{\otimes }\to \mathbf{Symp}_{\mathds{R}}^{\otimes } \; \textrm{ is homomorphic.}
\end{equation}
On C*-algebras, there is a directed set of possible monoidal structures, each one corresponding to a
different choice of C*-cross-norm \citep[see][for an introduction to tensor products on
C*-algebras]{takesaki1979I
}. For every monoidal structure on $\mathbf{C^{*}alg}$ the identity is $\mathds{C}\in
\mathbf{C^{*}alg}$. Clearly, we denote by $\mathbf{C^{*}alg}^{\otimes_{\alpha} }$ the monoidal category of
C*-algebras with a given choice $\otimes_{\alpha}$ of C*-cross-norm ($\alpha_{\mathrm{min}}\leq \alpha\leq
\alpha_{\mathrm{max}}$, where $\alpha_{\mathrm{min}}$ and $\alpha_{\mathrm{max}}$ are the injective and
projective C*-cross-norms respectively). This also induces a tensor structure on the space of states, and
the functor $\mathbb{D}_+$ is homomorphic. It is also not difficult to see that the functor $\mathbb{W}_h$
is homomorphic for any $h\in (0,1)$, and for any choice of monoidal structure
$\mathbf{C^{*}alg}^{\otimes_{\alpha} }$. Therefore the following lemma is true.
\begin{lemma}
  \label{lemma:2}
  If \eqref{eq:21} holds, then both $\mathbb{L}\!\mathbb{W}_h^{\otimes_{\alpha}}$ and
  $\mathbb{L}\mspace{-1mu}\mathbb{S}_h^{\otimes_{\alpha}}$ are homomorphic for any $\alpha_{\mathrm{min}}\leq
  \alpha\leq \alpha_{\mathrm{max}}$.
\end{lemma}
\noindent Let us remark that the formulation above is \emph{equivalent} to the requirement that space-like
separated local observables commute if the C*-cross-norm is the inductive one $\alpha_{\mathrm{min}}$
\citep{brunetti2014rmp
}.

The last axiom to be introduced is time-slice. The idea behind the time slice axiom is that a quantum
field theory should be a quantum evolution theory, \emph{i.e.}\ it should be determined completely by
fields ``at one fixed time'', and by the action on them of an evolution operator. In curved spacetimes,
the concept of space at a fixed time is given by \emph{Cauchy surfaces}, that is surfaces that are
intersected by every non-extensible, causal spacetime path exactly once. Since we do not use here the
properties of Cauchy surfaces, we do not discuss the conditions a spacetime has to satisfy to guarantee
their existence. Given the small category of local spacetimes $\mathbf{A}$, supposed to be based on sets,
we define a Cauchy surface of $\mathbf{A}$ to be any set $\Sigma$ with the following property
\begin{equation}
  \label{eq:22}
  \Sigma=\bigcap_{\Sigma \subset A \in \mathbf{A}} A\neq \emptyset \; .
\end{equation}
Let us also suppose that the set of objects of $\mathbf{A}$ that include $\Sigma$ is directed, when
ordered by inclusion. Then if we denote by $\mathfrak{i}_{AB}:A\hookleftarrow B$, $A\supseteq B$, the
canonical inclusion morphisms, it follows that
\begin{equation*}
  \bigl(A\supset \Sigma, \mathfrak{i}_{AB} \bigr)_{A\supseteq B \in \mathbf{A}}
\end{equation*}
is a projective family that satisfies
\begin{equation*}
  \varprojlim_{\Sigma \subset A \in \mathbf{A}} A= \bigcap_{\Sigma \subset A \in \mathbf{A}} A=\Sigma\; .
\end{equation*}
By \eqref{eq:20} and \eqref{eq:19}, it also follows that
$\mathbb{L}\!\mathbb{W}_h(\mathfrak{i}_{AB}):\mathbb{L}\!\mathbb{W}_h(A)\hookleftarrow
\mathbb{L}\!\mathbb{W}_h(B)$ are embeddings, whenever $A\supseteq B$. It appears then natural to define
the algebra of observables in the Cauchy surface $\Sigma$ as the C*-projective limit of the algebras of
observables on the local spacetimes containing $\Sigma$ \citep[see,
\emph{e.g.},][]{fredenhagen2011qftg
}\footnote{More generally, the algebra of observables in the Cauchy surface should contain $\varprojlim
  \mathbb{L}\!\mathbb{W}_h(a)$ as a subalgebra.}:
\begin{equation}
  \label{eq:23}
  \mathfrak{O}(\Sigma):=\varprojlim_{\Sigma \subset A \in \mathbf{A}} \mathbb{L}\!\mathbb{W}_h(A) \; .
\end{equation}
One should however be careful since such limit may be the empty set. By construction, the algebra of
observables $\mathfrak{O}(\Sigma)$ has natural projections $P_{A,\Sigma}:\mathfrak{O}(\Sigma)\to
\mathbb{L}\!\mathbb{W}_h(A)$. The time slice axiom then consists of two parts: the first is that there
exists a Weyl C*-algebra $\mathbb{W}_h(X_{\Sigma},\varsigma_{\Sigma})$ of time-sliced (or time-zero)
bosonic fields in the Cauchy surface (this \emph{a fortiori} also guarantees that the algebra of
observables in the Cauchy surface is not empty); and the second is that all projections $P_A$ are
*-isomorphisms.
\begin{align}
  \label{eq:24}\tag{TS}
  \begin{aligned}
      &\mathbb{W}_h(\,\!^{\exists }X_{\Sigma},\,\!^{\exists }\varsigma_{\Sigma})\overset{^{\exists }\mathrm{w}_{\Sigma}}{\hookrightarrow}\mathfrak{O}(\Sigma)\;\textrm{ and } P_{A,\Sigma} \textrm{ is a *-isomorphism} &\bigl(\forall \Sigma \in \mathrm{Cau}_{\mathbf{A}}\;,\; \\ &  & \forall \mathbf{A}\ni A\supset \Sigma \bigr)\;.
  \end{aligned}
\end{align}
In \crefrange{sec:semicl-covar-isot}{sec:semicl-time-slice} we discuss the semiclassical implications of
the four axioms above. For the reader's convenience, such implications are summarized schematically in
\cref{tab:1}.

\begin{table}[h]
  \centering
  \resizebox{\textwidth}{!}{
    \begin{tabular}{c l | r }
      \toprule
      \multicolumn{1}{c}{{\bfseries Quantum property}}&\multicolumn{1}{l}{{\bfseries Semiclassical consequence}}&\multicolumn{1}{c}{{\bfseries Ref.}}\\ 
      \bottomrule
      \noalign{\vskip 2mm}
      Local structure& Loc.\ functor of cl.\ fields; convergence of loc.\ states and linear maps (q.\ to cl.); surjectivity of convergence&~\cref{eq:18}\\[2mm]
      \eqref{eq:20}& Covariance of loc.\ cl.\ fields; convergence of spacetime transformations represented on states (q.\ to cl.)&\cref{sec:semicl-covar-isot}\\[2mm]
      \eqref{eq:19}& Projective structure of loc.\ cl.\ states; definition of global cl.\ states as projective limits&~\cref{sec:semicl-covar-isot}\\[2mm]
      \eqref{eq:21}& Tensoriality of cl.\ spacelike-separated fields; at most destruction of spacelike entanglement by cl.\ limit &\cref{prop:5}\\[2mm]
      \eqref{eq:24} & Definition of cl.\ time-sliced fields; (case-by-case) characterization of cl.\ flow, and Egorov theorems &~\cref{sec:semicl-time-slice}\\
      \noalign{\vskip 1.5mm}
      \midrule
      \noalign{\vskip 1mm}
    \end{tabular}}
  \caption{Semiclassical consequences of (locally covariant) QFT axioms for bosonic theories. List of relevant abbreviations: (Iso) -- isotony axiom; (Cov) -- covariance axiom; (Ein) -- Einstein causality; (TS) -- time-slice axiom; loc.\ -- local; cl.\ -- classical; q.\ -- quantum; stat.\ -- statistical.}
  \label{tab:1}
\end{table}

\subsubsection{Covariance and isotony}
\label{sec:semicl-covar-isot}

The covariance axiom \eqref{eq:20} has the semiclassical consequence that the classical limit fields
behave properly with respect to spacetime transformations. More precisely, consider the convergence
$\mathbb{L}\mspace{-1mu}\mathbb{S}_h\underset{h\to 0}{\rightarrowtail} \mathbb{L}\mspace{-1mu}\mathbb{S}_0$, applied to a
spacetime morphism $\mathfrak{a}:A\to B$, $A,B\in \mathbf{A}$. The latter could be interpreted,
\emph{e.g.}, as a spacetime symmetry transformation, that induces a covariant transformation
$\mathbb{L}(\mathfrak{a})$ of local test functions, and in turn by quantization a transformation of
quantum observables $\mathbb{L}\!\mathbb{W}_h(\mathfrak{a})$. By duality, it also induces a
contravariant transformation of quantum states $\mathbb{L}\mspace{-1mu}\mathbb{S}_h(\mathfrak{a})$. At the classical
level, it induces a transformation $\mathbb{L}\mspace{-1mu}\mathbb{S}_0(\mathfrak{a})$ of cylindrical Wigner
measures. The transformed semiclassical quantum states converge, as they should, to cylindrical Wigner
measures pushed forward by the corresponding classical fields' transformation:
\begin{equation*}
  \omega_{h_{\beta},B}\underset{h_{\beta}\to 0}{\overset{\mathfrak{P}_{\mathbb{L}(B)}\vee\mathfrak{T}_{\mathbb{L}(B)}}{\longrightarrow}} M_{B}\; \Longrightarrow \; \mathbb{L}\mspace{-1mu}\mathbb{S}_h(\mathfrak{a}) \omega_{h_{\beta},B}\underset{h_{\beta}\to 0}{\overset{\mathfrak{P}_{\mathbb{L}(A)}\vee\mathfrak{T}_{\mathbb{L}(A)}}{\longrightarrow}}\mathbb{L}\mspace{-1mu}\mathbb{S}_0(\mathfrak{a})\,_{*}\, M_{B}\; .
\end{equation*}
The classical transformation is contravariant since it pushes forward measures acting on fields, dual to
test functions (and $\mathbb{L}$ is covariant on test functions). Taking the concrete example of a Lorentz
transformation between Minkowski local spacetimes, this could be interpreted as the fact that the unitary
representation of Lorentz transformations on local quantum bosonic fields converges ``in the Schrödinger
picture'' to the representation of Lorentz transformations on the corresponding local classical fields.

The isotony axiom \eqref{eq:19} is used to provide a notion of \emph{global} observables that is in
accordance with the local structure. The quantum global structure induces a corresponding classical global
structure as well, together with the associated semiclassical properties. Let us define a partial ordering
on the set of local regions of spacetime $\mathrm{Obj}(\mathbf{A})$ (the set of objects of the small
category $\mathbf{A}$): $\preceq$ is defined as follows for any $A,B\in \mathrm{Obj}(\mathbf{A})$
\begin{equation*}
  A \preceq B\; \Longleftrightarrow \;  \,\!^{\exists}\mathfrak{i}_{AB}: A\hookrightarrow B \; .
\end{equation*}
Usually the set of spacetimes is \emph{directed}\footnote{A partially ordered set $(S,\preceq)$ is
  directed iff for any $s,t\in S$, there exists $u\in S$ such that $s\preceq u$ and $t\preceq u$.} by
$\preceq$, and we suppose this to be the case. From \eqref{eq:19} it then follows that
\begin{equation*}
  \Bigl(\mathbb{L}\!\mathbb{W}_h(A,\mathfrak{i}_{AB})\Bigr)_{A\preceq B\in \mathrm{Obj}(\mathbf{A})}
\end{equation*}
is an \emph{inductive} family of C*-algebras (or, using the terminology introduced by \citet[see,
\emph{e.g},][]{haag1992tmp
}, a local net of observables). In addition, both
\begin{equation*}
  \Bigl(\mathbb{L}\mspace{-1mu}\mathbb{S}_h(A,\mathfrak{i}_{AB})\Bigr)_{A\preceq B\in \mathrm{Obj}(\mathbf{A})}\; \textrm{ and }\; \Bigl(\mathbb{L}\mspace{-1mu}\mathbb{S}_0(A,\mathfrak{i}_{AB})\Bigr)_{A\preceq B\in \mathrm{Obj}(\mathbf{A})}
\end{equation*}
are \emph{projective} families, of Banach cones and vector spaces respectively. Since the partial order is
essentially the one given by inclusion (embedding) of one local spacetime into another, it is natural to
define the C*-algebra of quantum global observables as the inductive limit of the net of local
C*-algebras:
\begin{equation*}
  \varinjlim_{A\in \mathrm{Obj}(\mathbf{A})}\mathbb{L}\!\mathbb{W}_h(A)\; .
\end{equation*}
Analogously, is it possible to characterize the global classical states as a projective limit of the local
classical states? The answer is \emph{affirmative}, and it is a consequence of the fact that, since
$(\mathrm{Obj}(\mathbf{A}),\preceq)$ is directed,
\begin{equation*}
  \varprojlim_{A\in \mathrm{Obj}(\mathbf{A})} \Bigl(\mathds{C}^{\mathbb{L}(A)}\Bigr)\cong \Bigl(\,\mathds{C}\,\Bigr)^{\,\underset{A\in \mathrm{Obj}(\mathbf{A})}{\varinjlim} \mathbb{L}(A)} \quad .
\end{equation*}
Let us denote $\underline{X}:=\varinjlim \mathbb{L}(A)$ the \emph{global phase space}. From the above
property, it follows that the projective limit of the spaces of local cylindrical measures is isomorphic
to the space of cylindrical Wigner measures associated to the global phase space $\underline{X}$:
\begin{equation*}
  \varprojlim_{A\in \mathrm{Obj}(\mathbf{A})} \!\!\mathcal{M}_{\mathrm{cyl}}\bigl(\mathbb{L}(A)^{*},\mathbb{L}(A)\bigr)\cong \mathcal{M}_{\mathrm{cyl}}\bigl(\underline{X}^{*},\underline{X}\bigr)\; .
\end{equation*}
On the other hand, by an analogous reasoning, it is also possible to identify the projective limit of the
set of local regular quantum states as states on the global algebra of observables $\varinjlim
\mathbb{L}\!\mathbb{W}_h(A)$. Since in addition we have the convergence of functors $\mathbb{S}_h\underset{h\to 0}{\rightarrowtail}\mathbb{S}_0$ (see
\cref{th:de}), it is clear that global states defined by projective families of convergent semiclassical
states converge semiclassically to the corresponding measure on
$\mathcal{M}_{\mathrm{cyl}}\bigl(\underline{X}^{*},\underline{X}\bigr)$.

\subsubsection{Einstein causality}
\label{sec:semicl-einst-caus}

As we have discussed in the previous section, Einstein causality can be abstractly formulated as the fact
that the local functor preserves tensor structures \eqref{eq:21}. As a consequence, also the local bosonic
quantization functor and the functor of local quantum states preserve tensor structures
(\cref{lemma:2}). At the classical level, there is a tensor structure for cylindrical measures
$\mathbf{CylM}$, given by the product measures: for any $\mathcal{A}\subseteq \mathds{R}^A$,
$\mathcal{B}\subseteq \mathds{R}^B $, $A,B\in \mathbf{Set}$, then
\begin{equation}
  \label{eq:25}
  \mathcal{M}_{\mathrm{cyl}}(A,\mathcal{A})\otimes \mathcal{M}_{\mathrm{cyl}}(B,\mathcal{B}):= \mathcal{M}_{\mathrm{cyl}}(A\times B,\mathcal{A}\oplus \mathcal{B})\; .
\end{equation}
In addition, let
\begin{align*}
  \mspace{193mu}&\mathfrak{a}:A_1\to A_2\quad ,\quad \mathfrak{b}:B_1\to B_2 &(A_1,A_2,B_1,B_2\in \mathbf{Set})
\end{align*}
be compatible with
$\mathcal{M}_{\mathrm{cyl}}(A_1,\mathcal{A}_1),\mathcal{M}_{\mathrm{cyl}}(A_2,\mathcal{A}_2)$ and
$\mathcal{M}_{\mathrm{cyl}}(B_1,\mathcal{B}_1),\mathcal{M}_{\mathrm{cyl}}(B_2,\mathcal{B}_2)$
respectively. Then
\begin{equation*}
  \mathfrak{a}\otimes \mathfrak{b}:= \mathfrak{a}\times \mathfrak{b}\quad ,\quad \mathfrak{a}\times \mathfrak{b}: A_1\times B_1\to A_2\times B_2
\end{equation*}
is compatible with $\mathcal{M}_{\mathrm{cyl}}(A\times B,\mathcal{A}\oplus \mathcal{B})$. The identity is
$i_{\mathbf{CylM}}=\mathcal{M}_{\mathrm{cyl}}\bigl(\{\emptyset\},\{\mathbf{0}\}\bigl)$, where
$\mathbf{0}(\varnothing)=0$ is the zero function. Using this structure, $\mathbf{CylM}^{\otimes }$ is a
monoidal category. The classical counterpart of \cref{lemma:2} is then the following lemma.
\begin{lemma}
  \label{lemma:3}
  If \eqref{eq:21} holds, then $\mathbb{L}\mspace{-1mu}\mathbb{S}_0^{\otimes}$ is homomorphic.
\end{lemma}
The convergence of functors \eqref{eq:18} can then be rewritten as a convergence of homomorphic functors
$\mathbb{L}\mspace{-1mu}\mathbb{S}_h^{\otimes_{\alpha}}\underset{h\to
  0}{\rightarrowtail}\mathbb{L}\mspace{-1mu}\mathbb{S}_0^{\otimes }$, for any $\alpha_{\mathrm{min}}\leq \alpha\leq
\alpha_{\mathrm{max}}$. One important consequence of this convergence is the fact that spacelike
entanglement of local regions can only be destroyed in the classical limit. In other words, given a
quantum state describing the fields localized in two spacelike separated regions of spacetime, the
corresponding entanglement could only disappear in the classical limit $h\to 0$. The most convenient way
to prove that is to show that on one hand, given any semiclassical quantum state with no entanglement, all
its corresponding cylindrical Wigner measures are statistically independent as well; on the other hand,
that there exist entangled quantum states whose classical limit is not entangled. Let us remark that the
proof of what we have just discussed (that is given after the proposition below), could be immediately
adapted to prove the more general \cref{prop:2}.

\begin{prop}
  \label{prop:5}
  If \eqref{eq:21} holds, then
  \begin{equation*}
    \mathbb{L}\mspace{-1mu}\mathbb{S}_h^{\otimes_{\alpha}}\underset{h\to 0}{\rightarrowtail}\mathbb{L}\mspace{-1mu}\mathbb{S}_0^{\otimes
    }\; ,
  \end{equation*}
  for any $\alpha_{\mathrm{min}}\leq \alpha\leq \alpha_{\mathrm{max}}$, with both
  $\mathbb{L}\mspace{-1mu}\mathbb{S}_h^{\otimes_{\alpha}}$ and $\mathbb{L}\mspace{-1mu}\mathbb{S}_0^{\otimes}$ homomorphic. In
  addition, entanglement between spacelike separated regions can only be destroyed by the classical limit
  $h\to 0$.
\end{prop}
\begin{proof}
  The convergence of functors has already been discussed, let us prove the fact that entanglement can only
  be destroyed in the classical limit. Let us consider two (spacelike separated) local regions of
  spacetime $A,B\in \mathbf{A}^{\otimes }$. Then
  \begin{equation*}
    \mathbb{L}\!\mathbb{W}_h(A\otimes B)=\mathbb{W}_h\bigl(\mathbb{L}(A)\oplus \mathbb{L}(B)\bigr) \; .
  \end{equation*}
  Now, let $\omega_{h,A}$ be a semiclassical quantum state on $\mathbb{W}_h\bigl(\mathbb{L}(A)\bigr)$,
  $\omega_{h,B}$ a semiclassical quantum state on $\mathbb{W}_h\bigl(\mathbb{L}(B)\bigr)$, and
  $\omega_{h,A\otimes B}$ an \emph{entangled} semiclassical quantum state on
  $\mathbb{W}_h\bigl(\mathbb{L}(A)\oplus \mathbb{L}(B)\bigr)$, \emph{i.e.}\ the latter cannot be written
  as the tensor product of one state acting on $\mathbb{W}_h\bigl(\mathbb{L}(A)\bigr)$ and another acting
  on $\mathbb{W}_h\bigl(\mathbb{L}(B)\bigr)$. Then
  \begin{equation*}
    \omega_{h,A}\otimes \omega_{h,B} \in \mathbb{W}_h\bigl(\mathbb{L}(A)\oplus \mathbb{L}(B)\bigr)'
  \end{equation*}
  is a \emph{non-entangled} semiclassical quantum state. Let us consider one of its cluster points in
  either the $\mathfrak{P}_{\mathbb{L}(A)\oplus \mathbb{L}(B)}$ or the $\mathfrak{P}_{\mathbb{L}(A)\oplus
    \mathbb{L}(B)}\vee\mathfrak{T}_{\mathbb{L}(A)\oplus \mathbb{L}(B)}$ topology:
  \begin{equation*}
    M_{AB}\in \mathcal{M}_{\mathrm{cyl}}\Bigl(\mathbb{L}(A)^{*}\times \mathbb{L}(B)^{*},\mathbb{L}(A)\oplus  \mathbb{L}(B) \Bigr)\; .
  \end{equation*}
  The cluster point can indeed be written as a cylindrical measure on
  $\mathbb{L}(A)^{*}_{\mathbb{L}(A)}\times \mathbb{L}(B)^{*}_{\mathbb{L}(B)}$, \emph{i.e.}\ on the set
  $\mathbb{L}(A)^{*}\times \mathbb{L}(B)^{*}$ with the weak
  $\sigma\bigl(\mathbb{L}(A)^{*},\mathbb{L}(A)\bigr)\times
  \sigma\bigl(\mathbb{L}(B)^{*},\mathbb{L}(B)\bigr)$ topology, since
  \begin{equation*}
    \sigma\bigl(\mathbb{L}(A)^{*}\times \mathbb{L}(B)^{*},\mathbb{L}(A)\oplus \mathbb{L}(B) \bigr)=\sigma\bigl(\mathbb{L}(A)^{*},\mathbb{L}(A)\bigr)\times \sigma\bigl(\mathbb{L}(B)^{*},\mathbb{L}(B)\bigr)\;.
  \end{equation*}
  Now, from the fact that $\omega_{h,A}\otimes \omega_{h,B}$ is a tensor product it follows that
  \begin{equation*}
    M_{AB}=M_A\otimes M_B\quad ,\quad M_A\in \mathcal{M}_{\mathrm{cyl}}\bigl(\mathbb{L}(A)^{*},\mathbb{L}(A)\bigr)\; ,\; M_B\in \mathcal{M}_{\mathrm{cyl}}\bigl(\mathbb{L}(B)^{*},\mathbb{L}(B)\bigr)\; .
  \end{equation*}
  This is due to a general property of tensor products, however let us prove the result explicitly for the
  $\mathfrak{T}_{\mathbb{L}(A)\oplus \mathbb{L}(B)}$ convergence. Let
  $\mathcal{G}_{\omega_{h_{\beta},A}\otimes \omega_{h_{\beta},B}}(\cdot )$ be the family of generating
  functionals converging pointwise to $\hat{M}_{AB}(\cdot )$. Since the quantum state is a tensor product,
  it follows that for any $a_1,a_2\in \mathbb{L}(A)$ and $b_1,b_2\in \mathbb{L}(B)$,
  \begin{equation*}
    \begin{split}
      \mathcal{G}_{\omega_{h_{\beta},A}\otimes \omega_{h_{\beta},B}}(a_1,b_1)\mathcal{G}_{\omega_{h_{\beta},A}\otimes \omega_{h_{\beta},B}}(a_2,b_2)=\\\mathcal{G}_{\omega_{h_{\beta},A}\otimes \omega_{h_{\beta},B}}(a_1,b_2)\mathcal{G}_{\omega_{h_{\beta},A}\otimes \omega_{h_{\beta},B}}(a_2,b_1)
    \end{split}
  \end{equation*}
  Then, in the limit $h_{\beta}\to 0$,
  \begin{equation*}
    \hat{M}_{AB}(a_1,b_1)\hat{M}_{AB}(a_2,b_2)=\hat{M}_{AB}(a_1,b_2)\hat{M}_{AB}(a_2,b_1)\; .
  \end{equation*}
  Hence there exist two functions of positive type (continuous on finite dimensional subspaces)
  $\hat{M}_A:\mathbb{L}(A)\to \mathds{C}$ and $\hat{M}_B:\mathbb{L}(B)\to \mathds{C}$ such that
  \begin{equation*}
    \hat{M}_{AB}=\hat{M}_{A}\otimes \hat{M}_{B}\; .
  \end{equation*}
  Hence $\hat{M}_{AB}$ is the Fourier transform of the tensor product of two cylindrical measures.

  It remains to show that there exist entangled quantum states whose classical limit is not entangled. The
  easiest way to do it is ``perturbing'' a convergent non-entangled semiclassical state. Let
  $\omega_{h_{\beta},A}\otimes \omega_{h_{\beta},B}$ be the non-entangled semiclassical state converging
  in the $\mathfrak{P}_{\mathbb{L}(A)\oplus \mathbb{L}(B)}\vee\mathfrak{T}_{\mathbb{L}(A)\oplus
    \mathbb{L}(B)}$ topology to $M_A\otimes M_B$. Then let us consider the semiclassical quantum state
  \begin{equation*}
    \varpi_{h_{\beta}}:=\omega_{h_{\beta},A}\otimes \omega_{h_{\beta},B} + h_{\beta}\omega_{h_{\beta},A\otimes B}\; ;
  \end{equation*}
  where $\omega_{h_{\beta},A\otimes B}$ is the entangled state defined above at the beginning of this
  proof. Since $\omega_{h_{\beta},A\otimes B}$ is entangled, also $\varpi_{h_{\beta}}$ is entangled. In
  addition, it converges in the $\mathfrak{P}_{\mathbb{L}(A)\oplus
    \mathbb{L}(B)}\vee\mathfrak{T}_{\mathbb{L}(A)\oplus \mathbb{L}(B)}$ topology to $M_A\otimes M_B$
  (since the non-entangled perturbation becomes small as $h_{\beta}\to 0$).
\end{proof}

\subsubsection{Time slice and Egorov type theorems}
\label{sec:semicl-time-slice}

The time slice axiom \eqref{eq:24} defines a non-empty algebra of bosonic observables on any Cauchy
surface, and states that the knowledge of such algebra is sufficient to describe all local observables of
regions containing the aforementioned Cauchy surface. A Cauchy surface is a ``slice of time''. Therefore,
if the time slice axiom holds true, one should define the (Heisenberg-picture) evolution of the system as
a map between the observables of two different Cauchy surfaces. Let $\Sigma,\Theta\in
\mathrm{Cau}_{\mathbf{A}}$ such that $^{\exists }A\supset \Sigma\cup \Theta$. Then \eqref{eq:24} implies
that a $\Theta$-time generator $W_h(x_{\Theta})\in \mathbb{W}_h(X_{\Theta},\varsigma_{\Theta})$ is evolved
to some $\Sigma$-observable $O_h(x_{\Theta})=U^{(\Sigma\Theta)}_h\bigl(W_h(x_{\Theta})\bigr) \in
\mathfrak{O}(\Sigma)$ by the *-isomorphism $U^{(\Sigma\Theta)}_h$ defined by
\begin{equation*}
  U^{(\Sigma\Theta)}_h:= P_{A,\Sigma}^{-1} \circ P_{A,\Theta}\; .
\end{equation*}
For interacting theories \emph{we cannot expect} the latter to be a $\Sigma$-time-zero field, and it could
even not belong to the Weyl C*-algebra \citep[for interesting examples, in non-relativistic particle
quantum mechanics, see][]{fannes1974cmp
}. By duality, it also follows that
\begin{equation*}
  ^{\mathrm{t}}U_h^{(\Sigma\Theta)}: \mathfrak{O}(\Sigma)'\to \mathfrak{O}(\Theta)'
\end{equation*}
transforms states of $\mathbb{W}_h(X_{\Sigma},\varsigma_{\Sigma})$ in complex states of
$\mathbb{W}_h(X_{\Theta},\varsigma_{\Theta})$. Physically, we may think of
$^{\mathrm{t}}U_h^{(\Sigma\Theta)}$ as the \emph{evolution in the Schrödinger picture} (not necessarily
positivity preserving), mapping a quantum state at a given time (described by the Cauchy surface
$\Sigma$), to a state at another time (described by the Cauchy surface $\Theta$).

Consider now the semiclassical complex quantum state\footnote{A complex quantum state can always be
  decomposed in four positive states. A complex state is semiclassical if each of the four positive states
  of its decomposition are semiclassical.}
\begin{equation*}
  \omega_h^{(\Sigma)}\in \mathfrak{O}(\Sigma)' \overset{^{\mathrm{t}}\mathrm{w}_{\Sigma}}{\to} \mathbb{W}_h(X_{\Sigma},\varsigma_{\Sigma})'\; ,
\end{equation*}
and suppose that
\begin{equation*}
  ^{\mathrm{t}}U_h^{(\Sigma\Theta)}\omega_h^{(\Sigma)}\in\mathfrak{O}(\Theta)'\overset{^{\mathrm{t}}\mathrm{w}_{\Theta}}{\to} \mathbb{W}_h(X_{\Theta},\varsigma_{\Theta})'
\end{equation*}
is semiclassical as well. In addition, suppose that both $\omega_h^{(\Sigma)}$ and
$^{\mathrm{t}}U_h^{(\Sigma\Theta)}\omega_h^{(\Sigma)}$ converge\footnote{The convergence of a
  semiclassical complex state to a complex cylindrical measure has to be intended as follows: all the four
  positive states of its decomposition converge separately to positive cylindrical measures (whose
  combination is a complex cylindrical measure).}, as $h\to 0$, in the $\mathfrak{P}$ or
$\mathfrak{P}\vee\mathfrak{T}$ topologies, to $M_{\Sigma}\in \mathcal{M}_{\mathrm{cyl}}\bigl(\,\!^{\exists
}S_{\Sigma},\,\!  ^{\exists }e_{X_{\Sigma}}(X_{\Sigma}) \bigr)_{\mathds{C}}$ and $M_{\Theta}\in
\mathcal{M}_{\mathrm{cyl}}\bigl(\,\!^{\exists }S_{\Theta},\,\!  ^{\exists }f_{X_{\Theta}}(X_{\Theta})
\bigr)_{\mathds{C}}$ respectively. Hence it would be desirable to have a \emph{compatible} map
\begin{equation*}
  \Phi_{\Theta\Sigma}: S_{\Theta}\leftarrow S_{\Sigma}
\end{equation*}
such that
\begin{equation*}
  M_{\Theta}=\Phi_{\Theta\Sigma}\,_{*}\, M_{\Sigma}\; .
\end{equation*}
The map $\Phi_{\Theta\Sigma}$ could then interpreted as the \emph{classical flow}, evolving the fixed-time
classical fields from one Cauchy surface to the other. In particular, let us call an \emph{Egorov theorem}
the following sentence:
\begin{gather*}
  \bigl(\,\!^{\exists }\Phi_{\Theta\Sigma}:S_{\Theta}\leftarrow S_{\Sigma}\bigr)_{(\Sigma,\Theta)\in \mathrm{Cau}_{\mathbf{A}}^2}\textrm{ , with } \Phi_{\Theta\Sigma} \textrm{ compatible and } \Phi_{\Sigma\Sigma}= \mathrm{id}\textrm{ , such that}\\ \Bigl(\omega_h^{(\Sigma)}\underset{h\to 0}{\overset{(\mathfrak{P}\vee\mathfrak{T})_{\Sigma}}{\longrightarrow}} M_{\Sigma} \;\Longleftrightarrow \forall(\Sigma,\Theta)\in \mathrm{Cau}_{\mathbf{A}}^2\;,\; ^{\mathrm{t}}U_h^{(\Sigma\Theta)}\omega_{h}^{(\Sigma)}\underset{h\to 0}{\overset{(\mathfrak{P}\vee\mathfrak{T})_{\Theta}}{\longrightarrow}}\Phi_{\Theta\Sigma}\, _{*}\, M_{\Sigma} \Bigr)\; .
\end{gather*}
The terminology is borrowed from standard semiclassical analysis, although the Egorov theorem for finite
dimensional systems holds in a stronger form. In some semi-relativi\-stic and non-relativistic quantum
field theories, Egorov theorems have been proved rigorously \citep[see
e.g.][]{ammari2008ahp
  ,ammari2014jsp
  ,ammari2015asns
  ,ammari2016cms
  ,ammari2016barxiv
  ,ammari2017sima
}. They are proved to hold for a suitable set of Fock-normal families of states whose cylindrical Wigner
measures are concentrated as Radon measures (on a space $\mathcal{Z}$ with some structure, usually a
Hilbert space) \citep[see][for \emph{a priori} conditions that guarantee
concentration]{ammari2008ahp
  ,falconi2017ccm
}. The classical flow is the same for all families of states, and it is explicitly identified to be the
Hamiltonian flow $\Phi_{t-t_0,\mathcal{Z}}$ corresponding to a suitable classical nonlinear PDE, globally
well-posed on $\mathcal{Z}$. Interestingly, in
\citep{ammari2017sima
} the Egorov theorem holds for a \emph{renormalized} quantum system (the Nelson model) and the
corresponding naïf classical flow (Schrödinger-Klein-Gordon system with Yukawa coupling), so at least in
that case the renormalization does not affect the classical limit. On specific quantum states, in
particular on coherent states, it is possible to prove Egorov-type results more easily exploiting the
special semiclassical structure of such configurations \citep[see
\emph{e.g.}][]{hepp1974cmp
  ,eckmann1975lmp
  ,ginibre1979cmpI
  ,ginibre1979cmp2
  ,donald1981cmp
  ,ginibre2006ahp
}. Let us conclude remarking that it does not seem possible to infer Egorov theorems from the
\eqref{eq:24} axiom in a systematic fashion: they have to be proved on a case-by-case basis.

\subsubsection{Invariant, KMS, and Ground States}
\label{sec:semicl-kms-ground}

The notions of Invariant, KMS, and Ground States play an important role in statistical and relativistic
quantum mechanics (in Minkwoski spacetime). In order to introduce them we should define the concept of
quantum symmetries and quantum dynamical systems. We will not do it in full generality, but we will
consider only symmetries of quantum fields. Let $(X,\varsigma)\in \mathbf{Symp}_{\mathds{R}}$ be a phase
space, and let $(G,\cdot )$ be a (Lie) group that is represented by\footnote{As it is customary, we are
  interested in automorphic representations. Therefore, in the category of symplectic spaces, the
  representation is linear and symplectic.}  $\bigl(s(g)\bigr)_{g\in G}$ on $(X,\varsigma)$. This induces
a representation of the group on quantum fields (by action on generators, see~\cref{sec:maps-conv-prod}),
that is extended uniquely to a group of *-automorphisms $\bigl(\mathfrak{s}_h(g)\bigr)_{g\in
  G}=\Bigl(\mathbb{W}_h\bigl(s(g)\bigr) \Bigr)_{g\in G}$ on $\mathbb{W}_h(X,\varsigma)$. Clearly, we have
that
\begin{align*}
  \mspace{90mu}&\mathfrak{s}_h(g_1) \mathfrak{s}_h(g_2)^{-1} W_h(x) = W_h\bigl(s(g_1\cdot g_2^{-1})x\bigr) &\bigl(\forall x\in X,\forall g_1,g_2\in G \bigr)\; .
\end{align*}

\begin{definition} [Symmetry group of quantum fields]
  \label{def:7}
  Let $(X,\varsigma)\in \mathbf{Symp}_{\mathds{R}}$ ; and let $\bigl(\mathfrak{s}_h(g)\bigr)_{g\in G}$ be
  a representation of the group $G$ on $\mathbb{W}_h(X,\varsigma)$. Then
  $\bigl(\mathfrak{s}_h(g)\bigr)_{g\in G}$ is a \emph{group of symmetry transformations of quantum fields}
  iff $\bigr(\,\!^{\exists }s(g)\bigr)_{g\in G}$ representation of $G$ on $(X,\varsigma)$ such that
  \begin{align*}
    \mspace{170mu}&\mathfrak{s}_h(g)=\mathbb{W}_h\bigl(s(g)\bigr) &(\forall g\in G)\; .
  \end{align*}
\end{definition}

Symmetries of quantum fields play an important role in both nonrelativistic and relativistic quantum field
theories. In particular, it is one of the Wightman axioms that the proper, orthochronous Poincaré group is
a symmetry group of quantum fields. It is another axiom that the theory is set in the GNS representation
of an \emph{invariant state} of the aforementioned Poincaré group. The definition of an invariant state is
rather intuitive, and it is the following.

\begin{definition} [Invariant State]
  \label{def:8}
  Let $(X,\varsigma)\in \mathbf{Symp}_{\mathds{R}}$ ; and let $\bigl(\mathfrak{s}_h(g)\bigr)_{g\in G}$ be
  a representation of the group $G$ on $\mathbb{W}_h(X,\varsigma)$. Then a quantum state $\omega_h\in
  \mathbb{S}_h(X,\varsigma)$ is \emph{$G$-invariant} iff
  \begin{align*}
    \mspace{150mu}&^{\mathrm{t}}\mathfrak{s}_h(g) \omega_h:= \omega_h\circ \mathfrak{s}_h(g)=\omega_h &(\forall g\in G)\;.
  \end{align*}
\end{definition}
Analogously, it is possible to give a definition of invariant classical states, \emph{i.e.}\ invariant
cylindrical measures associated to a phase space.
\begin{definition} [Invariant measure]
  \label{def:9}
  Let $G$ be a group, represented by $\bigl(s(g)\bigr)_{g\in G}$ on $(X,\varsigma)\in
  \mathbf{Symp}_{\mathds{R}}$. Then a measure $M\in \mathcal{M}_{\mathrm{cyl}}(X^{*},X)$ is
  \emph{$G$-invariant} iff
  \begin{align*}
    \mspace{145mu}&^{\mathrm{t}}s(g)\,_{*}\, M = M &(\forall g\in G)\; .
  \end{align*}
\end{definition}
\begin{lemma}
  \label{lemma:4}
  A measure $M\in \mathcal{M}_{\mathrm{cyl}}(X^{*},X)$ is $G$-invariant iff
  \begin{align}
    \label{eq:26}
    \mspace{187mu}&\hat{M}\bigl(s(g)x\bigr)=\hat{M}(x) &(\forall x\in X,\forall g\in G)\; .
  \end{align}
  Using \cref{thm:2,eq:26} it is then possible to define invariant measures for any space of cylindrical
  measures $\mathcal{M}_{\mathrm{cyl}}\bigl(\,\!^{\exists }A,\,\!^{\exists }e_X(X) \bigr)$ associated to
  $X$.
\end{lemma}
\begin{proof}
  The lemma is an easy consequence of Bochner's theorem, \cref{thm:2}.
\end{proof}

From \crefrange{def:8}{def:9}, it is clear that we can apply \cref{th:de} to nicely characterize the
cylindrical Wigner measures of semiclassical invariant states with no loss of mass.

\begin{proposition}
  \label{prop:4}
  Let $\bigl(\mathfrak{s}_h(g)\bigr)_{g\in G}$ be a group of symmetry transformations of quantum fields on
  $\mathbb{W}_h(X,\varsigma)$, and $\bigl(s(g)\bigr)_{g\in G}$ the corresponding representation of $G$ on
  $(X,\varsigma)$. Then \emph{the cylindrical Wigner measure} of any $G$-invariant semiclassical quantum
  state with no loss of mass \emph{is $G$-invariant}.
\end{proposition}
\begin{proof}
  By \cref{eq:8}, it follows that
  \begin{align*}
    \omega_{h_{\beta}}\underset{h_{\beta}\to 0}{\overset{\mathfrak{P}\vee\mathfrak{T}}{\longrightarrow}}M \quad \Longrightarrow \quad (\forall g\in G)\; ^{\mathrm{t}}\mathfrak{s}_{h_{\beta}}(g)\omega_{h_{\beta}}\underset{h_{\beta}\to 0}{\overset{\mathfrak{P}\vee\mathfrak{T}}{\longrightarrow}} \,\!^{\mathrm{t}}s(g) \,_{*}\, M\; .
  \end{align*}
  Hence from $^{\mathrm{t}}\mathfrak{s}_{h_{\beta}}(g)\omega_{h_{\beta}}=\omega_{h_{\beta}}$ it follows
  that $^{\mathrm{t}}s(g) \,_{*}\, M=M$.
\end{proof}
In addition, by \cref{thm:5} we know that any $G$-invariant measure is the limit of at least one
semiclassical quantum state, but such quantum state \emph{may not be $G$-invariant}, and there may be
invariant classical measures that are not the limit of any invariant semiclassical state (in fact, there
are systems in which no invariant quantum state exists). For extremal (\emph{i.e.}\ pure) invariant
states, the so-called Haag's theorem \citep[see, \emph{e.g.},][Corollaries
5.3.41-42]{bratteli1997tmp2
} holds, and it plays a very important role in relativistic theories. In fact, one consequence of Haag's
theorem is that free and interacting theories must correspond to inequivalent representations of the
canonical commutation relations.
\begin{thm} [Haag]
  \label{thm:6}
  Let $\omega_h,\varpi_h\in \mathbb{W}_h(X,\varsigma)'_+$ be two $G$-Abelian pure normalized
  states\footnote{The $G$-abelian states are a subset of the set of $G$-invariant states that satisfy an
    additional property of commutativity, let us omit the precise definition here \citep[see,
    \emph{e.g.},][Definition
    4.3.6]{bratteli1987tmp1
    }.}. Then they are either equal or disjoint.
\end{thm}
At the classical level, however, it may happen that two different (therefore disjoint) $G$-Abelian pure
normalized regular (therefore semiclassical) states converge to the same cylindrical Wigner measure. It
would be interesting to find explicit examples in which it is the case, or examples in which they converge
to two mutually singular Wigner measures.

A special set of invariant states that has been extensively studied in statistical mechanics is the subset
of so-called \emph{KMS states}.
\begin{definition}[KMS state]
  \label{def:10}
  Let $\bigl(\mathfrak{s}_h(t)\bigr)_{t\in \mathds{R}}$ be a representation of the abelian group
  $\mathds{R}$ on $\mathfrak{B}_h\overset{^{\exists }\mathrm{w}_X}{\hookleftarrow}
  \mathbb{W}_h(\,\!^{\exists }X,\,\!^{\exists }\varsigma)$. Then a quantum state $\omega_h\in
  \mathfrak{B}'_h$ is a $\mathfrak{s}_h$-\emph{KMS} state at value $\beta\in \mathds{R}$ iff
  \begin{align}
    \label{eq:28}
    \mspace{175mu}&\omega_h(b_ha_h)= \omega_h\bigl(a_h \mathfrak{s}_h(ih\beta)b_h\bigr)&(\forall a_h,b_h\in \mathfrak{B}_h)\;.
  \end{align}
\end{definition}
\begin{remark}
  \label{rem:5}
  It is sufficient to test the KMS condition on a norm dense and $\mathfrak{s}_h$-invariant subalgebra.
\end{remark}
A KMS state is a state that is almost a trace (it is a trace only if $\beta=0$), the deviation being
measured by $\mathfrak{s}_h$. In this context, the abelian group $\mathds{R}$ is usually interpreted as
the group of time translations, and thus $\mathfrak{s}_h$ is the quantum dynamical map (and $\beta$ is the
inverse temperature). Semiclassically, one would like to prove that the cylindrical Wigner measures of
semiclassical quantum states satisfy an equation of the following type:
\begin{equation}
  \label{eq:27}
  \int_{}^{}\bigl\{a(z),b(z)\bigr\}  \mathrm{d}M(z) = \beta \int_{}^{}b(z)\bigl\{a(z),\mathfrak{h}(z)\bigr\}  \mathrm{d}M(z)\;,
\end{equation}
for any $a,b$ in a suitable set of classical observables, where $\mathfrak{h}$ is the classical
Hamiltonian observable, and $\bigl\{\cdot,\cdot \bigr\}$ is a Poisson bracket. This ``static''
semiclassical KMS condition has been studied for systems with finitely many degrees of freedom
\citep{gallavotti1975nc
  ,genovese2012jsp
}, but its origin from the quantum KMS condition was justified only formally. It is possible to derive
\cref{eq:27} from \cref{eq:28} in our framework; however, unless $X$ is finite dimensional, additional
properties should hold, and they have to be proved case-by-case. For any $x, y\in X$, $\varsigma(x,\cdot
)\in X^{*}$ and $\varsigma(x,\cdot )\neq \varsigma(y,\cdot )$ if $x\neq y$ (by the non-degeneracy of
$\varsigma$), and thus $X\overset{\mathrm{s}}{\hookrightarrow} X^{*}_X$. In addition, $\varsigma(\cdot
,\cdot )$ extends to a symplectic form $\tilde{\varsigma}(\cdot ,\cdot )$ on $\mathrm{s}(X)$ by
$\tilde{\varsigma}(\cdot ,\cdot )=\varsigma(\mathrm{s}^{-1}\,\cdot \,,\mathrm{s}^{-1}\,\cdot \,)$. The map
$\mathrm{s}$ is a bijection iff $X$ is finite dimensional. Since $X^{*}_X$ is locally convex, there is a
notion of smooth maps and therefore a Poisson bracket $\bigl\{\cdot ,\cdot \bigr\}$ can be defined on
$\bigl(\mathrm{s}(X),\tilde{\varsigma}\bigr)$. These two types of results can therefore be expected to be
provable in suitable systems:
\begin{itemize}
\item If the generator\footnote{The generator $\delta_h$ of a (strongly continuous) group
    $\bigl(\mathfrak{s}_h(t)\bigr)_{t\in \mathds{R}}$ (with suitable continuity properties) is a map from
    a dense domain $D(\delta_h)\subset \mathfrak{B}_h$ to $\mathfrak{B}_h$ such that for any $a_h\in
    \mathfrak{B}_h$, $\mathfrak{s}_h(t)a_h$ is differentiable with respect to $t$, and $\delta_h
    a_h=\frac{\mathrm{d}}{\mathrm{d}t}\mathfrak{s}_h(t)a_h\bigr\rvert_{t=0}$.}  of
  $\bigl(\mathfrak{s}_h(t)\bigr)_{t\in \mathds{R}}$ is $\delta_h=\frac{i}{h}[\mathfrak{h}_h,\,\cdot\, ]$,
  with $\mathfrak{h}_h= \mathrm{Op}_{\frac{1}{2}}^h(\mathfrak{h})$, $\mathfrak{h}:X^{*}_X\to \mathds{R}$ a
  cylindrical function with base $\Phi$, and $\omega_h$ is a semiclassical quantum state such that
  $\omega_h(\mathfrak{h}_h)\leq C$ uniformly with respect to $h$; then
  \begin{equation*}
    \begin{tikzcd}
      \omega_{h}\underset{h\to
        0}{\overset{\mathfrak{P}}{\longrightarrow}} M \arrow[d, Rightarrow]\\
      \begin{aligned}
        \omega_{h}\Bigl(b_{h}a_{h}\Bigr)= \omega_h\Bigl(a_{h} \mathfrak{s}_h(ih\beta)b_h\Bigr)
        \underset{h\to 0}{\overset{\mathfrak{P}}{\longrightarrow}}
        \int_{X^{*}_X/\Phi}^{}\bigl\{a_{\Phi}(z),b_{\Phi}(z)\bigr\} \mathrm{d}\mu_{\Phi}(z) \\= \beta
        \int_{X^{*}_X/\Phi}^{}b_{\Phi}(z)\bigl\{a_{\Phi}(z),\mathfrak{h}_{\Phi}(z)\bigr\}
        \mathrm{d}\mu_{\Phi}(z)
      \end{aligned}
    \end{tikzcd}
  \end{equation*}
  for any $a_h= \mathrm{Op}^h_{\frac{1}{2}}(a)$ and $b_h= \mathrm{Op}^h_{\frac{1}{2}}(b)$, with
  $a,b:X^{*}_X\to \mathds{R}$ smooth cylindrical functions with the same base $\Phi$ as $\mathfrak{h}$.
  
\item If the generator of $\bigl(\mathfrak{s}_h(t)\bigr)_{t\in \mathds{R}}$ is
  $\delta_h=\frac{i}{h}[\mathfrak{h}_h,\,\cdot\, ]$, with $\mathfrak{h}_h$ a suitable quantization of some
  symbol $\mathfrak{h}:\mathrm{s}(X)\to \mathds{R}$ (possibly only densely defined with domain
  $D(\mathfrak{h})\subseteq \mathrm{s}(X)$ carrying a suitable topology), and $\omega_h$ is a
  semiclassical quantum state\footnote{Let us denote by
    $(\mathscr{H}_{\omega_h},\pi_{\omega_h},\Omega_{\omega_h})$ the GNS representation given by
    $\omega_h$.} with no loss of mass such that $\bigl\langle \Omega_{\omega_h} ,
  \pi_{\omega_h}(\mathfrak{h}_h)\Omega_{\omega_h} \bigr\rangle_{\mathscr{H}_{\omega_h}}\leq C$ uniformly
  with respect to $h$ (more regularity of the state could be necessary); then
  \begin{equation*}
    \begin{tikzcd}
      \omega_{h}\underset{h\to
        0}{\overset{\mathfrak{P}\vee \mathfrak{T}}{\longrightarrow}} \mu\in \mathcal{M}_{\mathrm{rad}}\bigl(D(\mathfrak{h})\bigr) \arrow[d, Rightarrow]\\
      \begin{aligned}
        \omega_{h}\Bigl(b_{h}a_{h}\Bigr)= \omega_h\Bigl(a_{h} \mathfrak{s}_h(ih\beta)b_h\Bigr)
        \underset{h\to 0}{\overset{\mathfrak{P}\vee \mathfrak{T}}{\longrightarrow}}
        \int_{D(\mathfrak{h})}^{}\bigl\{a(z),b(z)\bigr\} \mathrm{d}\mu(z) \\= \beta
        \int_{D(\mathfrak{h})}^{}b(z)\bigl\{a(z),\mathfrak{h}(z)\bigr\} \mathrm{d}\mu(z)
      \end{aligned}
    \end{tikzcd}
  \end{equation*}
  for any $a_h,b_h$ that are suitable quantizations of smooth symbols $a,b:\mathrm{s}(X)\to \mathds{R}$.
\end{itemize}

A special subset of the set of KMS states is the set of ground states. There are general algebraic
definitions of ground states as KMS states with special properties, let us focus however on a more
physical definition. Let $(\mathfrak{B}_h,\mathfrak{s}_h)$, $\mathfrak{B}_h\overset{^{\exists
  }\mathrm{w}_X}{\hookleftarrow} \mathbb{W}_h(\,\!^{\exists }X,\,\!^{\exists }\varsigma)$, be an algebra
of bosonic observables with evolution group $\mathfrak{s}_h$ whose generator is $\delta_h$.

\begin{definition}[Ground state]
  \label{def:11}
  A state $\omega_h\in \mathfrak{B}_h'$ is a \emph{ground state} iff
  $\pi_{\omega_h}(\delta_h)=\frac{i}{h}[\,\!^{\exists }\mathfrak{h}_h,\,\cdot \,]$, with the Hamiltonian
  $\mathfrak{h}_h$ self-adjoint and bounded from below on $\mathscr{H}_{\omega_h}$, and such that
  \begin{equation*}
    \mathfrak{h}_h\Omega_{\omega_h}=\lambda_0\Omega_{\omega_h}\; ,\; \lambda_0=\min\bigl\{\lambda\in \mathrm{spec}(\mathfrak{h}_h)\bigr\}\; .
  \end{equation*}
\end{definition}
Taking into account the semiclassical properties of Invariant and KMS states, it is natural to expect that
semiclassical quantum ground states would converge to classical ground states. Suppose that $\omega_h$ is
a semiclassical quantum ground state, and that $\mathfrak{h}_h$ is the quantization of some bounded from
below and densely defined symbol $\mathfrak{h}:\mathrm{s}(X)\to \mathds{R}$ (the classical energy, with
domain $D(\mathfrak{h})\subseteq \mathrm{s}(X)$). Then it is interesting to prove whether the following
statement is true:
\begin{equation}
  \label{eq:29}
  \omega_h\underset{h\to 0}{\overset{\mathfrak{P}\vee \mathfrak{T}}{\longrightarrow}}\mu\in \mathcal{M}_{\mathrm{rad}}\bigl(D(\mathfrak{h})\bigr)\; \Longrightarrow \; \int_{D(\mathfrak{h})}^{}\mathfrak{h}(z)  \mathrm{d}\mu(z)=\inf_{z\in D(\mathfrak{h})}\mathfrak{h}(z)\; .
\end{equation}
Suppose that \eqref{eq:29} holds, and that $\mathfrak{h}$ has minimizers, \emph{i.e.}\ suppose that
$\varnothing \neq \,\!^{\exists }\mathrm{Min}(\mathfrak{h}) \subset D(\mathfrak{h})$ such that
\begin{align*}
  \mspace{223mu}&\mathfrak{h}(z_0)=\inf_{z\in D(\mathfrak{h})}\mathfrak{h}(z) &\bigl(\forall z_0\in \mathrm{Min}(\mathfrak{h})\bigr)\; .
\end{align*}
It is also interesting to see whether in this case it follows that
\begin{equation}
  \label{eq:30}
  \mu=\delta_{^{\exists }z_0}\; , \; z_0\in \mathrm{Min}(\mathfrak{h})\; ,
\end{equation}
\emph{i.e.}\ whether the cylindrical Wigner measure has to be concentrated on a classical minimizer. While
statements that yield \eqref{eq:29} has been proved at least in one suitable case \citep[see,
\emph{e.g.},][]{ammari2014jsp
}, to the author's knowledge there are no results of type \eqref{eq:30} in the literature.

\subsubsection{Other perspectives}
\label{sec:other-poss-appl}

There are many other interesting problems that could be studied within this framework. Let us mention very
briefly some of them. One is to study the convergence of quantum to classical ergodicity, and mixing. The
idea is that under suitable conditions the cylindrical Wigner measures of ergodic (mixing) semiclassical
quantum states should be ergodic (mixing) as well. Another is to study the classical limit of Haag-Ruelle
scattering theory for quantum fields. Scattering theory in the classical limit has only been studied for
coherent states in some specific cases \citep{eckmann1975lmp,ginibre1979cmpI,ginibre1979cmp2%
}, and it would be interesting to make a more systematic study for relativistic theories in the
Haag-Ruelle framework using cylindrical Wigner measures and semiclassical techniques. Finally, to
complement for curved spacetimes what we discussed in
\crefrange{sec:semicl-covar-isot}{sec:semicl-time-slice}, it would be interesting to study the
semiclassical behavior of Hadamard states, that play the role of ground states in curved spacetimes
\citep[see][and references thereof contained for more information about Hadamard
states]{kay1991pr,verch1994cmp,radzikowski1996cmpI,dappiaggi2011cmp,gerard2014cmp}. It would also be
interesting to develop a similar framework for the semiclassical analysis of fermionic quantum field
theories, such as Dirac quantum fields. However, for fermions the semiclassical description is rather different,
and it is better captured by a multiscale analysis \citep[see][for additional details]{ammari2017arxiv}.

\section{Weyl C*-algebras}
\label{sec:class-char-stat}

In this section we introduce the Weyl C*-algebra of canonical commutation relations corresponding to an
infinite dimensional Heisenberg group. We then describe some of the properties of composite quantum
systems, consisting of a semiclassical and a fixed part.

\subsection{Infinite-dimensional Heisenberg groups and the algebra of CCR}
\label{sec:non-comm-sett}

Let $(X,\varsigma)\in \mathbf{Symp}_{\mathds{R}}$ be of arbitrary (infinite) dimension. Let us recall that
the symplectic form $\varsigma:X\times X\to \mathds{R}$ is non-degenerate, bilinear and antisymmetric. The
Heisenberg group $\mathbb{H}(X,\varsigma)$ associated to $(X,\varsigma)$ is the space $X\times \mathds{R}$
endowed with the group structure
\begin{equation*}
  \bigl(x,t\bigr)\cdot\bigl(y,s\bigr)=\bigl(x+y,t+s-\varsigma(x,y)
  \bigr)\;. 
\end{equation*}
The set of elements $Z=\bigl\{(0,t),t\in \mathds{R}\bigr\}$ is the center of the group. Among the
representations of the Heisenberg group, one plays a prominent role in quantum theories and semiclassical
analysis: the \emph{Weyl C*-algebra}. From a physical standpoint, this algebra encodes in a natural way
the canonical commutation relations of (bosonic) quantum systems. In semiclassical analysis, it is the
starting point to define quantizations and pseudodifferential calculus. The Weyl C*-algebra is uniquely
defined, up to *-isomorphisms, as the smallest C*-algebra containing the set
\begin{equation*}
  \bigl\{W(x),x\in X\bigr\}\; ,
\end{equation*}
together with the following three properties for its elements
\begin{align*}
  &\bullet\; W(x)\neq 0 &(\forall x\in X)\\
  &\bullet\; W(-x)= W(x)^{*} &(\forall x\in X)\\
  &\bullet\; W(x)W(y)=e^{-i\varsigma(x,y)}W(x+y) &(\forall x,y\in X)
\end{align*}
As in \cref{sec:introduction-1}, we adopt the notation
$\mathbb{W}_1(X,\varsigma)=\textrm{C*}\bigl(\bigl\{W(x),x\in X\bigr\}\bigr)$. From the definition, it
follows that $W(0)=1$ (identity element), and that each $W(x)$ is unitary. Therefore
$\mathbb{W}_1(X,\varsigma)$ is a non abelian unital C*-algebra generated by unitaries. The map
\begin{equation*}
  \begin{aligned}
    X\times \mathds{R} &\longrightarrow \mathbb{W}_1(X,\varsigma)\\
    (x,t)&\longmapsto e^{i t}\,W(x)\\
  \end{aligned}\; ,
\end{equation*}
together with the identification of the group product with the C*-algebra product, provides the unitary
representation of the Heisenberg group in the Weyl C*-algebra. Therefore from now on we will focus on the
Weyl C*-algebra, keeping in mind the underlying Heisenberg group structure.

The non abelian nature of the Heisenberg group, or equivalently of the Weyl algebra, is given by the
symplectic factors $-i\varsigma(x,y)$ in the product. Borrowing an idea from deformation theory, it is
natural to ``measure'' the noncommutativity of the Weyl algebra introducing a real parameter $h\geq 0$
such that when $h>0$ the algebra is non-abelian, and when $h=0$ it becomes abelian. This justifies the
following definition.

\begin{definition}[Weyl deformation]\label{def:17}
  Let $(X,\varsigma)\in \mathbf{Symp}_{\mathds{R}}$. Then the \emph{Weyl deformation}
  $\bigl(\mathbb{W}_h(X,\varsigma)\bigr)_{h\geq 0}$ is a family of C*-algebras. For any $h\geq 0$, the
  algebra $\mathbb{W}_h(X,\varsigma)$ is generated by the set
  \begin{equation*}
    \bigl\{W_h(x),x\in X\bigr\}\; ,
  \end{equation*}
  together with the three properties of its elements
  \begin{align*}
    &\bullet\; W_h(x)\neq 0 &(\forall x\in X)\\
    &\bullet\; W_h(-x)= W_h(x)^{*} &(\forall x\in X)\\
    &\bullet\; W_h(x)W_h(y)=e^{-i h \varsigma(x,y)}W_h(x+y) &(\forall x,y\in X)
  \end{align*}
\end{definition}

\noindent For any $h>0$, the algebra $\mathbb{W}_h(X,\varsigma)$ is *-isomorphic to
$\mathbb{W}_1(X,\varsigma)$, since $W_h(x)=W(h^{1/2}x)$. When $h=0$ however, the algebra
$\mathbb{W}_0(X,\varsigma)$ is an abelian unital C*-algebra of almost periodic functions \citep[see
\cref{sec:states-c-algebra}; and][for an introduction to almost periodic
functions]{bohr1947apf
}. In other words, a Weyl deformation contains infinitely many identical copies of the Weyl C*-algebra,
and a single abelian algebra of almost periodic functions.

\subsection{Tensor product of C*-algebras and partial evaluation}
\label{sec:tensor-product-c}

For physical reasons \citep[see, \emph{e.g.},][]{correggi2017arxiv,correggi2017ahp%
}, we couple the Weyl algebra with another C*-algebra that represent some additional degrees of freedom
that do not behave semiclassically (and therefore do not depend on $h$). Instead of
$\bigl(\mathbb{W}_h(X,\varsigma)\bigr)_{h\geq 0}$, we consider the deformation with the additional degrees
of freedom given by a C*-algebra of physical (either quantum or classical) observables $\mathfrak{A}\in
\mathbf{C^{*}alg}$.
\begin{equation}
  \label{eq:112}
  (\mathfrak{W}_h)_{h\geq 0}=\bigl(\mathbb{W}_h(X,\varsigma)\otimes_{\gamma_h}\mathfrak{A} \bigr)_{h\geq 0}\; ;
\end{equation}
where the index $\gamma_h$ stands for a suitable choice of cross norm for the tensor product C*-algebra
\citep[see, \emph{e.g.},][]{takesaki1979I
}. There are some differences depending on what norm is chosen, as it will be highlighted in the
following. In applications, it is sometimes important to consider the enveloping von Neumann algebra
$\mathbb{W}_h(X,\varsigma)''$ in place of $\mathbb{W}_h(X,\varsigma)$, or other algebras that embed
$\mathbb{W}_h(X,\varsigma)$. The majority of our results extend to any deformation
$(\mathfrak{B}_h\otimes_{\gamma_h} \mathfrak{A})_{h\geq 0}$ such that $\mathbb{W}_h(\,\! ^{\exists }X,\,\!
^{\exists }\varsigma)\overset{^{\exists }\mathrm{w}_h}{\hookrightarrow}\mathfrak{B}_h$. It will be pointed
out explicitly in the text whether a result extends or not to the aforementioned case. Finally, if one is
interested only in the deformation $\bigl(\mathbb{W}_h(X,\varsigma)\bigr)_{h\geq 0}$, it suffices to take
$\mathfrak{A}$ to be the trivial C*-algebra generated by a single element. On $\mathfrak{W}_h$, there are
two natural maps that play an important role, and we call them partial evaluations. We recall that for any
topological linear space $V\in \mathbf{TVS}$, we denote by $V'$ the continuous dual of $V$ (while $V^{*}$
stands for the algebraic dual).

\begin{definition}[Partial evaluation]\label{def:18}
  For any $h\geq 0$, define the \emph{partial evaluation map}
  \begin{equation*}
    \mathbb{E}_{h,1}^{(\cdot)}: \mathfrak{W}_h'\to \mathcal{B}\bigl(\mathbb{W}_h(X,\varsigma),\mathfrak{A}'\bigr)\; ,
  \end{equation*}
  by its action
  \begin{align*}
    \mspace{60mu}&\mathbb{E}_{h,1}^{\omega_h}(w_h)(a)=\omega_h(w_h\otimes a)&\bigl(\forall \omega_h\in \mathfrak{W}_h'\,,\, \forall w_h\in \mathbb{W}_h(X,\varsigma)\,,\, \forall a\in \mathfrak{A}\bigr)\; .
  \end{align*}
  The \emph{partial trace} of the complex quantum state $\omega_h\in \mathfrak{W}_h'$ is the partial
  evaluation of the identity element $\mathbb{E}_{h,1}^{\omega_h}(1)\in \mathfrak{A}'$. The partial
  evaluation $E_{h,2}^{(\cdot )}: \mathfrak{W}_h'\to \mathcal{B}\bigl(\mathfrak{A},
  \mathbb{W}_h(X,\varsigma)'\bigr)$ is defined in a symmetric fashion.
\end{definition}
\noindent In the definition above, $\mathcal{B}(X,Y)$ stands for the space of continuous linear maps from
$X$ to $Y$. We chose to emphasize (perhaps with an heavy notation) the dependence on the semiclassical
parameter $h$, for it will play an important role. The partial evaluation does what its name suggests:
given a (complex) state on the tensor algebra, it evaluates any observable of the first algebra and gives
as output a (complex) state acting on the second algebra
alone. 
The partial evaluation map has some important properties that are summarized in the following proposition
\citep[see, \emph{e.g.},][for a proof]{takesaki1979I
}.
\begin{proposition}\label{prop:26}
  For any $h\geq 0$, the evaluation map $\mathbb{E}_{h,1}$ is an isometry of $\mathfrak{W}_h'$ into
  $\mathcal{B}\bigl(\mathbb{W}_h(X,\varsigma),\mathfrak{A}'\bigr)$. In addition, an element $\omega_h\in
  \mathfrak{W}_h'$ is a state of total mass $m_h$ -- \emph{i.e.}\ $\omega_h\in (\mathfrak{W}_h)'_+$ and
  $\lVert \omega_h \rVert_{\mathfrak{W}_h'}^{}=m_h$ -- if the resulting evaluation
  $\mathbb{E}_{h,1}^{\omega_h}:\mathbb{W}_h(X,\varsigma)\to \mathfrak{A}'$ is completely positive and the
  partial trace $\mathbb{E}_{h,1}^{\omega_h}(1)\in \mathfrak{A}'_+$ satisfies
  \begin{equation*}
    \Bigl\lVert \mathbb{E}_{h,1}^{\omega_h}(1)  \Bigr\rVert_{\mathfrak{A}'}^{}=m_h\; .
  \end{equation*}
  If $\gamma_h$ is the maximal norm $\gamma_{\mathrm{max}}$, then $\mathbb{E}_{h,1}$ is \emph{onto}
  $\mathcal{B}\bigl(\mathbb{W}_h(X,\varsigma),\mathfrak{A}'\bigr)$ and the converse of the second
  statement holds.
\end{proposition}

\subsection{The generating map and regular states.}
\label{sec:generating-map}

Given a state on the Weyl C*-algebra, it is possible to define its generating
functional~\citep[see][]{segal1961cjm
} or noncommutative Fourier transform; in our framework it is not a functional, but a map from $X$ to
$\mathfrak{A}'$. Throughout this section, we take $h>0$ if not specified otherwise.
\begin{definition}[Generating map]\label{def:19}
  Let $\omega_h\in (\mathfrak{W}_h)'_+$ be a state, we define the \emph{generating map}
  $\mathcal{G}_{\omega_h}:X\to \mathfrak{A}'$ by
  \begin{align*}
    \mspace{210mu}&\mathcal{G}_{\omega_h}(x)=\mathbb{E}_{h,1}^{\omega_h}\bigl(W_h(x)\bigr)\qquad&(\forall x\in X)\; .
  \end{align*}
\end{definition}
\noindent The generating map is used to define a very important class of states (and hence its name), the
so-called regular states. Regular states are those that allow for a natural semiclassical description in
term of Wigner measures.
\begin{definition}[Regular states]\label{def:20}
  Let $\omega_h\in (\mathfrak{W}_h)'_+$ be a state, $\mathcal{G}_{\omega_h}$ its generating map. Then
  $\omega_h$ is \emph{regular} iff for any $x\in X$, the $\mathds{R}$-action
  \begin{equation*}
    \mathcal{G}_{\omega_h}(\,\cdot \, x):\mathds{R}\to \mathfrak{A}'
  \end{equation*}
  is continuous when $\mathfrak{A}'$ is endowed with the $\sigma(\mathfrak{A}',\mathfrak{A})$ topology
  (ultraweakly continuous).
\end{definition}
\noindent There are many equivalent definitions of regular states. We also make use of the following one,
that can be proved, \emph{e.g.}, using the properties of the map $\mathbb{E}_{h,2}$ and the equivalent
result for trivial $\mathfrak{A}$
\citep[see][\textsection~5.2.3]{bratteli1997tmp2
}. Let $(R,\varsigma)$ be a finite dimensional real symplectic vector space. We say that a state
$\varrho_h$ on $\mathbb{W}_h(R,\varsigma)\otimes_{\gamma_h} \mathfrak{A}$ is \emph{normal} iff for any
$a\in \mathfrak{A}_+$, $\mathbb{E}_{h,2}^{\varrho_h}(a)$ is a (positive) trace class operator in the
unique irreducible representation of $\mathbb{W}_h(R,\varsigma)$ (the uniqueness up to unitary equivalence
of such representation is guaranteed by Stone-von Neumann's theorem).
\begin{proposition}\label{prop:105}
  A state $\omega_h$ is regular iff for any finite dimensional $R\subset X$ its restriction $\varrho_h$ to
  $\mathbb{W}_h(R,\varsigma)\otimes_{\gamma_h} \mathfrak{A}$ is a normal state. In particular, it follows
  that the generating map of a regular state is ultraweakly continuous when restricted to any finite
  dimensional subspace of $X$.
\end{proposition}

The following result is an extension to our setting of the main result of the aforementioned paper of
\citet{segal1961cjm
}. The idea is that regular states are uniquely determined by the generating map, and the latter is
``almost'' completely positive (up to a complex phase factor) and ultraweakly continuous on finite
dimensional subspaces.
\begin{proposition}\label{prop:27}
  For any $h> 0$, a map $\mathcal{G}_h:X\to \mathfrak{A}'$ is the generating map of a regular state
  $\omega_h\in (\mathfrak{W}_h)'_+$ of partial trace $\alpha_h\in \mathfrak{A}'_+$ only if all the
  restrictions of $\mathcal{G}_h$ to finite dimensional subspaces of $X$ are ultraweakly continuous,
  $\mathcal{G}_h(0)=\alpha_h$ and
  \begin{equation*}
    \sum_{j,k\in J}^{}\mathcal{G}_h(x_j-x_k)e^{ih\varsigma(x_j,x_k)}(a_k^{*}a_j)\geq 0\; ;
  \end{equation*}
  where the $x_j\in X$ are arbitrary as well as the $a_j\in \mathfrak{A}$, and $J$ is any finite index
  set. If in addition $\gamma_h$ is the maximal norm, the converse holds and the map $\mathcal{G}_h$
  uniquely determines $\omega_h$.
\end{proposition}
\begin{remark}\label{rem:22}
  If in $\mathfrak{W}_h$, with maximal norm, we replace $\mathbb{W}_h(X,\varsigma)$ by its enveloping von
  Neumann algebra or any algebra that embeds the Weyl C*-algebra as a subalgebra, $\mathcal{G}_h$ does not
  determine $\omega_h$ uniquely.
\end{remark}
\begin{proof}
  Let us start with the easy ``only if'' part for generic cross norms. Ultraweak continuity follows from
  Proposition~\ref{prop:105}, the other two properties follow from Proposition~\ref{prop:26}: in fact
  $W_h(0)=1$;
  \begin{equation*}
    \sum_{j,k\in F}^{}\mathbb{E}_{h,1}^{\omega_h}\bigl(W_h(x_k)^{*}W_h(x_j)\bigr) (a_k^{*}a_j)\geq 0
  \end{equation*}
  by complete positivity of $\mathbb{E}_{h,1}^{\omega_h}$; and
  $W_h(-x)W_h(y)=e^{ih\varsigma(x,y)}W_h(x-y)$ by definition of the Weyl algebra. To prove the ``if'' part
  and uniqueness when $\gamma_h$ is the maximal cross norm, we act with the generating map on an arbitrary
  $a\in \mathfrak{A}_+$. Since
  $\mathfrak{A}=\mathfrak{A}_+-\mathfrak{A}_++i(\mathfrak{A}_+-\mathfrak{A}_+)$, this suffices to
  characterize the map $\mathcal{G}_h:X\to \mathfrak{A}'$ by linearity. Let us denote by
  $\mathcal{G}^a_h(\cdot )=\mathcal{G}_h(\cdot )(a):X\to \mathds{C}$. By Theorem 1
  of~\citep{segal1961cjm
  }, to $\mathcal{G}^a_h$ there corresponds a unique regular state $\varrho^a_h\in
  \bigl(\mathbb{W}_h(X,\varsigma)\bigr)'_+$ such that $\mathcal{G}_h^a(\cdot )=\varrho_h^a\bigl(W_h(\cdot
  )\bigr)$. By the last property of $\mathcal{G}_h$ this defines a unique completely positive map
  $\varrho_h^{(\cdot )}:\mathfrak{A}\to \mathbb{W}_h(X,\varsigma)'$. Therefore the analogous of
  Proposition~\ref{prop:26} for $\mathbb{E}_{h,2}$ yields that
  $\omega_h=\mathbb{E}_{h,2}^{-1}(\varrho_h^{(\cdot )})$ is a positive regular measure of total mass
  $\mathcal{G}_h(0)$, uniquely determined by $\mathcal{G}_h$.
\end{proof}

\section{Semiclassical compactness for families of regular states}
\label{sec:semicl-analys-regul}

In this section we study the semiclassical behavior $h\to 0$ of families of regular states. In particular,
we prove compactness, in a suitable topology, of families of generating maps.

\subsection{Finite dimensional Weyl algebra}
\label{sec:finite-dimens-result}

By Proposition~\ref{prop:105}, if $X$ is finite dimensional then any state that is normal with respect to
the Schrödinger representation is regular. In other words, the semiclassical analysis of regular states on
$\mathbb{W}_h(X,\varsigma)$, with $X\cong \mathds{R}^d\times (\mathds{R}^d)'$, reduces to the
semiclassical analysis of families in $\mathfrak{S}^1\bigl(L^2 (\mathds{R}^d )\bigr)_+$, the cone of
positive trace class operators. By linearity, it is sufficient to study families of vectors
$(u_h)_{h>0}\subset L^2 (\mathds{R}^d )$ as $h\to 0$ (with norms uniformly bounded with respect to
$h$). The commutative objects corresponding to the family $(u_h)_{h>0}$ are called Wigner measures, and
are finite Borel measures on $\mathds{R}^d\times (\mathds{R}^d)'$ (equivalently on $X$). In addition, to
any family $(u_h)_{h>0}$ with uniformly bounded norms there corresponds at least one Wigner measure. We
formulate this well-known result \citep[][Théorème III.1]
{lions1993rmi
} in a slightly different way, that is better suited for our algebraic approach.

Let $T$ be a topological space, and $S$ a set; we denote by $T^S_{\mathrm{s}}$ the space of functions from
$S$ to $T$ with the topology of simple (pointwise) convergence, that coincides with the product
topology. We use the following notation for semiclassical generalized sequences (nets):
$(\diamond_{h_{\beta}})_{\beta\in B}$, $h_{\beta}\to 0$, is a net of objects indexed by the directed set
$B$, such that $\diamond_{h_{\beta}}=\diamond_{h_{\beta'}}$ whenever $h_{\beta}=h_{\beta'}$,
$(h_{\beta})_{\beta\in B}$ is relatively compact in the topology of $\mathds{R}$, and
\begin{equation*}
  \lim_{\beta\in B}h_{\beta}=0\; .
\end{equation*}
A point $x$ is a cluster point for $(x_{\beta})_{\beta\in B}$ if there exists a subnet $(x_b)_{b \in
  \underline{B}}$ such that $\lim_{b\in \underline{B}}x_b=x$. A point $x$ is a sequential cluster point
for $(x_{\beta})_{\beta\in B}$ if it is a cluster point and the corresponding subnet is a subsequence
(\emph{i.e.}\ if $\underline{B}=\mathds{N}$).

A uniformly bounded net of regular states $(\varrho_{h_{\beta}})_{\beta\in B}$ on
$\mathbb{W}_{h_{\beta}}(X,\varsigma)\otimes_{\gamma_{h_{\beta}}} \mathfrak{A}$, with $X\cong X'$ finite
dimensional, defines a net of linear maps $(H_{\varrho_{h_{\beta}}})_{\beta\in B}$ from $\mathfrak{A}$ to
the dual $C_0^{\infty}(X')'$ of smooth compactly supported functions by the Bochner integral
\begin{align}
  \label{eq:122}
  &H_{\varrho_{h}}(a)(\varphi)=\varrho_h\Bigl(\int_{X}^{}\hat{\varphi}(x) W_{h}(\pi x) \mathrm{d}L_X(x)\otimes a\Bigr) &\bigl(\forall a\in \mathfrak{A}\,,\, \forall \varphi\in C_0^{\infty}(X')\bigr)\; ,
\end{align}
where $L_X$ is the Lebesgue measure of $X$. In addition, each $H_{\varrho_{h_{\beta}}}$ is also a
continuous linear map from $\mathfrak{A}$ to $C_0(X')'$, the dual of the compactly supported continuous
functions.

\begin{proposition}
  \label{prop:102}

  Let $(\varrho_{h_{\beta}})_{\beta\in B}$, $h_{\beta}\to 0$, be a net of regular states that acts on
  $\bigl(\mathbb{W}_{h_{\beta}}(X,\varsigma)\otimes_{\gamma_{h_{\beta}}} \mathfrak{A}\bigr)_{\beta\in B}$
  with $X$ \emph{of finite dimension}. If
  \begin{equation*}
    \sup_{\beta\in B}\,\varrho_{h_{\beta}}\bigl(W_{h_{\beta}}(0)\bigr)<\infty\; ,
  \end{equation*}
  then for any $a\in \mathfrak{A}$, $\Bigl\{H_{\varrho_{h_{\beta}}}(a),\beta\in B\Bigr\}$ is relatively
  compact in the $\sigma\bigl(C_0(X')',C_0(X)\bigr)$ topology, and the cluster points of any subnet
  $\bigl(H_{\varrho_{h_b}}(a_+)\bigr)_{b\in \underline{B}}$, $h_b\to 0$, with $a_+\in \mathfrak{A}_+$ are
  \emph{positive} linear functionals. Thus by Riesz-Markov's theorem they can be identified with finite
  Radon measures on $X'$.

  Let us fix now $a\in \mathfrak{A}_+$, and a subnet $\bigl(H_{\varrho_{h_b}}(a)\bigr)_{b\in \underline{B}}$, $h_b\to
  0$, converging to $\mu_a\in \mathcal{M}_{\mathrm{rad}}(X')$. In addition, denote by $(\psi_{\varepsilon})_{\varepsilon\in (0,1)}\subset
  C_0^{\infty}(X')$ an approximated identity function on $X'$. The following four statements are equivalent:
  \begin{align}
    \label{eq:31}\tag{i}
    &\bigl(\mathcal{G}_{\varrho_{h_{b}}}(\,\cdot\,)(a)\bigr)_{b\in \underline{B}}\overset{\mathds{C}^X_{\mathrm{s}}}{\underset{h_b\to 0}{\longrightarrow}}\hat{\mu}_a(\,\cdot \,)&; \\
    \label{eq:32}\tag{ii}
    &\lim_{\varepsilon\to 0}\lim_{b\in \underline{B}}\varrho_{h_b}\bigl(W_{h_b}(x)\otimes a-\mathrm{Op}_{\frac{1}{2}}^{h_b}(^{\exists }\psi_{\varepsilon}(\cdot )\, e^{2ix(\cdot )})\otimes a\bigr)=0&\bigl(\forall x\in X\bigr)\,;\\
    \label{eq:33}\tag{iii}
    &\lim_{\varepsilon\to 0}\lim_{b\in \underline{B}}\varrho_{h_b}\bigl(a-\mathrm{Op}_{\frac{1}{2}}^{h_b}(^{\exists }\psi_{\varepsilon})\otimes a\bigr)=0&;\\
    \label{eq:34}\tag{iv}
    &\lim_{b\in \underline{B}}\varrho_{h_b}\bigl(W_{h_b}(0)\otimes a\bigr)=\mu_a(X')&.
  \end{align}
\end{proposition}

The cluster points $\mu_a$ of \cref{prop:102}, considered as Radon measures, are the so-called Wigner or
semiclassical measures.

\subsection{Compactness and convergence in infinite-dimensional Weyl algebras}
\label{sec:compactness}

Let $(X,\varsigma)\in \mathbf{Symp}_{\mathds{R}}$ be of arbitrary dimension, and let
\begin{equation*}
  (X_{\lambda})_{\lambda\in F}\subset 2^X
\end{equation*}
be a collection of finite dimensional symplectic subspaces, indexed by a set $F$. In addition, let
$\omega_h$ be a regular state on $\mathfrak{W}_h$. In analogy with \eqref{eq:122}, for any $\lambda\in F$,
define the map $H^{(\lambda)}_{\omega_h}$ from $\mathfrak{A}$ to $C_0^{\infty}(X'_{\lambda})'$ by
\begin{align}
  \label{eq:123}
  &H^{(\lambda)}_{\omega_h}(a)(\varphi_{\lambda})=\omega_h\Bigl(\int_{X_{\lambda}}^{}\hat{\varphi}_{\lambda}(x) W_{h}(\pi x) \mathrm{d}L_{X_{\lambda}}(x)\otimes a\Bigr) &\bigl(\forall a\in \mathfrak{A}\, , \forall \varphi_{\lambda}\in C_0^{\infty}(X'_{\lambda})\bigr)\; .
\end{align}
Let us consider now the set
\begin{equation}
  \label{eq:124}
  H_B=\prod_{\lambda\in F}\prod_{a\in \mathfrak{A}_+} \Bigl\{H^{(\lambda)}_{\omega_{h_{\beta}}}(a), \beta\in B\Bigr\}\; ,
\end{equation}
and denote by $\mathfrak{P}_F$ the product topology on $H_B$, with each
$\Bigl\{H^{(\lambda)}_{\omega_{h_{\beta}}}(a), \beta\in B\Bigr\}$ endowed with the weak
$\sigma\bigl(C_0(X_{\lambda}')',C_0(X_{\lambda})\bigr)$ topology.
\begin{lemma}
  \label{lemma:12}
  Let $(\omega_{h_{\beta}})_{\beta\in B}$, $h_{\beta}\to 0$, be a net of regular states on
  $(\mathfrak{W}_{h_{\beta}})_{\beta\in B}$. If
  \begin{equation*}
    \sup_{\beta\in B}\,\omega_{h_{\beta}}\bigl(W_{h_{\beta}}(0)\bigr)<\infty\; ,
  \end{equation*}
  then $H_B$ is $\mathfrak{P}_F$-relatively compact, and the cluster points of any subnet
  \begin{equation*}
    \Bigl(\mathfrak{h}_{h_b}\Bigr)_{b\in\underline{B}}\subset H_B\; ,
  \end{equation*}
  $h_b\to 0$, can be identified with a unique family of vector valued measures
  \begin{equation*}
    (\mu_{\lambda})_{\lambda\in F}\;, \; \forall \lambda\in F\;,\; \mu_{\lambda}\in \mathcal{M}_{\mathrm{rad}}(X_{\lambda}';\mathfrak{A}'_+)\; .
  \end{equation*}
\end{lemma}
\begin{proof}
  Compactness in the product topology follows immediately from Proposition~\ref{prop:102}. In addition,
  the cluster points $\mathfrak{h}$ of $(\mathfrak{h}_{h_b})_{b\in \underline{B}}$, $h_b\to 0$, have the
  form
  \begin{equation*}
    a\mapsto \mathfrak{h}(a)\;,\; \mathfrak{h}(a)=(\mathfrak{h}_{\lambda}(a))_{\lambda\in F}\; ,
  \end{equation*}
  with each $\mathfrak{h}_{\lambda}(a)$ a positive element of $C_0(X_{\lambda}')'$. Therefore by
  Riesz-Markov's theorem, for any $a\in \mathfrak{A}_+$ there exists a unique $\mu_{\lambda}(a)\in
  \mathcal{M}_{\mathrm{rad}}(X_{\lambda}';\mathds{R}^+)$ whose action on $C_0(X_{\lambda}')$ agrees with
  $\mathfrak{h}_{\lambda}(a)$. In other words, there exists an injective map $\mathcal{R}$ such that
  \begin{equation}
    \label{eq:125}
    \bigl(a\mapsto \mathfrak{h}_{\lambda}(a)\bigr)_{\lambda\in F}\overset{\mathcal{R}}{\longmapsto} \bigl(a\mapsto \mu_{\lambda}(a)\bigr)_{\lambda\in F}\; .
  \end{equation}
  In addition, by linearity it follows that for any $\lambda\in F$,
  \begin{equation*}
    \bigl(a\mapsto \mu_{\lambda}(a)\bigr)\in \mathrm{Hom}_{\mathrm{mon}}\bigl(\mathfrak{A}_+,\mathcal{M}_{\mathrm{rad}}(X_{\lambda}';\mathds{R}^+) \bigr)\; .
  \end{equation*}
  Therefore by Theorem~\ref{thm:101} there exists a bijection $\mathcal{P}$ such that, with
  $\mathcal{Q}=\mathcal{P}\circ\mathcal{R}$,
  \begin{equation}
    \label{eq:126}
    \bigl(a\mapsto \mathfrak{h}_{\lambda}(a)\bigr)_{\lambda\in F}\overset{\mathcal{Q}}{\longmapsto} \bigl(\mu_{\lambda}\bigr)_{\lambda\in F}\; ;\; \forall \lambda\in F\; ,\; \mu_{\lambda}\in \mathcal{M}_{\mathrm{rad}}(X_{\lambda}';\mathfrak{A}'_+)\; .
  \end{equation}
\end{proof}

Let us now turn attention to the converging nets whose cluster points have no loss of mass. For these
nets, there is an easier characterization of cluster points by means of the generating maps, introduced in
\cref{sec:generating-map}. The set of generating maps is always relatively compact with respect to the
topology of simple convergence. In fact, let $(\mathfrak{W}_h)_{h\geq 0}$ be the tensor Weyl deformation
introduced in \cref{sec:tensor-product-c}, with $X$ of arbitrary (infinite) dimension. We are interested
in \emph{semiclassical states} (see \cref{def:1}), \emph{i.e.}\ generalized bounded sequences of regular
states
\begin{equation*}
  (\omega_{h_{\beta}})_{\beta\in B}\quad,\quad \lim_{\beta\in B}h_{\beta}=0\quad ,\quad \sup_{\beta\in B}\, \omega_{h_{\beta}}(1)=m<\infty\; .
\end{equation*}
Let $\mathfrak{A}'_{\mathfrak{A}}$ be the continuous dual of $\mathfrak{A}$ endowed with the ultraweak
$\sigma(\mathfrak{A}',\mathfrak{A})$ topology. Let us denote by $G_B\subset
(\mathfrak{A}'_{\mathfrak{A}})^X$ and $G_{B}(x)\subset \mathfrak{A}'_{\mathfrak{A}}$ the following sets:
\begin{equation*}
  G_{B}=\Bigl\{\mathcal{G}_{\omega_{h_{\beta}}},\beta\in B\Bigr\}\; ;\; G_{B}(x)=\Bigl\{\mathcal{G}_{\omega_{h_{\beta}}}(x),\beta\in B\Bigr\}\;,\; x\in X\; .
\end{equation*}
The first result is that the family of images $G_B(x)$ is pointwise compact for any $x\in X$, and
therefore $G_B$ is relatively compact as a subset of the space of functions from $X$ to
$\mathfrak{A}'_{\mathfrak{A}}$.
\begin{lemma}\label{lemma:9}
  Let $(\omega_{h_{\beta}})_{\beta\in B}$ be a semiclassical state in the Weyl deformation
  $(\mathfrak{W}_{h_{\beta}})_{\beta\in B}$. Then $G_{B}(x)$ is relatively compact for any $x\in X$. It
  then follows that $G_{B}$ is relatively compact as a subset of
  $(\mathfrak{A}'_{\mathfrak{A}})^X_{\mathrm{s}}$.
\end{lemma}
\begin{proof}
  It follows from \cref{def:19} of the generating map -- since the Weyl operators are unitary -- that for
  any $x\in X$, $\beta\in B$ and $a\in \mathfrak{A}$
  \begin{equation*}
    \bigl\lvert  G_{\omega_{h_{\beta}}}(x)(a) \bigr\rvert_{}^{}\leq \lVert\omega_{h_{\beta}}  \rVert_{\mathfrak{W}_{h_{\beta}}'}^{}\lVert a  \rVert_{\mathfrak{A}}^{}\leq m \lVert a  \rVert_{\mathfrak{A}}^{}\; .
  \end{equation*}
  Therefore $G_{B}(x)$ is contained in the ball of radius $m$ of $\mathfrak{A}'$, and therefore it is
  relatively compact in the ultraweak topology by Banach-Alaoglu's theorem.
\end{proof}

Let us now consider a converging net $(\mathcal{G}_{\omega_{h_{b}}})_{b\in \underline{B}}$, $h_{b}\to 0$,
in the relatively compact set $G_B$, and denote its restriction to any finite dimensional subspace
$R\subset X$ by $\bigl(\mathcal{G}_{\omega_{h_{b}}}\bigr\rvert_R \bigr)_{b\in \underline{B} }$. Now, let
$a\in \mathfrak{A}$; then
\begin{equation*}
  a=a^+_{\mathrm{R}}-a^-_{\mathrm{R}}+i(a^+_{\mathrm{I}}-a^-_{\mathrm{I}})\;,\quad a^+_{\mathrm{R}},a^-_{\mathrm{R}},a^+_{\mathrm{I}},a^-_{\mathrm{I}}\in \mathfrak{A}_+\; . 
\end{equation*}
A function $f\in (\mathfrak{A}'_{\mathfrak{A}})^R$ is continuous iff $f(\cdot )(a):R\to \mathds{C}$ is
continuous for any $a\in \mathfrak{A}$. By linearity it follows, using the decomposition above, that $f\in
(\mathfrak{A}'_{\mathfrak{A}})^R$ is continuous iff $f(\cdot )(a_+)$ is continuous for any $a_+\in
\mathfrak{A}_+$. For any $b\in \underline{B}$, define
\begin{equation*}
  \mathcal{G}_{\omega_{h_{b}}}^{a_+}\bigr\rvert_R: R\to \mathds{C}\quad , \quad \mathcal{G}_{\omega_{h_{b}}}^{a_+}\bigr\rvert_R(r)= \mathcal{G}_{\omega_{h_{b}}}(r)(a_+)\; (\forall r\in R)\; .
\end{equation*}

\begin{lemma}
  \label{lemma:11}
  For any $a_+\in \mathfrak{A}_+$ and $(R_0,\varsigma)\subset (X,\varsigma)$ a symplectic subspace,
  $\mathcal{G}_{\omega_{h_{b}}}^{a_+}\bigr\rvert_{R_0}$ is the generating functional of the normal state
  $\mathbb{E}_{h_{b},2}^{\omega_{h_{b}}}(a_+)\bigr\rvert_{\mathbb{W}_{h_{b}}(R_0,\varsigma)}$.
\end{lemma}
\begin{proof}
  By Proposition~\ref{prop:105}, since $\omega_{h_b}$ is regular then the restricted states
  \begin{equation*}
    \mathbb{E}_{h_b,2}^{\omega_{h_b}}(a_+)\bigr\rvert_{\mathbb{W}_{h_b}(R_0,\varsigma)}
  \end{equation*}
  are normal with respect to the Schrödinger representation. It is easy to check that, for any $r\in R_0$,
  \begin{equation*}
    \mathcal{G}_{\omega_{h_b}}(r)(a_+)=\mathcal{G}_{\mathbb{E}_{h_b,2}^{\omega_{h_b}}(a_+)}(r)\; .
  \end{equation*}
\end{proof}

As in Proposition~\ref{prop:102}, let us always denote by
\begin{equation}
  \label{eq:128}
  (\psi_{\varepsilon}^{(\lambda)})_{\varepsilon\in (0,1)}\subset C_0^{\infty}(X'_{\lambda})=:\mathcal{D}(X'_{\lambda})
\end{equation}
a mollifying sequence, \emph{i.e.}\ a sequence of smooth functions converging in the
$\sigma\bigl(\mathcal{D}'(X'_{\lambda}),\mathcal{D}(X'_{\lambda})\bigr)$ topology to the delta
distribution.
\begin{proposition}
  \label{prop:7}
  Let $(X_{\lambda})_{\lambda \in F}\subset 2^X$ be a family of finite dimensional symplectic subspaces
  such that $\bigcup_{\lambda\in F}X_{\lambda}=X$, and $(\omega_{h_{\beta}})_{\beta\in B}$, a
  semiclassical quantum state on $(\mathfrak{W}_{h_{\beta}})_{\beta\in B}$. In addition, suppose that
  $(\mathfrak{h}_{h_b})_{b\in \underline{B}}\subset H_B$, $h_b\to 0$, is $\mathfrak{P}_F$-convergent, with
  limit $(\mu_{\lambda})_{\lambda\in F}$ (see Lemma~\ref{lemma:12}). Then the following four statements
  are equivalent:
  \begin{align}
    \label{eq:35}\tag{i}
    &\mathcal{G}_{\omega_{h_b}}\underset{h_b\to 0}{\overset{\mathds{C}^X_{\mathrm{s}}}{\longrightarrow}} g\quad,\quad g\bigr\rvert_{X_{\lambda}}=\hat{\mu}_{\lambda}\in C(X_{\lambda},\mathfrak{A}_{\mathfrak{A}}') &; \\
    \label{eq:36}\tag{ii}
    &\begin{aligned}
      &\lim_{\varepsilon\to 0}\lim_{b\in \underline{B}}\omega_{h_b}\bigl(W_{h_b}(x)\otimes a-\mathrm{Op}_{\frac{1}{2}}^{h_b}(^{\exists }\psi_{\varepsilon}^{(\lambda)}(\cdot )\, e^{2ix(\cdot )})\otimes a\bigr)=0&\bigl(\forall \lambda\in F\,,\\ & &\mspace{-60mu}\forall x\in X_{\lambda}\,,\,\forall a\in \mathfrak{A}_+\bigr)
    \end{aligned}
    &;\\
    \label{eq:37}\tag{iii}
    &\begin{aligned}
      &\lim_{\varepsilon\to 0}\lim_{b\in \underline{B}}\omega_{h_b}\bigl(a-\mathrm{Op}_{\frac{1}{2}}^{h_b}(^{\exists }\psi_{\varepsilon}^{(\lambda)})\otimes  a\bigr)=0&\mspace{88mu}\bigl(\forall \lambda\in F\,,\,\forall a\in \mathfrak{A}_+\bigr)
    \end{aligned}
    &;\\
    \label{eq:38}\tag{iv}
    &\lim_{b\in \underline{B}}\omega_{h_b}\bigl(W_{h_b}(0)\bigr)=\mu_{\lambda}(0)&.
  \end{align}
  A semiclassical quantum state that satisfies \crefrange{eq:35}{eq:38} is called a state with \emph{no
    loss of mass}.
\end{proposition}
\begin{proof}
  For any $\lambda\in F$, $g\bigr\rvert_{X_{\lambda}}$ is continuous iff
  $g\bigr\rvert_{X_{\lambda}}(a_+)\in C(X_{\lambda},\mathds{C})$ for any $a_+\in \mathfrak{A}_+$. Let us
  consider the net $(\mathcal{G}_{\omega_{h_{b}}}^{a_+}\rvert_{X_{\lambda}})_{b\in
    \underline{B}}$. 
  By Lemma~\ref{lemma:11} it is the family of generating functionals of the net of normal states
  $\bigl(\mathbb{E}_{h_{b},2}^{\omega_{h_{b}}}(a_+)\rvert_{\mathbb{W}_{h_{b}}(X_{\lambda},\varsigma)}\bigr)_{b\in
    \underline{B}}$ on $\bigl(\mathbb{W}_{h_b}(X_{\lambda},\varsigma)\bigr)_{b\in
    \underline{B}}$. Therefore by \cref{prop:102} it follows that the statements \eqref{eq:35} to
  \eqref{eq:38} are equivalent.
\end{proof}

Hence we have proved that given a semiclassical net of quantum states, and a family
$(X_{\lambda})_{\lambda\in F}$ of finite dimensional subspaces of $X$, there is always at least one family
of Radon vector measures, on each $X_{\lambda}'$, associated to it. In addition, if $\bigcup_{\lambda\in
  F}X_{\lambda}=X$, and each measure $\mu_{\lambda}$ does not lose mass, then the corresponding net of
generating maps converges, and the limit is ultraweakly continuous when restricted to the finite
dimensional subsets $(X_{\lambda})_{\lambda\in F}$. Let us formulate this generalization of
\cref{prop:102} to spaces of arbitrary dimension as a theorem. This theorem is a generalization of
\cref{thm:1,thm:4}, as it will become clearer after the introduction of the topologies of semiclassical
convergence in \cref{sec:cylindr-wign-meas}.

\begin{thm}\label{thm:14}
  Let $(X,\varsigma)\in \mathbf{Symp}_{\mathds{R}}$, $\mathfrak{A}\in \mathbf{C^{*}alg}$, and
  $(\mathfrak{W}_h)_{h\geq 0}$ the corresponding Weyl deformation \eqref{eq:112}. For any semiclassical
  quantum state $(\omega_{h_{\beta}})_{\beta\in B}$, and for any collection $(X_{\lambda})_{\lambda\in F}$
  of finite dimensional symplectic subspaces of $X$, there exists a nonempty set of cluster points for
  $H_B$ in the $\mathfrak{P}_F$ product topology, and any cluster point is identified with a family
  $(\mu_{\lambda})_{\lambda\in F}$ of finite $\mathfrak{A}'_+$-valued Radon measures on each
  $X'_{\lambda}$.

  In addition, there exists a nonempty set of cluster points of the family of generating maps
  $(\mathcal{G}_{\omega_{h_{\beta}}})_{\beta\in B}$. Each cluster point $g$ satisfies:
  \begin{equation*}
    \sum_{j,k\in J}^{}g(x_j-x_k)(a_k^{*}a_j)\geq 0\; ;
  \end{equation*}
  where the $x_j\in X$ are arbitrary as well as the $a_j\in \mathfrak{A}$, and $J$ is any finite index
  set. For a convergent net $\mathfrak{h}_{h_b}\to_{\mathfrak{P}_F} (\mu_{\lambda})_{\lambda\in F} $, if
  $\,\bigcup_{\lambda\in F}X_{\lambda}=X$ then the following four statements are equivalent:
  \begin{align}
    \label{eq:39}\tag{i}
    &\mathcal{G}_{\omega_{h_b}}\underset{h_b\to 0}{\overset{\mathds{C}^X_{\mathrm{s}}}{\longrightarrow}} g\quad,\quad g\bigr\rvert_{X_{\lambda}}=\hat{\mu}_{\lambda}\in C(X_{\lambda},\mathfrak{A}_{\mathfrak{A}}') &; \\
    \label{eq:40}\tag{ii}
    &\begin{aligned}
      &\lim_{\varepsilon\to 0}\lim_{b\in \underline{B}}\omega_{h_b}\bigl(W_{h_b}(x)\otimes a-\mathrm{Op}_{\frac{1}{2}}^{h_b}(^{\exists }\psi_{\varepsilon}^{(\lambda)}(\cdot )\,e^{2ix(\cdot )})\otimes a\bigr)=0&\bigl(\forall \lambda\in F\,,\\ & & \mspace{-60mu} \forall x\in X_{\lambda}\,,\,\forall a\in \mathfrak{A}_+\bigr)
    \end{aligned}
    &;\\
    \label{eq:41}\tag{iii}
    &\begin{aligned}
      &\lim_{\varepsilon\to 0}\lim_{b\in \underline{B}}\omega_{h_b}\bigl(a-\mathrm{Op}_{\frac{1}{2}}^{h_b}(^{\exists }\psi_{\varepsilon}^{(\lambda)})\otimes  a\bigr)=0&\mspace{88mu}\bigl(\forall \lambda\in F\,,\,\forall a\in \mathfrak{A}_+\bigr)
    \end{aligned}
    &;\\
    \label{eq:42}\tag{iv}
    &\lim_{b\in \underline{B}}\omega_{h_b}\bigl(W_{h_b}(0)\bigr)=\mu_{\lambda}(0)&.
  \end{align}
\end{thm}
\begin{corollary}
  \label{cor:3}
  If $F=\mathds{N}$, and the resulting sequence $(X_n,\varsigma)_{n\in \mathds{N}}$ of symplectic
  subspaces of $(X,\varsigma)$ satisfies $\bigcup_{n\in \mathds{N}}X_n=X$, then from any net
  $(\omega_{h_{\beta}})_{\beta\in B}$ it is possible to extract a \emph{subsequence} $(\omega_{h_j})_{j\in
    \mathds{N}}$ such that
  \begin{equation*}
    \lim_{j\to \infty}\mathfrak{h}_{h_j}=\mathfrak{h}\;,\;     \lim_{j\to \infty}\mathcal{G}_{\omega_{h_j}}=g\; .
  \end{equation*}
\end{corollary}
\begin{proof}
  The proof is the standard diagonal trick. Let us prove it for the generating map $\mathcal{G}$. On
  $X_0$, it is possible to extract the subsequence $h_{\beta_k}^0$ such that
  $\mathcal{G}_{\omega_{h_{\beta_k}^0}}\bigr\rvert_{X_0}\to g\bigr\rvert_{X_0}\in
  C(X_0,\mathfrak{A}'_{\mathfrak{A}})$. On $X_1$ it is possible to extract a subsequence $h_{\beta_k}^1$
  of $h_{\beta_k}^0$ such that $\mathcal{G}_{\omega_{h_{\beta_k}^1}}\bigr\rvert_{X_1}\to
  g\bigr\rvert_{X_1}\in C(X_1,\mathfrak{A}'_{\mathfrak{A}})$, and so on. Therefore the diagonal sequence
  $h_j=h_{\beta_j}^j$ yields the pointwise convergence to $g\in (\mathfrak{A}'_{\mathfrak{A}})^{X}$, such
  that $g\rvert_{X_n}\in C(X_n,\mathfrak{A}'_{\mathfrak{A}})$ for any $n\in \mathds{N}$.
\end{proof}

\section{Cylindrical Wigner measures and topologies of semiclassical convergence}
\label{sec:cylindr-wign-meas}

\cref{thm:14} above proves that any semiclassical quantum state has at least a classical
counterpart. However, it would be suitable to identify such classical object more explicitly. If $X$ is
finite dimensional, the classical characterization given by \cref{prop:102} is satisfactory, and we would
like to obtain an analogous one for spaces $X$ of arbitrary dimension. If $X$ is infinite dimensional
however, the structure is richer, and the situation is more involved. First of all, there are many
inequivalent topologies that are admissible on $X$, and it is therefore possible to put $X$ in duality
with many spaces $V$ in a way such that $X\cong V'$. Analogously, it is possible to interpret $X$ as a
subset $\mathcal{A}$ of the set of functions on some set $A$ in many ways. In this section, we discuss
mostly how duality in linear spaces can be used to identify the classical counterparts of semiclassical
states with \emph{cylindrical measures}. Let us recall that a cylindrical measure on a topological vector
space $V\in \mathbf{TVS}$ can be seen equivalently as a finitely additive inner regular measure on the
algebra of cylinders (induced by the duality of $V$ and $V'$), or as a projective family of Borel measures
in the space of finite dimensional quotients of $V$. A concise definition of vector valued cylindrical
measures in topological vector spaces is given in \cref{sec:locally-conv-spac}. In the same section,
vector-valued Bochner's theorem is also discussed (\cref{thm:103}). Thanks to \cref{thm:103}, it is
possible to identify in a unique fashion cylindrical measures with their Fourier transform (or
characteristic maps), \emph{i.e.} with completely positive maps on $V'$ that are continuous when
restricted to any finite dimensional subspace of $V'$. In the light of \cref{thm:14}, it is clear why they
can be considered as the natural classical counterpart of semiclassical states. Since, as it is discussed
in \cref{sec:cylindrical-measures} (~\cref{prop:1}), the semiclassical description in terms of cylindrical
measures always exists and it is unique up to isomorphisms, we focus on measures acting on topological
vector spaces. The corresponding results for measures defined on more general sets can be easily deduced
using the isomorphisms yielded by Bochner's theorem.

\subsection{Dual topologies and the definition of semiclassical measures}
\label{sec:locally-conv-spac-1}

Let $W,Y\in \mathbf{VS}_{\mathds{R}}$ be two real vector spaces; $W$ and $Y$ are in $B$-duality if there
exists a bilinear form $B:W\times Y\to \mathds{R}$. Given $W$ and $Y$ in $B$-duality, the duality
separates points in $W$ (respectively $Y$) iff $W\ni w\mapsto B(w ,\cdot )\in Y^{*}$ (respectively $Y\ni
y\mapsto B(\cdot ,y )\in W^{*}$) is injective. Let us recall that $Z^{*}$ denotes the algebraic dual of
the linear space $Z$. The duality between linear spaces induces many important and useful
topologies. Before discussing the definition of cylindrical Wigner measures, let us briefly outline some
classical results about duality in linear spaces. Given a linear space $W$, the algebraic dual $W^{*}$ is
in duality with $W$, by the canonical bilinear form $B(w,\xi )=\xi(w)$ on $W\times W^{*}$. In addition,
the duality separates points in both $W$ and $W^{*}$. Using the injective maps $d_B: y\mapsto B(\cdot ,y
)$ and $s_B:w\mapsto B(w,\cdot )$, it is often convenient to identify $Y$ or $W$ (in separating duality
with $W$ or $Y$) with a subspace of $W^{*}$ or $Y^{*}$.

In our setting, given a dual pair $W$ and $Y$ it is useful to characterize the topologies on $W$ that are
\emph{compatible with the duality}, \emph{i.e.}\ the topologies on $W$ such that $d_B$ is a bijection of
$Y$ onto the continuous dual $W'\subset W^{*}$. It turns out that there are coarsest and finest locally
convex topologies on $W$ compatible with the duality. Such topologies are the weak $\sigma(W,Y)$ and the
Mackey $\tau(W,Y)$ topology respectively. Mackey's theorem proves that any locally convex topology on $W$
that is compatible with the duality is finer than the weak and coarser than the Mackey topology. The
strong topology $\beta(W,Y)$ is in general finer than the Mackey topology and therefore not compatible
with the duality. However if $W$ has a given (initial) topology, $\beta(W,W')=\tau(W,W')$ whenever $W$ is
barrelled (\emph{e.g.}, if $W$ is a Fréchet space). If the topological spaces $W$ and $Y$ are in
compatible duality, the bijection $Y\cong W'$ is implicitly understood.

Let us now consider the setting of \cref{thm:14}. As it will become clear in the following, if $X\cong V'$
for some topological real linear space $V$, then we can identify the classical counterpart of regular
quantum states with cylindrical measures on $V$. This can be done, however, in many different ways. First
of all, there may be (infinitely) many spaces $V$ for which there exists a $B: V\times X\to \mathds{R}$
that puts $V$ and $X$ in $B$-duality, separating points in $X$ (the latter requirement is necessary to
have compatible topologies on $V$). For example $X^{*}$, or any of its subspaces that separates points in
$X$. Once we have fixed such $V$, all the locally convex topologies finer than $\sigma(V,X)$ and coarser
than $\tau(V,X)$ are compatible with the identification $X\cong V'$. From the physical perspective, this
indicates that there may exist many equivalent effective descriptions of a quantum field theory as a
classical theory, due to the rich structure of infinite dimensional linear spaces. Since the starting
point is the quantum Heisenberg group associated to $(X,\varsigma)$, it could seem natural to choose
$V=X^{*}$ with the coarsest locally convex topology $\sigma(X^{*},X)$ as the classical space for the
definition of Wigner measures (and \emph{a fortiori} this choice proves that there always exists a space
of classical fields). As we will see, there are situations in which such a choice implies that all the
Wigner cylindrical measures are in fact true (Borel) Radon measures on $V$.

Let us now introduce the appropriate topologies to prove convergence of regular states to Wigner
measures. Let $V$ be a topological space in compatible duality with $X$. Let us denote by
$\mathcal{R}_h\subset (\mathfrak{W}_h)'_+$ the set of regular states on $\mathfrak{W}_h$, and by
\begin{equation*}
  \mathcal{M}_{\mathrm{cyl}}(V;\mathfrak{A}'_+):=\mathcal{M}_{\mathrm{cyl}}(V,V';\mathfrak{A}'_+)
\end{equation*}
the cone-valued cylindrical measures on $V$ -- see \cref{sec:definitions,sec:class-concr-exampl} for a
precise definition of cone-valued measures and for the admissibility of $\mathfrak{A}'_+$ as a cone. Let
us denote the union of the aforementioned sets by
\begin{equation}
  \label{eq:1011}
  S(V,X,\mathfrak{A})=\mathcal{M}_{\mathrm{cyl}}(V;\mathfrak{A}'_+)  \cup \bigsqcup_{h>0}\mathcal{R}_h\; .
\end{equation}
Let us denote by $F(V)$ the set of $\sigma(V,V')$-closed subspaces of $V$ with finite codimension. For any
$\Phi\in F(V)$, $V'\supset \Phi^{\circ}\cong (V/\Phi)'$ is its finite dimensional polar. Finally, given a
locally compact real vector space $L$ let us denote by $C_0(L)'_{C_0(L)}$ the continuous dual of the space
of its compactly supported continuous functions endowed with the $\sigma\bigl(C_0(L)', C_0(L)\bigr)$
topology. It is possible to define a map $\mathscr{H}$ from $S(V,X,\mathfrak{A})$ to $\prod_{\Phi\in
  F(V)}\bigl(C_0(V/\Phi)'_{C_0(V/\Phi)}\bigr)^{\mathfrak{A}_+}$ in the following way:
\begin{equation*}
  \mathscr{H}(s)=\left\{
    \begin{aligned}
      &\Bigl(a\mapsto\int_{V/\Phi}^{}\;\cdot\;   \mathrm{d}\mu_{\Phi,a}\Bigr)_{\Phi\in F(V)} &\quad s=M\in \mathcal{M}_{\mathrm{cyl}}(V;\mathfrak{A}'_+)\\
      &\Bigl(a\mapsto H^{(\Phi^{\circ})}_{\omega_h}(a)(\cdot)\Bigr)_{\Phi\in F(V)}&\quad s=(\omega_h,h)\in \mathcal{R}_h\times \mathds{R}^+_*
    \end{aligned}\right . \; ,
\end{equation*}
where $H^{(\Phi^{\circ})}_{\omega_h}(a)$ is defined by~\eqref{eq:123}, and
$\mathds{R}^+_*=\mathds{R}^+\smallsetminus \{0\}$. The product topology $\mathfrak{P}_{F(V)}$ on
\begin{equation*}
  \prod_{\Phi\in F(V)}\bigl(C_0(V/\Phi)'_{C_0(V/\Phi)}\bigr)^{\mathfrak{A}_+}_{\mathrm{s}}\; ,
\end{equation*}
where each $\bigl(C_0(V/\Phi)'_{C_0(V/\Phi)}\bigr)^{\mathfrak{A}_+}_{\mathrm{s}}$ is endowed with the
topology of simple convergence, induces a topology on $S(V,X,\mathfrak{A})$: the preimage topology
$\mathfrak{P}=\mathscr{H}^{-1}\mathfrak{P}_{F(V)}$.

On the other hand, by Bochner's theorem~\ref{thm:103} and Proposition~\ref{prop:27}, there is a map
$\mathscr{F}:S(V,X,\mathfrak{A})\to (\mathfrak{A}'_{\mathfrak{A}})^X$ defined as follows:
\begin{equation*}
  \mathscr{F}(s)=\left\{
    \begin{aligned}
      &\hat{M}(\cdot )&\quad s=M\in \mathcal{M}_{\mathrm{cyl}}(V;\mathfrak{A}'_+)\\
      &\mathcal{G}_{\omega_h}(\cdot )&\quad s=(\omega_h,h)\in \mathcal{R}_h\times \mathds{R}^+_*
    \end{aligned}\right . \; .
\end{equation*}
The topology $\mathfrak{T}_{\mathrm{s}}$ of simple convergence on $(\mathfrak{A}'_{\mathfrak{A}})^X$
induces therefore the topology $\mathfrak{T}=\mathscr{F}^{-1}\mathfrak{T}_{\mathrm{s}}$ on
$S(V,X,\mathfrak{A})$. Finally, let us denote by $\mathfrak{P}\vee \mathfrak{T}$ the join topology,
\emph{i.e.}\ the coarsest topology on $S(V,X,\mathfrak{A})$ that is finer than both $\mathfrak{P}$ and
$\mathfrak{T}$. Hence \cref{thm:14} can be reformulated as a convergence result, in suitable topologies,
for semiclassical quantum states, as $h\to 0$. In order to do that, we have only to take into account the
projective structure induced by $F(V)$.

Let $F_{\varsigma}(X)$ be the set of all finite dimensional symplectic subspaces of $X$. First of all, let
us remark that given any finite dimensional vector subspace $R$ of $X$, there is always an $X_R\in
F_{\varsigma}(X)$ such that $R$ is a (closed) subspace of $X_R$. Therefore if $g\in
(\mathfrak{A}'_{\mathfrak{A}})^X$ is continuous when restricted to every symplectic subspace in
$F_{\varsigma}(X)$, then it is continuous when restricted to every finite dimensional subspace. In
addition, let $V$ be a topological real vector space in compatible duality with $X$; then the set
$F^{\circ}_{\varsigma}(X)=\bigl\{\Xi^{\circ}, \Xi\in F_{\varsigma}(X)\bigr\}$, where the polar is taken
with respect to the duality between $X$ and $V$, is the subset of $F(V)$ consisting of the weakly closed
subspaces with \emph{even} codimension. If we partially order the set $F(V)$ by $\supset $, then
$F^{\circ}_{\varsigma}(X)$ is a \emph{cofinal} subset of $F(V)$. In \cref{sec:compactness}, given a
semiclassical quantum state $(\omega_{h_{\beta}})_{\beta\in B}$ on $(\mathfrak{W}_{h_{\beta}})_{\beta\in
  B}$, we defined the family of relatively compact and nonempty sets
\begin{equation*}
  \bigl(H_B(\Phi)\bigr)_{\Phi\in F^{\circ}_{\varsigma}(X) }=\Bigl(\prod_{a\in \mathfrak{A}_+} \bigl\{H^{(\Phi)}_{\omega_{h_{\beta}}}(a),\beta\in B\bigr\}\Bigr)_{\Phi\in F^{\circ}_{\varsigma}(X)}\; ,
\end{equation*}
by \cref{eq:123,eq:124}, with $X_{\Phi}=\Phi^{\circ}$; and since the definition of
$H^{(\Phi)}_{\omega_{h_{\beta}}}(a)$ can be extended to any $\Phi\in F(V)$, we can also define the family
\begin{equation*}
  \bigl(H_B(\Phi)\bigr)_{\Phi\in F(V) }=\Bigl(\prod_{a\in \mathfrak{A}_+} \bigl\{H^{(\Phi)}_{\omega_{h_{\beta}}}(a),\beta\in B\bigr\}\Bigr)_{\Phi\in F(V)}\; .
\end{equation*}
In addition, if $R\subset X$ is a finite dimensional vector subspace of $X$, we define the subalgebra
$\mathbb{W}_h(R)\subset \mathbb{W}_h(X,\varsigma)$ by
\begin{equation*}
  \mathbb{W}_h(R)=\textrm{C*}(\bigl\{W_h(x),x\in R\bigr\})\; .
\end{equation*}
By this definition, it follows that for any $a\in \mathfrak{A}_+$ and $\Phi\in F(V)$,
\begin{equation*}
  H^{(\Phi)}_{\omega_{h_{\beta}}}(a)=H_{\omega_{h_{\beta}}\bigr\rvert_{\mathbb{W}_{h_{\beta}}(\Phi^{\circ})\otimes _{\gamma_{h_{\beta}}}\mathfrak{A}}}(a)\; ,
\end{equation*}
where the latter is defined by \cref{eq:122}. For any $\Phi\supset \Psi \in F(V)$, let us now introduce
the (continuous) maps $\pi_{\Phi\Psi}:H_B(\Psi)\to H_B(\Phi)$:
\begin{equation*}
  \pi_{\Phi\Psi}\Bigl(a\mapsto H^{(\Psi)}_{\omega_{h_{\beta}}}(a)\Bigr)=a\mapsto H^{(\Phi)}_{\omega_{h_{\beta}}}(a) \; ,
\end{equation*}
where the definition is justified by the fact that $\Phi^{\circ}$ is injected canonically in
$\Psi^{\circ}$ by the map $^{\mathrm{t}}p_{\Psi}\circ \, ^{\mathrm{t}}p_{\Phi\Psi}\circ\,
^{\mathrm{t}}p_{\Phi}^{-1}$, where $p_{\Phi}:V\to V/\Phi$ is the canonical map and $p_{\Phi\Psi}:
V/\Psi\to V/\Phi$ is the map obtained from the identity passing to the quotients. It is easy to see that
$p_{\Phi\Omega}=p_{\Phi\Psi}\circ p_{\Psi\Omega}$ for any $\Phi\supset \Psi\supset \Omega\in
F(V)$. Therefore the families
\begin{equation*}
  \bigl(H_B(\Phi),\pi_{\Phi\Psi}\bigr)_{\Phi\supset \Psi\in F^{\circ}_{\varsigma}(X) }\; ,\; \bigl(H_B(\Phi),\pi_{\Phi\Psi}\bigr)_{\Phi\supset \Psi\in F(V) }
\end{equation*}
are projective systems with cofinal index set and thus they share the same projective limit. In addition,
the limit is nonempty and relatively compact since $H_B(\Phi)$ is nonempty and relatively compact for any
$\Phi\in F^{\circ}_{\varsigma}(X)$. Therefore for $F=F_{\varsigma}(X)$, the set $H_B$ of \cref{lemma:12}
can be restricted to the projective limit
\begin{equation*}
  H_B=\varprojlim H_B(\Phi)\; .
\end{equation*}
It then follows by \cref{thm:14} that the $\mathfrak{P}$-cluster points of semiclassical quantum states
are \emph{projective families} of Radon measures on the finite dimensional quotients of $V$, and therefore
cylindrical
measures. 
A convergent net $\omega_{h_b}\underset{h_b\to 0}{\overset{\mathfrak{P}}{\longrightarrow}}M $, that
satisfies one of the following four equivalent conditions has no loss of mass:
\begin{align}
  \label{eq:43}\tag{i}
  &\omega_{h_b}\underset{h_b\to 0}{\overset{\mathfrak{P}\vee \mathfrak{T}}{\longrightarrow}} M &; \\
  \label{eq:44}\tag{ii}
  &\begin{aligned}
    &\lim_{\varepsilon\to 0}\lim_{b\in \underline{B}}\omega_{h_b}\bigl(W_{h_b}(x)\otimes a\!-\!\mathrm{Op}_{\frac{1}{2}}^{h_b}(^{\exists }\psi_{\varepsilon,\Phi}(\cdot )\, e^{2ix(\cdot )})\otimes a\bigr)\!=\!0&\bigl(\forall \Phi\in F_{\varsigma}^{\circ}(X)\,,\\ & & \mspace{-57mu}\forall x\in \Phi^{\circ}\,,\,\forall a\in \mathfrak{A}_+\bigr)
  \end{aligned}
  &;\\
  \label{eq:45}\tag{iii}
  &\begin{aligned}
      &\lim_{\varepsilon\to 0}\lim_{b\in \underline{B}}\omega_{h_b}\bigl(a-\mathrm{Op}_{\frac{1}{2}}^{h_b}(^{\exists }\psi_{\varepsilon,\Phi})\otimes  a\bigr)=0&\mspace{54mu}\bigl(\forall \Phi\in F_{\varsigma}^{\circ}(X)\,,\,\forall a\in \mathfrak{A}_+\bigr)
  \end{aligned}
  &;\\
  \label{eq:46}\tag{iv}
  &\lim_{b\in \underline{B}}\omega_{h_b}\bigl(W_{h_b}(0)\bigr)=M(0)&.
\end{align}
From a practical point of view, it is easier to study the semiclassical behavior of quantum states with no
loss of mass, since it is very convenient to characterize the limit measure by means of the noncommutative
and commutative Fourier transforms.

\begin{thm}\label{thm:107}
  Let $(X,\varsigma)\in \mathbf{Symp}_{\mathds{R}}$, $\mathfrak{A}\in \mathbf{C^{*}alg}$, and
  $(\mathfrak{W}_h)_{h\geq 0}$ the corresponding Weyl deformation \eqref{eq:112}. In addition, let $V$ be
  any topological space in compatible duality with $X$. Consider a semiclassical state
  $(\omega_{h_{\beta}})_{\beta\in B}$ on $(\mathfrak{W}_{h_{\beta}})_{\beta\in B}$. Then there exists a
  subnet $(\omega_{h_b})_{b\in \underline{B}}$ such that
  \begin{equation}
    \label{eq:1013}
    \omega_{h_b}\underset{h_b\to 0}{\overset{\mathfrak{P}}{\longrightarrow}} M\; ,
  \end{equation}
  where $M\in \mathcal{M}_{\mathrm{cyl}}(V;\mathfrak{A}'_+)$ is a \emph{cylindrical Wigner measure}. In
  other words, if we denote by $\mathscr{W}\bigl(\omega_{h_{\beta}},\beta\in B\bigr)$ the set of cluster
  points of $(\omega_{h_{\beta}})_{\beta\in B}\subset S(V,X,\mathfrak{A})_{\mathfrak{P}}$, then $s\in
  \mathscr{W}\bigl(\omega_{h_{\beta}},\beta\in B\bigr)$ implies $s\in
  \mathcal{M}_{\mathrm{cyl}}(V;\mathfrak{A}'_+)$. If in addition $(\omega_{h_b})_{b\in \underline{B}}$ has
  no loss of mass then
  \begin{equation}
    \label{eq:1030}
    \omega_{h_b}\underset{h_b\to 0}{\overset{\mathfrak{P}\vee \mathfrak{T}}{\longrightarrow}} M\; .
  \end{equation}
\end{thm}

The cylindrical measures are not countably additive Borel measures, even if each countably additive Borel
measure defines a cylindrical measure. A cylindrical measure on the topological space $V$ is a projective
system of Borel measures on finite dimensional vector spaces, indexed by the weakly closed subspaces of
$V$ of finite codimension. Let us denote by $F(V)$ the set of $\sigma(V,V')$-closed subspaces with finite
codimension. The set $F(V)$ is identified uniquely by the finite dimensional subspaces of $V'$. In fact,
$\Phi\in F(V)$ iff there exists a finite dimensional subspace $F\subset V'$ such that $F^{\circ}=\Phi$,
where $F^{\circ}$ is the polar (orthogonal) of $F$ (with respect to the canonical duality between $V$ and
$V'$). The \emph{if} part is proved as follows: the bipolar $F^{\circ \circ}$ is isomorphic to $V/\Phi$,
and it is also the closure of $F$ with respect to the $\sigma(V',V)$ topology; since the duality between
$V$ and $V'$ separates points in $V'$, it follows that the finite dimensional subspace $F$ is closed with
respect to $\sigma(V',V)$, and therefore $F=F^{\circ \circ}$, thus proving $\Phi\in F(V)$. The \emph{only
  if} part is proved as follows: let $\Phi\in F(V)$, and choose $F=\Phi^{\circ}$; then
$F^{\circ}=\Phi^{\circ \circ}$ is the closure of $\Phi$ with respect to the $\sigma(V,V')$ topology, and
therefore $\Phi^{\circ \circ}=\Phi$ since $\Phi\in F(V)$. This further motivates the fact that it is the
phase space $X$ (of test functions) that determines the cylindrical measures, notwithstanding the
possibility of choosing the measurable space of classical fields $V$ in many different ways.

To sum up, the semiclassical quantum states always have classical counterparts, the cylindrical Wigner
measures, equivalently acting on some suitably defined topological space of fields. On the other hand,
families of non-regular states cannot be identified, in general, with cylindrical measures, because their
generating map could fail to be ultraweakly continuous on every finite dimensional subspace. Therefore
non-regular states are not suitable for a classical effective description. From a physical standpoint this
is not unreasonable, since non-regular states appear in relation to typically quantum physical behaviors,
such as infrared divergence \citep[see, \emph{e.g.},][]{acerbi1993rmp,acerbi1993jmp%
}.

\subsection{Phase space and Wigner measures}
\label{sec:phase-spaces-physics}

The symplectic space $(X,\varsigma)$ has a natural interpretation in physics, especially when finite
dimensional, as the phase space of a given theory. The symplectic structure is extremely useful to study
dynamical properties (Hamiltonian flows). The usual Wigner measures for finite-dimensional Heisenberg
groups can be interpreted straightforwardly as (probability) measures on the phase space, since
$(X,\varsigma)$ is in compatible duality with itself whenever $\dim X <\infty$. For infinite-dimensional
Heisenberg groups, it is possible only in few favorable cases, that may be called \emph{admissible phase
  spaces}. Explicitly, Wigner measures are cylindrical measures on the phase space if $X$ is a topological
space in compatible duality with itself. Interestingly enough, this allows to consider the geometric phase
spaces constructed from reflexive spaces (seen as topological manifolds). Given a \emph{reflexive} locally
convex space $\Sigma$, its cotangent bundle $\mathrm{T}^{*}\Sigma\cong \Sigma\oplus \Sigma'$ is a
symplectic space with the canonical symplectic form. Since $\Sigma$ is reflexive, $\Sigma\oplus \Sigma'$
is in compatible duality with itself, and therefore $\mathrm{T}^{*}\Sigma$ is an admissible phase space.

As already discussed in \cref{sec:phys-appl}, in quantum field theories the phase space $(X,\varsigma)$ is
usually interpreted as a space of \emph{test functions}, and the space of classical fields is in duality
with it (it is therefore a space of distributions). The Wigner measures act on the latter, and this is
physically reasonable, since the measures should describe the configuration of classical fields and not of
test functions. The symplectic structure of test functions is carried out on the space of classical
fields, if the latter is taken to be $X^{*}_X$: for any $x, y\in X$, $\varsigma(x,\cdot )\in X^{*}$ and
$\varsigma(x,\cdot )\neq \varsigma(y,\cdot )$ if $x\neq y$ (by the non-degeneracy of $\varsigma$), and
thus $X\overset{\mathrm{s}}{\hookrightarrow} X^{*}_X$, and there is a symplectic form
$\tilde{\varsigma}(\cdot ,\cdot )$ on $\mathrm{s}(X)$ defined by $\tilde{\varsigma}(\cdot ,\cdot
)=\varsigma(\mathrm{s}^{-1}\,\cdot \,,\mathrm{s}^{-1}\,\cdot \,)$. In order to define an Hamiltonian
dynamics on the space of classical fields it is therefore necessary to either restrict to measures
concentrated on $\mathrm{s}(X)$ (but as discussed in \cref{sec:surj-union-sets}, this point of view may be
too restrictive), or to extend $\tilde{\varsigma}$ to the whole $X^{*}_X$.

\subsection{Cylindrical Wigner measures as Radon measures}
\label{sec:wigner-measures-x}


Cylindrical measures are not completely satisfactory from a practical standpoint. Integration of functions
with respect to cylindrical measures is possible only if the functions are cylindrical as well. It is
possible to interpret cylindrical measures on $V$ as Radon measures only on a space that is ``bigger''
than $V$. This is done exploiting Prokhorov's tightness of the projective family of measures
\citep[][]{schwartz1973tirsm
}. For any $\Phi\in F(V)$, let us consider the \v{C}ech compactification $\overline{V/\Phi}$ of
$V/\Phi$. From the canonical injection $j_{\Phi}:V/\Phi\to \overline{V/\Phi}$ it is possible to construct
a canonical injection, let us call it $j$, of $V$ into the product space
\begin{equation*}
  \overline{V}=  \prod_{\Phi\in F(V)}\overline{V/\Phi}\; .
\end{equation*}
We endow $\overline{V}$ with the product topology. The push-forward $\overline{M}= j\, _{*}\, M$ of any
$M\in \mathcal{M}_{\mathrm{cyl}}(V;\mathfrak{A}'_+)$ is defined as the family
$\overline{M}=(j_{\Phi}\,_{*}\, \mu_{\Phi})_{\Phi\in F(V)}$. It is a tight projective family of measures
on $\overline{V}$, and therefore it originates from a finite Borel Radon measure $\mu_{\overline{M}}$. It
follows that the set of cylindrical measures on $V$ can be identified with a subset of finite Radon
measures on $\overline{V}$. If $M$ does not originate from a Radon measure on $V$ then
$\mu_{\overline{M}}(j(V))\neq m$ (where $m=\mu(\overline{V})$ is the total mass of the measure), and there
are cylindrical measures for which $\mu_{\overline{M}}(j(V))= 0$.

As it will be proved in the next section, every cylindrical measure, on any space $V$ in compatible
duality with $X$, is the Wigner measure of at least one generalized sequence of regular quantum
states. Therefore cylindrical measures that are not Radon measures play an important role in the
semiclassical description. In other words, the spaces $\overline{V}$, rather than $V$ or $X$, are the most
suitable candidates to accommodate a complete classical description of field theories (as classical
probability theories). The non-uniqueness of $\overline{V}$ may seem on one hand disappointing, on the
other hand it allows for additional freedom in the choice of classical effective description. We have
already highlighted that there seems to be one ``preferred'' choice of $V$, \emph{i.e.}\ choosing it to be
the algebraic dual $X^{*}$ of $X$, with the $\sigma(X^{*},X)$ topology. A classical result \citep[][II.54,
Proposition 9]{bourbaki1981evt1-5
} lets us identify $X^{*}_X$ with the Hausdorff completion of $\bigl(V,\sigma(V,X)\bigr)$ for any $V$ and
$X$ in a duality that separates points of both $V$ and $X$. Finally, for any $V$ in a compatible duality
with $X$ that separates points in $V$, there is a natural bijection
\citep[][]{schwartz1973tirsm
} (equivalently yielded by Bochner's theorem)
\begin{equation}
  \label{eq:1014}
  \mathcal{M}_{\mathrm{cyl}}(V;\mathfrak{A}'_+) \cong \mathcal{M}_{\mathrm{cyl}}(X^{*}_X;\mathfrak{A}'_+)\; .
\end{equation}
The bijection~\eqref{eq:1014} becomes of particular importance whenever $X$ is a second countable
topological space. In fact, in such case any cylindrical measure originates from a Radon measure:
\begin{equation*}
  \mathcal{M}_{\mathrm{cyl}}(X^{*}_X;\mathfrak{A}'_+)\cong \mathcal{M}_{\mathrm{rad}}(X^{*}_X;\mathfrak{A}'_+)\; ,
\end{equation*}
where the latter is the set of all finite Borel Radon measures. The Wigner measures are then exactly the
finite Borel Radon measures on $X^{*}_X$. Let us remark that also for a second countable $X$ Wigner
measures can be seen as Radon measures on $\overline{X^{*}_X}$, however they are precisely the set of
measures that are concentrated as Borel Radon measures on $X^{*}_X$.

\subsection{Characterization of the set of all Wigner measures}
\label{sec:surj-union-sets}

As anticipated in \cref{thm:5,sec:wigner-measures-x}, any cylindrical measure is the Wigner measure of
some semiclassical quantum state. To prove this very interesting result, let us define the set
$\mathscr{W}$ of all possible Wigner measures associated to the Weyl deformation \emph{with maximal cross
  norm}. Let us denote the aforementioned Weyl deformation by
\begin{equation*}
  (\mathfrak{W}_h^{\mathrm{m}})_{h>0}=\bigl(\mathbb{W}_h(X,\varsigma)\otimes_{\gamma_{\mathrm{max}}}\mathfrak{A}\bigr)_{h>0}\; .
\end{equation*}
Then, using the notation introduced in \cref{thm:107}, the set of all Wigner measures is defined as
\begin{equation*}
  \mathscr{W}=\bigcup_{\underline{h}\in \mathds{R}^+_*}\; \bigcup_{(\omega_{h})_{h\in (0,\underline{h})}\in  2^{S(V,X,\mathfrak{A})}} \mathscr{W}\bigl(\omega_{h},  h\in (0,\underline{h})\bigr)\; .
\end{equation*}
The aim is to prove, for any $V$ in compatible duality with $X$,
$\mathscr{W}=\mathcal{M}_{\mathrm{cyl}}(V;\mathfrak{A}'_+)$. By \cref{thm:107}, it suffices to prove that
$\mathscr{W}\supset\mathcal{M}_{\mathrm{cyl}}(V;\mathfrak{A}'_+)$, \emph{i.e.}\ that for any cylindrical
measure $M$ there exists a semiclassical quantum state $(\omega^{(M)}_{h})_{h\in (0,\underline{h})}$ such
that $\omega^{(M)}_{h}\underset{h\to 0}{\overset{\mathfrak{P}\vee \mathfrak{T}}{\longrightarrow}} M$ on
$S(V,X,\mathfrak{A})$.

We make use of the following result of standard semiclassical analysis. Its proof relies on the fact that
squeezed coherent states on $L^2(\mathds{R}^d)$ converge to measures concentrated on a point of
$\mathds{R}^{2d}$; and that linear combinations of point measures are dense in the space of finite
measures $\mathcal{M}_{\mathrm{rad}}(\mathds{R}^{2d})_{\mathds{C}}$, endowed with the weak topology
\citep[see][for additional details]{falconi2016bpmas
}. The extension to general finite dimensional symplectic spaces $(R,\varsigma)$ and to tensor products
$\mathbb{W}_{h}(R,\varsigma)\otimes_{\gamma_{\mathrm{max}}}\mathfrak{A}$ does not present difficulties.
\begin{lemma}\label{lemma:106}
  Let $(R_0,\varsigma)\in \mathbf{Symp}_{\mathds{R}}$. For any $\mu\in
  \mathcal{M}_{\mathrm{rad}}(R'_0;\mathfrak{A}'_+)$, there exists a semiclassical state
  $(\tilde{\varrho}_{h})_{h\in (0,1)}$ such that for any $h\in (0,1)$, $\tilde{\varrho}_{h}\in
  \bigl(\mathbb{W}_{h}(R_0,\varsigma)\otimes_{\gamma_{\mathrm{max}}}\mathfrak{A} \bigr)'_+$ is normal,
  \begin{equation*}
    \lVert \varrho_h  \rVert_{R_0}^{}:=\tilde{\varrho}_h\bigl(W_{h}(0)\otimes 1\bigr)=\lVert \hat{M}(0)  \rVert_{\mathfrak{A}'}^{}\; ,
  \end{equation*}
  and on $S(R_{0}',R_0,\mathfrak{A})$
  \begin{equation*}
    \tilde{\varrho}_{h} \underset{h\to0}{\overset{\mathfrak{P}\vee \mathfrak{T}}{\longrightarrow}} \mu\; .
  \end{equation*}
\end{lemma}
It is then sufficient to combine Lemma~\ref{lemma:106} with a projective argument, obtaining the following
theorem for any topological vector space $V$ in compatible duality with the real symplectic space
$(X,\varsigma)$.
\begin{thm}\label{thm:1016}
  $\mathscr{W}=\mathcal{M}_{\mathrm{cyl}}(V;\mathfrak{A}'_+):$ $\!\,^{\exists }\omega_h^{(M)}
  \underset{h\to0}{\overset{\mathfrak{P}\vee \mathfrak{T}}{\longrightarrow}}M\; \bigl(\forall M\in
  \mathcal{M}_{\mathrm{cyl}}(V;\mathfrak{A}'_+)\bigr)$.
\end{thm}
\begin{proof}
  Let us start with the simpler case where $X=\bigcup_{m\in \mathds{N}}X_m$, for a countable directed
  sequence of finite dimensional symplectic subspaces $(X_m)_{m\in \mathds{N}}$, ordered by inclusion. The
  family $(X_m)_{m\in \mathds{N}}$, together with the identities $i_{mn}:X_m\to X_n$, $m\leq n$ is an
  inductive family of finite dimensional vector spaces. Therefore $\varinjlim X_m= X$, and
  \begin{equation}
    \label{eq:1015}
    (\mathfrak{A}_{\mathfrak{A}}')^X=\varprojlim\, (\mathfrak{A}_{\mathfrak{A}}')^{X_m}\; .
  \end{equation}
  Now given a cylindrical measure $M=(\mu_{\Phi})_{\Phi\in F(V)}$, the subfamily
  $(\mu_{X_m^{\circ}})_{m\in \mathds{N}}$ is sufficient to characterize it completely (the polar is
  intended with respect to the duality between $V$ and $X$). In fact, $\hat{M}: X\to
  \mathfrak{A}_{\mathfrak{A}}'$ can be characterized uniquely by the projective family of restrictions
  $(\hat{M}\rvert_{X_m})_{m\in \mathds{N}}$.

  Consider $X_0$. We define the following set of generating maps, corresponding to \emph{regular} states
  on the maximal cross product algebra
  $\mathbb{W}_{h}(X_0,\varsigma)\otimes_{\gamma_{\mathrm{max}}}\mathfrak{A}$:
  \begin{equation*}
    G(\mu_{X_0^{\circ}})=\Bigl\{\mathcal{G}_{\varrho_h}\;, \; h\in (0,1)\;,\; \varrho_h\underset{h\to0}{\overset{\mathfrak{P}\vee
        \mathfrak{T}}{\longrightarrow}}\mu_{X_0^{\circ}}\Bigr\}\; . 
  \end{equation*}
  From $G(\mu_{X_0^{\circ}})$, construct (recursively) the compatible projective family of generating maps
  of regular states $\bigl(G_m,\imath_{mn}\bigr)_{m\leq n\in \mathds{N}}$, defined as follows.
  $G_0=G(\mu_{X_0^{\circ}})$, and for any $m>0$
  \begin{equation*}
    G_m=\Bigl\{\mathcal{G}_{\varrho_h}\;,h\in (0,1)\;, \mathcal{G}_{\varrho_h}\bigr\rvert_{X_{m-1}}\in \mathcal{G}_{m-1}\Bigr\}\; .
  \end{equation*}
  The \emph{surjective} maps $\imath_{mn}:G_n\to G_m$\; , $m\leq n$, are defined by
  $\imath_{mn}(\mathcal{G}_{\varrho_h})=\mathcal{G}_{\varrho_h}\circ\, i_{mn}$, where $i_{mn}$ are the
  identities defined above for the inductive family $(X_m)_{m\in \mathds{N}}$. It is easy to prove that
  for any $m\in \mathds{N}$, there exists $(\varrho_h^{(m)})_{h\in (0,1)}$,
  $\bigl(\mathcal{G}_{\varrho_h^{(m)}}\bigr)_{h\in (0,1)}\subset G_m$, such that
  $\varrho_h^{(m)}\underset{h\to0}{\overset{\mathfrak{P}\vee \mathfrak{T}}{\longrightarrow}}
  \mu_{X_m^{\circ}}$. For $m=0$ it is true by definition, and for any $m>0$ it follows from the fact that
  for any $h\in (0,1)$, $\mathcal{G}_{\tilde{\varrho}_h}\bigr\rvert_{X_{m-1}}\in G_{m-1}$, where
  $(\tilde{\varrho}_h)_{h\in (0,1)}$ is the convergent net of \cref{lemma:106}. Therefore, \emph{a
    fortiori}, all $G_m$ are non-empty. By \cref{eq:1015} and \citep[][III.58 Proposition
  5]{bourbaki1970ens
  }, the projective limit $G=\varprojlim G_m$ is a non-empty set of nets of generating maps
  $(\mathcal{G}_{\omega_h})_{h\in (0,1)}$ of regular states on $(\mathfrak{W}_h^{\mathrm{m}})_{h\in
    (0,1)}$ (with maximal cross norm on the tensor products). In addition, by the reasoning above, there
  exists $(\mathcal{G}_{\omega_h^{(M)}})_{h\in (0,1)}\subset G$ such that
  $\omega_h^{(M)}\underset{h\to0}{\overset{\mathfrak{P}\vee \mathfrak{T}}{\longrightarrow}} M$.

  In the general case, it is possible to follow the same guidelines. Given a finite dimensional symplectic
  subspace $R_0\subset X$, it is always possible to find a family $(X_{\lambda})_{\lambda\in J(V)}$ of
  finite dimensional symplectic subspaces -- where $J(V)$ is some directed set with partial order $\leq $
  and least element $0$ -- such that $X_0=R_0$, $X_{\eta}\subset X_{\lambda}$ for any $\eta\leq \lambda\in
  J(V)$, and $\bigcup_{\lambda\in J(V)}X_{\lambda}=X$. Therefore $(X_{\lambda})_{\lambda\in J(V)}$ is an
  inductive family, and the following identity holds:
  \begin{equation}
    \label{eq:1016}
    (\mathfrak{A}_{\mathfrak{A}}')_{\mathrm{s}}^X=\varprojlim\, (\mathfrak{A}_{\mathfrak{A}}')_{\mathrm{s}}^{X_{\lambda}}\; .
  \end{equation}
  In addition, any cylindrical measure $M=(\mu_{\Phi})_{\Phi\in F(V)}$ can be equivalently characterized
  by the family of restricted Fourier transforms $(\hat{M}\rvert_{X_{\lambda}})_{\lambda\in
    J(V)}$. Similarly to the countable case, we define the set of generating maps of regular states
  \begin{equation*}
    G_0=\Bigl\{\mathcal{G}_{\varrho_h}\;, \; h\in (0,1)\;,\; \varrho_h\underset{h\to0}{\overset{\mathfrak{P}\vee
        \mathfrak{T}}{\longrightarrow}}\mu_{X_0'}\;, \lVert \varrho_h  \rVert_{X_0}^{}=\lVert \hat{M}(0)  \rVert_{\mathfrak{A}'}^{}\Bigr\}\; ,
  \end{equation*}
  where $\lVert \,\cdot \,_h \rVert_{X_0}^{}$ is the norm of the Banach space
  $\bigl(\mathbb{W}_h(X_0,\varsigma)\otimes_{\gamma_{\mathrm{max}}}\mathfrak{A}\bigr)'$. The set
  $G_0\subset (\mathfrak{A}'_{\mathfrak{A}})^{X_0}_{\mathrm{s}}$ is relatively compact and non-empty, by
  Banach-Alaoglu's theorem and \cref{lemma:106}. Recursively, we define the sets $G_{\lambda}\subset
  (\mathfrak{A}'_{\mathfrak{A}})^{X_{\lambda}}_{\mathrm{s}}$, $\lambda\in J(V)$, by
  \begin{equation*}
    G_{\lambda}=\Bigl\{\mathcal{G}_{\varrho_h}\;,h\in (0,1)\;, \forall \eta\leq \lambda\;,\; \mathcal{G}_{\varrho_h}\bigr\rvert_{X_{\eta}}\in G_{\eta}\Bigr\}\; .
  \end{equation*}
  The sets $G_{\lambda}$ are all relatively compact and non-empty as well. In addition, the projections
  $\imath_{\lambda\eta}:G_{\eta}\to G_{\lambda}$, $\lambda\leq \eta$ defined by
  $\imath_{\lambda\eta}(\mathcal{G}_{\varrho_h})=\mathcal{G}_{\varrho_h}\circ\, i_{X_{\lambda}X_{\eta}}$
  are continuous. Therefore the projective system $(G_{\lambda},\imath_{\lambda\eta})_{\lambda\leq \eta\in
    J(V)}$ has a non-empty limit \citep[][]{bourbaki1970ens
  }, and $G=\varprojlim G_{\lambda}$ is a set of generating maps $(\mathcal{G}_{\omega_h})_{h\in (0,1)}$
  of regular states on $(\mathfrak{W}_h^{\mathrm{m}})_{h\in (0,1)}$. In addition, as in the countable
  case, there exists $(\mathcal{G}_{\omega_h^{(M)}})_{h\in (0,1)}\subset G$ such that
  $\omega_h^{(M)}\underset{h\to0}{\overset{\mathfrak{P}\vee \mathfrak{T}}{\longrightarrow}} M$.
\end{proof}

\section{Push-forward and convolution of Wigner measures}
\label{sec:group-homom-push}

In this section we study transformations on the cylindrical Wigner measures induced by transformations of
the regular quantum states, and prove the results outlined in \cref{sec:maps-conv-prod}. Throughout the
section, we set $\mathfrak{A}=\mathds{C}$ for simplicity.

\subsection{Group homomorphisms}
\label{sec:centr-heis-group}
Endomorphisms of $\mathbb{W}_h(X,\varsigma)$ -- and perhaps more interestingly endomorphisms of
$\pi\bigl(\mathbb{W}_h(X,\varsigma)\bigr)''$ (its bicommutant in some irreducible representation) -- play
a crucial role in defining quantum dynamical systems, as briefly discussed in
\cref{sec:semicl-time-slice,sec:semicl-kms-ground}. Therefore, it is interesting to study the
semiclassical action induced by *-homomorphisms in Weyl algebras. A systematic study of *-homomorphisms
would require, at least to some extent, the development of pseudodifferential calculus for infinite
dimensional phase spaces, and it is beyond the scope of this work. We restrict our attention to
*-homomorphisms induced by a class of central Heisenberg group homomorphisms (field morphisms). In this
way it is possible to characterize the induced semiclassical action on cylindrical Wigner measures in a
natural way.

Let $(X,\varsigma),(Y,\tau)\in \mathbf{Symp}_{\mathds{R}}$. For any $\tilde{f}\in X^{Y\times \mathds{R}}$
and $\tilde{F}\in \mathds{R}^{Y\times \mathds{R}}$,
\begin{equation*}
  (\tilde{f} ,\tilde{F})\in \mathrm{Hom}_{\mathrm{gr}}\bigl(\mathbb{H}(Y,\tau),\mathbb{H}(X,\varsigma)\bigr)
\end{equation*}
iff for any $\varphi,\psi\in Y\times \mathds{R}$:
\begin{align*}
  \tilde{f}&(\varphi\cdot \psi)=\tilde{f}(\varphi)+\tilde{f}(\psi)\\
  \tilde{F}&(\varphi\cdot \psi)=\tilde{F}(\varphi)+\tilde{F}(\psi)-\varsigma(\tilde{f}(\varphi),\tilde{f}(\psi))\; .
\end{align*}
In addition, $(\tilde{f},\tilde{F})$ is central if for any $t\in \mathds{R}$, $\tilde{f}(0,t)=0$. Among
the central homomorphisms, we restrict our attention to the ones of the following form: there exists $f\in
X^{Y}$ 
such that $\tilde{f}(y,t)=f(y)$.
Let us denote such homomorphisms by $\tilde{\mathfrak{F}}$. On the Weyl C*-algebra representation
$\mathbb{W}_h(Y,\tau)$ of the Heisenberg group, the induced action of $\tilde{\mathfrak{F}}$ is then
\begin{align*}
  & \tilde{\mathfrak{F}}\bigl(e^{it}W_h(y)\bigr)=e^{i \tilde{F}(y,t)}W_h(f(y)) \, ,\quad \forall y\in Y\;,\;\forall t\in \mathds{R}\; .
\end{align*}
As it will become clear in the following, it is useful to consider $h$-dependent homomorphisms of the form
$\tilde{\mathfrak{F}}_h=(f,\tilde{F}_h)$. In order to obtain a linear transformation on the Weyl
C*-algebra, we ``forget'' the $t$-dependence of $\tilde{\mathfrak{F}}_h$ setting
$F_h(x)=\tilde{F}_h(x,0)$, and consider the associated homomorphism $\mathfrak{F}_h=(f,F_h)$. Since
$\mathfrak{F}_h$ does not act on scalars anymore, 
it can be extended to a (isometric) *-homomorphism
\begin{equation*}
  \mathfrak{F}_h\in \,  ^{*}\mathrm{Hom}\bigl(\mathbb{W}_h(Y,\tau),\mathbb{W}_h(X,\varsigma) \bigr)\; ,
\end{equation*}
with the following conditions on $f$ and $F_h$:
\begin{align*}
  &f(y+y')=f(y)+f(y')&(\forall y,y'\in Y)\\
  &F_h(y+y')=F_h(y)+F_h(y')-ih \bigl(\varsigma(f(y),f(y')-\tau(y,y')\bigr)&(\forall y,y'\in Y)
\end{align*}
By duality, the adjoint map $^{\mathrm{t}}\mathfrak{F}_h\in
\mathcal{B}\bigl(\mathbb{W}_h(X,\varsigma)',\mathbb{W}_h(Y,\tau)'\bigr)$. The action of
$^{\mathrm{t}}\mathfrak{F}_h$ on a regular state $\omega_h\in \mathbb{W}_h(X,\varsigma)'_+$ is defined by
the generating map
\begin{equation*}
  \mathcal{G}_{^{\mathrm{t}}\mathfrak{F}_h (\omega_h)}(\cdot )=e^{i F_h(\cdot)}\mathcal{G}_{\omega_h}(f(\cdot ))\; .
\end{equation*}
In general $\mathcal{G}_{^{\mathrm{t}}\mathfrak{F}_h (\omega_h)}$ fails to be of almost positive type, and
thus $^{\mathrm{t}}\mathfrak{F}_h\bigl(\mathbb{W}_h(X,\varsigma)'_+\bigr)\nsubseteq
\mathbb{W}_h(Y,\tau)'_+$.

Nonetheless, consider a semiclassical quantum state $(\omega_{h_{\beta}})_{\beta\in B}$ such that
\begin{equation*}
  \omega_{h_{\beta}}\underset{h_{\beta}\to 0}{\overset{\mathfrak{P}\vee\mathfrak{T}}{\longrightarrow}} M\in \mathcal{M}_{\mathrm{cyl}}(V)\; .
\end{equation*}
If $\lim_{\beta\in B} F_{h_{\beta}}=F_0$ (pointwise), then
\begin{align}
  \label{eq:1017}
  \mspace{165mu}&g(y):= \lim_{\beta\in B} \mathcal{G}_{^{\mathrm{t}}\mathfrak{F}_{h_{\beta}} (\omega_{h_{\beta}})}(y)=e^{i F_0(y)}\hat{M}\bigl(f(y)\bigr)&(\forall y\in Y)\; .
\end{align}
This can also be extended to other maps that are not homomorphisms. Let $\xi\in X$ be fixed; then the map
\begin{equation*}
  \mathfrak{F}_h^{(\xi)}\bigl(W_h(x)\bigr)=W_h(\xi)W_h(x)=e^{-ih \varsigma(\xi,x)}W_h(x+\xi)
\end{equation*}
extends by linearity to a map $\mathfrak{F}_h^{(\xi)}:\mathbb{W}_h(X,\varsigma)\to
\mathbb{W}_h(X,\varsigma)$. For $\mathfrak{F}_h^{(\xi)}$, \cref{eq:1017} becomes
\begin{equation}
  \label{eq:1018}
  g^{(\xi)}(x)=\hat{M}(x+\xi)\;,\quad \forall x\in X\; .
\end{equation}
\cref{eq:1018} defines a (complex) signed cylindrical measure on $V$, that we may denote (with a slight
abuse of notation) by $e^{i\xi(z)}\mathrm{d}M(z)$ or $e^{i\xi(\cdot )}M$. In other words, we have the
following characterization of the Wigner measures associated to
$\bigl(\,\!^{\mathrm{t}}\mathfrak{F}_{h_{\beta}}^{(\xi)}\omega_{h_{\beta}}\bigr)_{\beta\in B}$:
\begin{equation}
  \label{eq:1019}
  \mathscr{W}_{\mathfrak{P}\vee\mathfrak{T}}\Bigl(\omega_{h_{\beta}}(W_{h_{\beta}}(\xi)\,\cdot \,),\beta\in B\Bigr)=\Bigl\{e^{i\xi(\cdot )}M\; ,\; M\in\mathscr{W}_{\mathfrak{P}\vee\mathfrak{T}}(\omega_{h_{\beta}},\beta\in B)\Bigr\}\; . 
\end{equation}
We are interested in giving a similar characterization for other mappings. The part $f^{(\xi)}$ of
$\mathfrak{F}_h^{(\xi)}$ is non-linear, and for more general non-linear transformations it may be
difficult to obtain an explicit formula of type \eqref{eq:1019}. We restrict to linear transformations
$f$. Let us remark that if $F_h=0$, then $f$ should be a symplectic transformation for $\mathfrak{F}_h$ to
be a *-homomorphism; in that case \cref{thm:1013} is equivalent to \cref{eq:8}.

\begin{proposition}
  \label{thm:1013}
  Let $(X,\varsigma),(Y,\tau)\in \mathbf{Symp}_{\mathds{R}}$, $V,W\in \mathbf{TVS}_{\mathds{R}}$ in
  compatible duality separating points with $X$ and $Y$ respectively. For any weakly
  continuous\footnote{Weak continuity here means continuity with respect to the $\sigma(Y,W)$ and
    $\sigma(X,V)$ topologies.} linear map $\mathrm{u}:Y\to X$, and for any $F_h\in \mathds{R}^{Y}$ such
  that for any $y\in Y$, $\lim_{h\to 0}F_h(y)=0$, define
  $\mathfrak{F}_h^{(\mathrm{u})}=(\mathrm{u},F_h)$. Then for any semiclassical quantum state
  $(\omega_{h_{\beta}})_{\beta\in B}$ on $\bigl(\mathbb{W}_{h_{\beta}}(X,\varsigma)\bigr)_{\beta\in B}$
  with no loss of mass:
  \begin{equation}
    \label{eq:1020}
    \mathscr{W}_{\mathfrak{P}\vee\mathfrak{T}}\Bigl(\, ^{\mathrm{t}}\mathfrak{F}_{h_{\beta}}^{(\mathrm{u})}(\omega_{h_{\beta}})\,,\,\beta\in B\,\Bigr)=\Bigl\{\,\!^{\mathrm{t}}\mathrm{u}\, _{*}\, M\, ,\, M\in \mathscr{W}_{\mathfrak{P}\vee \mathfrak{T}}(\omega_{h_{\beta}},\beta\in B)\Bigr\}\; ,
  \end{equation}
  where $^{\mathrm{t}}\mathrm{u}:V\to W$ is the transposed map of $\mathrm{u}$.
\end{proposition}
\begin{proof}
  The definition of push-forward (image) of a cylindrical measure is briefly recalled in
  \cref{sec:maps-conv-prod,sec:locally-conv-spac} \citep[see][for additional
  details]{bourbaki1969int9
    ,schwartz1973tirsm
    ,vakhania1987ma
  }; let us recall that with the assumptions above the transposed map $^{\mathrm{t}}\mathrm{u}$ is
  continuous with respect to the $\sigma(V,X)$ and $\sigma(W,Y)$ topologies \citep[see, \emph{e.g.},][EVT
  II.50]{bourbaki1981evt1-5
  }. \cref{eq:1017} for $\mathfrak{F}_h^{(\mathrm{u})}$ becomes
  \begin{align*}
    \mspace{170mu}&g(y)=\hat{M}\bigl(\, u(y)\bigr)=\widehat{\!\,^{\mathrm{t}}\mathrm{u}\,_{*}\, M}(y)&(\forall y\in Y)\; .
  \end{align*}
\end{proof}
\begin{remark}
  \label{rem:7}
  Choosing $V=X^{*}_X$ and $W=Y^{*}_Y$, \cref{thm:1013} holds \emph{for any linear map} $\mathrm{u}:Y\to
  X$, since all linear maps are weakly continuous with respect to $\sigma(Y,Y^{*})$ and $\sigma(X,X^{*})$.
\end{remark}

\subsection{Quantum convolution and convolution of measures}
\label{sec:free-conv-regul}

From Theorem~\ref{thm:107} it also follows that the Wigner measures associated to the quantum convolution
of regular quantum states are convoluted measures. The convolution $\circledast$ of two cylindrical
measures $M,N\in \mathcal{M}_{\mathrm{cyl}}(V)$ is defined as the pushforward of the product cylindrical
measure $M\otimes N\in \mathcal{M}_{\mathrm{cyl}}(V\times V)$ by means of the addition map $\alpha:V\times
V\to V$ defined by $(v,w)\mapsto v+w$. In other words,
\begin{equation*}
  M\circledast N=\alpha\,_{*}\, (M\otimes N)\; .
\end{equation*}
It is a well-known result that
\begin{equation*}
  \widehat{M\circledast N\,}=\hat{M}\cdot \hat{N}\; ,
\end{equation*}
where the dot stands for the multiplication of complex-valued functions.

On the other hand, we define the \emph{quantum convolution} as the quantum counterpart of the convolution
of cylindrical measures. Let $\omega_h,\varpi_h\in \mathbb{W}_h(X,\varsigma)'_+$ be regular quantum
states, their quantum convolution $\omega_h\star \varpi_h\in \mathbb{W}_h(X,\varsigma)'$ is the state
defined as follows. On generators, define the action
\begin{equation*}
  (\omega_h\star \varpi_h)\bigl(W_h(\xi)\bigr)=\omega_h\bigl(W_h(\xi)\bigr)\varpi_h\bigl(W_h(\xi)\bigr)\; .
\end{equation*}
On finite linear combinations one has
\begin{equation*}
  (\omega_h\star \varpi_h)\Bigl(\sum_{j=1}^n \lambda_j W_h(\xi_j) \Bigr)=\sum_{j=1}^n\lambda_j \, (\omega_h\star \varrho_h)\bigl(W_h(\xi_j)\bigr)\; .
\end{equation*}
Then $\omega_h\star \varpi_h$ extends to a state on $\mathbb{W}_h(X,\varsigma)$.

\begin{proposition}
  \label{prop:108}
  Let $(X,\varsigma)\in \mathbf{Symp}_{\mathds{R}}$, $V\in \mathbf{TVS}_{\mathds{R}}$ in compatible
  duality with $X$. Consider two semiclassical quantum states with no loss of mass
  $(\omega_{h_{\beta}})_{\beta\in B}$ and $(\varpi_{h_{\beta}})_{\beta\in B}$ on
  $\bigl(\mathbb{W}_{h_{\beta}}(X,\varsigma)\bigr)_{\beta\in B}$. If
  \begin{equation*}
    \omega_{h_{\beta}}\underset{h_{\beta}\to 0}{\overset{\mathfrak{P}\vee \mathfrak{T}}{\longrightarrow}}M \; ,\; \varpi_{h_{\beta}}\underset{h_{\beta}\to 0}{\overset{\mathfrak{P}\vee \mathfrak{T}}{\longrightarrow}}N\; ; 
  \end{equation*}
  then
  \begin{equation*}
    \omega_{h_{\beta}}\star \varpi_{h_{\beta}}\underset{h_{\beta}\to 0}{\overset{\mathfrak{P}\vee \mathfrak{T}}{\longrightarrow}} M \circledast N\; .
  \end{equation*}
\end{proposition}
\begin{proof}
  By definition of quantum convolution, the generating functional of $\omega_h\star \varpi_h$ satisfies:
  \begin{equation*}
    \mathcal{G}_{\omega_h\star \varpi_h}(\cdot )=\mathcal{G}_{\omega_h}(\cdot )\mathcal{G}_{\varpi_h}(\cdot )\; .
  \end{equation*}
  Therefore pointwise
  \begin{equation*}
    \lim_{\beta\in B}\mathcal{G}_{\omega_{h_{\beta}}\star \varpi_{h_{\beta}}}(\cdot )=\hat{M}(\cdot )\hat{N}(\cdot )=\widehat{M \circledast N\,}(\cdot )\; .
  \end{equation*}
\end{proof}
\section{States on the C*-algebra of almost periodic functions}
\label{sec:states-c-algebra}

In \cref{sec:cylindr-wign-meas} we discussed how it is possible to identify the limit points of
semiclassical quantum states, as $h\to 0$, with cylindrical measures on some topological vector space. In
this section we provide an alternative identification, that is perhaps not so interesting for
semiclassical analysis, but fits with the ideas of deformation theory.

In the definition of Weyl deformation $(\mathfrak{W}_h)_{h\geq 0}$, the C*-algebra at $h=0$ is the tensor
product of $\mathfrak{A}$ with an abelian C*-algebra of almost periodic functions. Let $G$ be a
topological abelian group, let us denote by $\mathbb{A}\!\mathbb{P}(G)$ the algebra of (continuous) almost
periodic functions. It has the following characterization:
\begin{equation*}
  \mathbb{A}\!\mathbb{P}(G)=\textrm{C*}\bigl\{\hat{G} \bigr\}\; ,
\end{equation*}
where the completion is intended with respect to the supremum norm. $\hat{G}$ is here the character group
of $G$, \emph{i.e.}\ it is the set of continuous group homomorphisms from $G$ to the multiplicative group
$\mathds{C}_1^{\times }=\bigl(\bigl\{z\in \mathds{C}, \lvert z \rvert_{}^{}=1\bigr\},\times
\bigr)$. Consider now a topological vector space $V\in \mathbf{TVS}_{\mathds{R}}$ as an abelian
topological group with respect to the addition operation (more precisely, apply the suitable forgetful
functor $\mathbb{F}\!\mathbb{F}$). There is a natural subalgebra
$\mathbb{L}\mspace{-1mu}\mathbb{A}\!\mathbb{P}(\mathbb{F}\!\mathbb{F}V)\subset
\mathbb{A}\!\mathbb{P}(\mathbb{F}\!\mathbb{F}V)$, defined as
\begin{equation*}
  \mathbb{L}\mspace{-1mu}\mathbb{A}\!\mathbb{P}(V) :=\mathbb{L}\mspace{-1mu}\mathbb{A}\!\mathbb{P}(\mathbb{F}\!\mathbb{F}V)=\textrm{C*}\bigl\{e^{2i \xi(\cdot )}\; ,\; \xi\in V' \bigr\}\; .
\end{equation*}

Therefore the $h=0$ algebra of the Weyl deformation defined in~\eqref{eq:112} can be identified as
follows: $\mathfrak{W}_0\cong \mathbb{L}\mspace{-1mu}\mathbb{A}\!\mathbb{P}(V)\otimes_{\gamma_0} \mathfrak{A}$ for
any topological vector space $V$ in compatible duality with $X$. The idea is to show that the cylindrical
measures $\mathcal{M}_{\mathrm{cyl}}(V;\mathfrak{A}'_+)$ can be identified with algebraic states of
$\mathfrak{W}_0$ whenever $\gamma_0=\gamma_{\mathrm{max}}$, \emph{i.e.}\ they can be identified as
(positive) elements of $( \mathbb{L}\mspace{-1mu}\mathbb{A}\!\mathbb{P}(V)\otimes_{\gamma_{\mathrm{max}}}
\mathfrak{A})'$.

\begin{lemma}\label{lemma:105}
  Let $V$ be a locally convex space in compatible duality with $X$. In addition, let $\mathfrak{A}\in
  \mathbf{C^{*}alg}$; and $M\in \mathcal{M}_{\mathrm{cyl}}(V;\mathfrak{A}'_+)$. Then there exists a
  completely positive map $\mathscr{F}_M\in
  \mathcal{B}\bigl(\mathbb{L}\mspace{-1mu}\mathbb{A}\!\mathbb{P}(V),\mathfrak{A}'\bigr)$.
\end{lemma}
\begin{proof}
  Let us define the map $\mathscr{F}_M$ by its action on the generators: let $x\in \Phi^{\circ} \subset
  X$,
  \begin{equation*}
    \mathscr{F}_M(e^{2ix(\cdot )})=\int_{V/\Phi}^{}e^{2ix(v)}  \mathrm{d}\mu_{\Phi}(v) =\hat{M}(x) \; .
  \end{equation*}
  It then follows that $\mathscr{F}_M$ is linear. Let
  \begin{equation*}
    f_n(\cdot )=\sum_{j=1}^nz_j e^{2i x_j(\cdot )}
  \end{equation*}
  be a complex linear combination of generators, then
  \begin{equation*}
    \mathscr{F}_M(f_n)=\sum_{j=1}^nz_j\hat{M}(x_j)=\int_{V/\Phi_n}^{}\sum_{j=1}^n z_je^{2i x_j(v)}  \mathrm{d}\mu_{\Phi_n}(v) \; ,
  \end{equation*}
  where $\Phi_n\in F(V)$ is such that $\bigl\{x_1,\dotsc,x_n\bigr\}\subset \Phi_n^{\circ}$. Using the
  corresponding result for standard measures $\mu_{\Phi_n,\kappa}$ (see \cref{sec:measures-with-values}
  for the notation), it is not difficult to prove that
  \begin{equation*}
    \Bigl\lVert \mathscr{F}_M(f_n) \Bigr\rVert_{\mathfrak{A}'}^{}\leq \lVert\, f_n  \rVert_{\infty}^{}\lVert \hat{M}(0)  \rVert_{\mathfrak{A}'}^{}\; ,
  \end{equation*}
  where $\lVert \,\cdot\, \rVert_{\infty}^{}$ denotes the supremum norm. Now let $(f_n)_{n\in \mathds{N}}$
  be a Cauchy sequence with respect to the supremum norm, that converges to $f\in
  \mathbb{L}\mspace{-1mu}\mathbb{A}\!\mathbb{P}(V)$. Then it is possible to define $\mathscr{F}_M(f)$ as the strong
  limit of $\mathscr{F}_M(f_n)$ in $\mathfrak{A}'$, and
  \begin{equation*}
    \Bigl\lVert\mathscr{F}_M(f)  \Bigr\rVert_{\mathfrak{A}'}^{}\leq \lVert\, f  \rVert_{\infty}^{}\lVert \hat{M}(0)  \rVert_{\mathfrak{A}'}^{}\; .
  \end{equation*}
  Finally, $\mathscr{F}_M$ is completely positive since it is completely positive on linear combinations
  of generators by Bochner's theorem,\cref{thm:103}.
\end{proof}
\begin{corollary}\label{cor:102}
  If $\gamma_0=\gamma_{\mathrm{max}}$, then by means of the map $\mathscr{F}_M$ it is possible to
  associate a state $\Omega_M\in \bigl(\mathbb{L}\mspace{-1mu}\mathbb{A}\!\mathbb{P}(V)
  \otimes_{\gamma_{\mathrm{max}}}\mathfrak{A} \bigr)'_+$ to any cylindrical measure $M\in
  \mathcal{M}_{\mathrm{cyl}}(V;\mathfrak{A}'_+)$:
  \begin{equation*}
    \Omega_M=(\mathbb{E}_{0,1})^{-1}\mathscr{F}_M \; .
  \end{equation*}
\end{corollary}
\begin{proof}
  By Lemma~\ref{lemma:105}, we can associate to $M$ the map $\mathscr{F}_M$, and the latter is a
  completely positive element of
  $\mathcal{B}\bigl(\mathbb{L}\mspace{-1mu}\mathbb{A}\!\mathbb{P}(V),\mathfrak{A}'\bigr)$. Therefore by
  \cref{prop:26}, $\Omega_M=(\mathbb{E}_{0,1})^{-1}\mathscr{F}_M$ is a state.
\end{proof}
Finally, the application of Tietze's extension theorem yields an extension $\tilde{\Omega}_M\in
\bigl(\mathbb{A}\!\mathbb{P}(\mathbb{F}\!\mathbb{F}V) \otimes_{\gamma_{\mathrm{max}}}\mathfrak{A}
\bigr)'_+$ of $\Omega_M$ to the algebra of almost periodic functions.
\begin{proposition}\label{thm:1015}
  Let $(X,\varsigma)\in \mathbf{Symp}_{\mathds{R}}$, and $V$ a locally convex space in compatible duality
  with $X$. In addition, let $\mathfrak{A}\in \mathbf{C^{*}alg}$, and $(\mathfrak{W}_h)_{h\geq 0}$ the
  associated Weyl deformation \eqref{eq:112}. Consider a semiclassical quantum state
  $(\omega_{h_{\beta}})_{\beta\in B}$, then for any $M\in \mathscr{W}(\omega_{h_{\beta}},\beta\in
  B)\subset \mathcal{M}_{\mathrm{cyl}}(V;\mathfrak{A}'_+) $,
  $\Omega_{M}=(\mathbb{E}_{0,1})^{-1}\mathscr{F}_M$ is a state on $\mathbb{L}\mspace{-1mu}\mathbb{A}\!\mathbb{P}(V)
  \otimes_{\gamma_{\mathrm{max}}}\mathfrak{A}$. In addition, the state $\Omega_M$ can be extended
  continuously to a state $\tilde{\Omega}_M \in \bigl(\mathbb{A}\!\mathbb{P}(\mathbb{F}\!\mathbb{F}V)
  \otimes_{\gamma_{\mathrm{max}}}\mathfrak{A} \bigr)'_+$.
\end{proposition}
\begin{remark}
  \label{rem:101}
  In general, the state $\Omega_M$ may have a mass defect, \emph{i.e.}\ $\Omega_M(1)\neq \lim_{h_b}
  \omega_{h_b}(1)$. However, if the converging sequence $(\omega_{h_b})_{b\in \underline{B}}$ has no loss
  of mass, then $\Omega_M(1)= \lim_{h_b} \omega_{h_b}(1)$.
\end{remark}
\begin{remark}
  \label{rem:6}
  One could define the ``no-quantization'' functor $\mathbb{W}_0: \mathbf{Symp}_{\mathds{R}}\to
  \mathbf{C^{*}alg}$ by $\mathbb{W}_0(X,\varsigma)=\mathbb{A}\!\mathbb{P}(\mathbb{F}\!\mathbb{F}X^{*}_X)$.
\end{remark}

\appendix

\section{Elements of cone-valued measure theory}
\label{sec:measures-with-values}

In this appendix, we outline some results of vector integration as a reference. For our purpose, vector
measures with values in cones behave essentially as standard measures, and in particular Bochner's theorem
is valid for cylindrical cone-valued measures.

\subsection{Definition of cone-valued measures}
\label{sec:definitions}

\begin{itemize}
\item Given a measurable space $E$, we denote by $\Sigma$ its $\sigma$-algebra.
\item We will always denote by $X$ a real vector space, and by $C$ a \emph{pointed and generating} convex
  cone in $X$ containing $0$. This means that
  \begin{equation*}
    C \cap -C=\{0\}\; ; \; C- C=X\; .
  \end{equation*}
\item As before, we denote by $X^{*}$ the algebraic dual of $X$, and for any $X'\subset X^{*}$ we denote
  by $C'$ the dual cone of $C$ defined by $C'=\bigl\{\kappa \in X', \kappa(C)\subseteq
  \mathds{R}^+\bigr\}$. If $X$ is locally convex, $X'$ its continuous dual and $C$ is closed, then the
  Hahn-Banach separation theorem yields
  \begin{equation}\label{eq:101}\tag{c1}
    C=C''=\bigl\{x\in X, x(C')\subseteq \mathds{R}^+\bigr\}\; .
  \end{equation}
  We will consider only triples $(X,X',C)$ satisfying~\eqref{eq:101}.
\item We denote by $\mathds{R}^+_{\infty}=[0,\infty]$ the extended real semi-line considered as an
  additive semigroup with the additional rule
  \begin{equation*}
    (\forall x\in \mathds{R}^+_{\infty})\;\infty+x=x+\infty=\infty\; .
  \end{equation*}
  We also denote by $\mathds{R}\cup \{-\infty,+\infty\}=[-\infty,+\infty]$ the (compact) complete lattice of extended reals, and by
  $\mathds{C}\cup \{\infty\}$ the extended complex numbers (one-point compactification of $\mathds{C}$).
\item We denote by $C_{\infty}=\mathrm{Hom}_{\mathrm{mon}}(C',\mathds{R}^+_{\infty})$ the subset of
  $(\mathds{R}^+_{\infty})^{C'}$ consisting of monoid homomorphisms. $C_{\infty}$ is a monoid with respect
  to pointwise addition.
\item We denote by $\mathfrak{i}_C$ the natural monoid morphism
  \begin{equation*}
    \mathfrak{i}_C:C\to C_{\infty}\; ,\; (\forall c\in C)(\forall \kappa\in C')\mathfrak{i}_C(c)(\kappa)=\kappa(c)\; .
  \end{equation*}
  $\mathfrak{i}_C(c_1)=\mathfrak{i}_C(c_2)$ yields $(\forall \kappa\in
  C')\kappa(c_1-c_2)=\kappa(c_2-c_1)=0$.  Therefore \eqref{eq:101} implies $c_1-c_2\in C \cap - C$ and by
  the pointedness of $C$ we have $c_1-c_2=0$. Thus $\mathfrak{i}_C$ is injective and $C\cong
  \mathfrak{i}_C(C)$ is a submonoid of $C_\infty$.
\item The next condition is important to define cone-valued measures:
  \begin{equation}
    \label{eq:102}\tag{c2}
    \mathfrak{i}_C(C)=\mathrm{Hom}_{\mathrm{mon}}(C',\mathds{R}^+)\; .
  \end{equation}
  We discuss later some explicit example of triples that satisfy \eqref{eq:101} and \eqref{eq:102}.
\end{itemize}

Cone valued measures are vector measures generalizing the concept of positive measures. They can be seen
as suitable collections of the latter, and therefore they share many interesting properties with ``usual''
positive measures.
\begin{definition}[$C$-valued measures]\label{def:101}
  Let $(X,X',C)$ be a triple that satisfies \eqref{eq:101}-\eqref{eq:102}, and $E$ a measurable space.
  Then $\mu\in (C_{\infty})^{\Sigma}$ is a $C$-valued measure on $E$ iff it is countably additive and
  $\mu(\emptyset)=0$.
\end{definition}
\begin{remark}
  In the definition above, $0$ is the trivial monoid morphism that maps every $\kappa\in C'$ to $0\in
  \mathds{R}^+_{\infty}$. In addition, countable additivity is intended in the following sense. Let
  $\bigl\{K_j\bigr\}_{j\in \mathds{N}}\subset C_{\infty}$ be a subset of $C_{\infty}$; then the countable
  combination $\sum_{j\in \mathds{N}}^{}K_j\in C_{\infty}$ is defined by pointwise convergence -- in the
  topology of extended reals -- of partial sums, i.e. by convergence of the sequences
  \begin{equation*}
    \mathds{R}^+_{\infty}\supset (w_n^{\kappa})_{n\in\mathds{N}}=\Bigl(\sum_{j=0}^nK_j(\kappa)\Bigr)_{n\in \mathds{N}}\;,\; \kappa\in C'\; .
  \end{equation*}
  Therefore a function $\mu\in (C_{\infty})^{\Sigma}$ is countably additive iff for any collection
  $\bigl\{b_j\bigr\}_{j\in \mathds{N}}\subset \Sigma$ of mutually disjoint measurable sets,
  \begin{equation*}
    \mu\Bigl(\bigcup_{j\in \mathds{N}}b_j\Bigr)=\sum_{j\in \mathds{N}}^{}\mu(b_j)\; .
  \end{equation*}
\end{remark}
If $C=\mathds{R}^+$, \cref{def:101} corresponds to that of positive Borel measures. As it was stated
before, a key feature of $C$-valued measures is that they are in fact families of positive measures,
indexed by the dual cone $C'$. The precise statement is the following, whose proof follows almost directly
from \cref{def:101} above.

\begin{thm}[\citet{neeb1998mm
}]\label{thm:101}
There is a bijection between $C$-valued measures $\mu$ on $E$ and families of positive measures
$(\mu_{\kappa})_{\kappa\in C'}$ on $E$ such that for any $b\in \Sigma$, $\bigl(\kappa\mapsto
\mu_{\kappa}(b) \bigr)\in \mathrm{Hom}_{\mathrm{mon}}(C',\mathds{R}^+_{\infty})$, \emph{i.e.}\ the map
$\kappa\mapsto \mu_{\kappa}(b)$ belongs to $C_{\infty}$.
\end{thm}

\noindent In light of \cref{thm:101}, we define a $C$-valued measure $\mu$ \emph{finite} if $\mu_{\kappa}$
is a finite positive measure for any $\kappa\in C'$.

We turn now to integration of (scalar) functions with respect to cone-valued measures. As usual, it is
convenient to start with the integration of non-negative functions. \cref{thm:101} is very convenient in
this context, since we can simply define cone-valued integration by means of usual integration. Let
$f:E\to \mathds{R}_{\infty}^+$ be a non-negative measurable function with values on the extended
reals. Let $b\in\Sigma$; then we define for any $\kappa\in C'$,
\begin{equation*}
  \mathds{R}_{\infty}^+\ni I_{\kappa}=\int_b^{}f(x) \mathrm{d}\mu_{\kappa}(x)\; .  
\end{equation*}
The map $\kappa\mapsto I_{\kappa}$ is a monoid morphism, and therefore an element of $C_{\infty}$, that we
denote by $\mu_b(f)$. This leads to the following natural definition.
\begin{samepage}
  \begin{definition}[$\mu$-integrable functions]\label{def:104}
    Let $(X,X',C)$ be a triple that satisfies \eqref{eq:101}-\eqref{eq:102}; and $\mu$ a $C$-valued
    measure on a measurable space $E$. The \emph{measure} of a \emph{non-negative} measurable function
    $f\in (\mathds{R}_{\infty}^+)^E$ is defined by
    \begin{equation*}
      C_{\infty}\ni \mu_b(f) = \biggl(\kappa\mapsto \int_b^{}f(x) \mathrm{d}\mu_{\kappa}(x)\biggr)\; .
    \end{equation*}
    A non-negative measurable function $f\in (\mathds{R}_{\infty}^+)^E$ is \emph{$\mu$-integrable} on the
    measurable set $b\in \Sigma$ iff $\mu_b(f)\in
    \mathfrak{i}_C(C)=\mathrm{Hom}_{\mathrm{mon}}(C',\mathds{R}^+)$. In this case, we denote the
    \emph{integral} by
    \begin{equation*}
      C\ni \int_b^{}f(x)  \mathrm{d}\mu(x)=\mathfrak{i}_C^{-1}\bigl(\mu_b(f)\bigr)\; .
    \end{equation*}

    A \emph{complex function} $f\in \mathds{C}^E$ is $\mu$-integrable on the measurable set $b\in \Sigma$
    iff $\lvert f \rvert_{}^{}$ is $\mu$-integrable, and
    \begin{equation*}
      \begin{split}
        X_{\mathds{C}}\ni \int_b^{}f(x)  \mathrm{d}\mu(x)= \int_b^{}(\Re f)_+(x)  \mathrm{d}\mu(x)-\int_b^{}(\Re f)_-(x)  \mathrm{d}\mu(x)\\+i\biggl(\int_b^{}(\Im f)_+(x)  \mathrm{d}\mu(x)-\int_b^{}(\Im f)_-(x)  \mathrm{d}\mu(x)\biggr) ;
      \end{split}
    \end{equation*}
    where $X_{\mathds{C}}$ is the complexification of $X$, and $f=(\Re f)_+- (\Re f)_-+i\Bigl((\Im f)_+-
    (\Im f)_-\Bigr)$ with
    \begin{equation*}
      \bigl\{(\Re f)_+,(\Re f)_-,(\Im f)_+,(\Im f)_-\bigr\}\subset (\mathds{R}^+)^E\; .
    \end{equation*}
  \end{definition}
\end{samepage}
\begin{remark}
  If $\mu$ is finite, then any $f\in \bigcap_{\kappa\in C'} L^{\infty} (E,\mathrm{d}\mu_{\kappa})$ is
  integrable. In particular, any continuous and bounded function is Borel integrable on a topological
  space $E$.
\end{remark}

In the next proposition we state the important linearity property of the integral $\int_b^{}f(x)
\mathrm{d}\mu(x)$. The proof is trivial.
\begin{proposition}\label{prop:101}
  The mapping $f\mapsto \int_b^{}f(x) \mathrm{d}\mu(x)$ is a linear, cone-homomorphism. In other words,
  for any complex-valued $\mu$-integrable functions $f_1,f_2$ and $z\in \mathds{C}$:
  \begin{equation*}
    \int_b^{}\bigl(f_1(x)+zf_2(x)\bigr)  \mathrm{d}\mu(x)=\int_b^{}f_1(x)  \mathrm{d}\mu(x)+z\int_b^{}f_2(x)  \mathrm{d}\mu(x)\; .
  \end{equation*}
  In addition, the cone of non-negative $\mu$-integrable functions is mapped into the cone
  $(C+i\{0\})\subset X_{\mathds{C}}$.
\end{proposition}

\subsection{Bochner's theorem}
\label{sec:bochners-theorem}

We are now ready to prove a result that is crucial in our framework: Bochner's theorem for finite
$C$-valued measures. To prove the theorem we follow closely
\citep{neeb1998mm
  ,glockner2003mams
}.

\subsubsection{Locally compact abelian groups}
\label{sec:finite-dimens-vect}

In this subsection -- if not specified otherwise -- we take as measure space $G$ a locally compact abelian
group with character group $\hat{G}$; and $(X,X',C)$ a triple satisfying~\eqref{eq:101}-\eqref{eq:102} and
some or all of the following additional conditions. Let $K$ be a pointed and generating cone in a real
vector space $A$. Then an involution $^{*}$ on $A_{\mathds{C}}$ agrees with $K$ if: $a^{*}a\in K+i\{0\}$
for any $a\in A_{\mathds{C}}$, and for any $k\in K$ there exists an $a_k\in A_{\mathds{C}}$ such that
$k=a_k^{*}a_k$. Then we define the following conditions:
\begin{gather}
  \label{eq:103}\tag{c3}
  X\textrm{ and }X'\textrm{ locally convex}\; ;\\
  \label{eq:105}\tag{c4}
  C'\textrm{ pointed and generating in }X'\; ;\\
  \label{eq:109}\tag{c5}
  (X')_{\mathds{C}}\textrm{ is an involutive algebra with involution agreeing with }C'\; ;\\
  \label{eq:104}\tag{c6}
  X_{\mathds{C}}\textrm{ is the continuous dual of } (X')_{\mathds{C}}\; .
\end{gather}
Given a locally convex real vector space $T$, there are (infinitely) many ways to endow $T_{\mathds{C}}$
with a topology in a ``natural'' way (\emph{i.e.}\ satisfying some suitable properties). Therefore one may
ask if it is always possible to endow $X_{\mathds{C}}$ and $(X')_{\mathds{C}}$ with suitable topologies
such that \eqref{eq:104} is satisfied. If $X'$ is a Banach space and $X$ its continuous dual, the answer
is that for any so-called natural complexification of $X'$ there is a so-called reasonable
complexification of $X$ such that \eqref{eq:104} is satisfied.

\begin{definition}[Completely positive functions]\label{def:103}
  Let $G$ be an abelian group, $(X,X',C)$ a triple satisfying~\eqref{eq:101}
  and~\eqref{eq:103}-\eqref{eq:104}. A function $f\in (X_{\mathds{C}})^G$ is \emph{completely positive}
  iff for any $n\in \mathds{N}_{*}$, for any $\bigl\{g_i\bigr\}_{i=1}^n\subset G$ and
  $\bigl\{\tilde{\kappa}_i\bigr\}_{i=1}^n\in (X')_{\mathds{C}}$:
  \begin{equation*}
    \sum_{i,j=1}^n\tilde{\kappa}_{j}^{*}\tilde{\kappa}_i\Bigl(f\bigl(g_ig_j^{-1}\bigr)\Bigr)\geq 0\; .
  \end{equation*}
\end{definition}
\noindent The definition above is the analogous of positive-definiteness for the cone $C$. In order to
study completely positive functions, it is convenient to introduce a slight generalization. Let $f\in
(X_{\mathds{C}})^G$, where $G$ is an abelian group. Then there exist an associated kernel $F_f(\cdot
,\cdot):G\times G\to X_{\mathds{C}}$ defined by $F_f(g_1,g_2)=f(g_1g_2^{-1})$. Hence it is natural to have
the following definition.
\begin{definition}[Completely positive kernels]\label{def:108}
  Let $A$ be a set, $(X,X',C)$ a triple satisfying~\eqref{eq:101} and~\eqref{eq:103}-\eqref{eq:104}. A
  kernel $F:A\times A\to X_{\mathds{C}}$ is completely positive iff for any $n\in \mathds{N}_{*}$, for any
  $\bigl\{a_i\bigr\}_{i=1}^n\subset A$ and $\bigl\{\tilde{\kappa}_i\bigr\}_{i=1}^n\in (X')_{\mathds{C}}$:
  \begin{equation*}
    \sum_{i,j=1}^n\tilde{\kappa}_{j}^{*}\tilde{\kappa}_i\bigl(F(a_i,a_j)\bigr)\geq 0\; .
  \end{equation*}
\end{definition}
\noindent The equivalence of the two definitions for groups is given by the following trivial result.
\begin{lemma}\label{lemma:102}
  Let $G$ be an abelian group, $(X,X',C)$ a triple satisfying~\eqref{eq:101}
  and~\eqref{eq:103}-\eqref{eq:104}. A function $f\in (X_{\mathds{C}})^G$ is completely positive iff the
  associated kernel $F_f\in (X_{\mathds{C}})^{G\times G}$ is completely positive.
\end{lemma}
\noindent In order to prove Bochner's theorem, we prove a couple of preliminary results related to
complete positivity.
\begin{lemma}\label{lemma:104}
  Let $(X,X',C)$ be a triple satisfying~\eqref{eq:101}-\eqref{eq:104}, and $\mu$ a finite $C$-valued Borel
  measure on a topological space $E$. If we denote by $L^2(E,\mathrm{d} \mu)\subset \mathds{C}^E$ the
  space of $\mu$-a.e.\ square integrable functions, \emph{i.e.}\
  \begin{equation*}
    L^2(E,\mathrm{d} \mu)=\bigcap_{\kappa\in C'}L^2(E, \mathrm{d} \mu_k)\; ;
  \end{equation*}
  then the integral map $I_{\mu}:L^2(E,\mathrm{d}\mu)\times L^2(E,\mathrm{d}\mu)\to X_{\mathds{C}}$,
  defined by
  \begin{equation*}
    I_{\mu}(f,g)=\int_E^{}f(x)\bar{g}(x)  \mathrm{d}\mu(x)\; ,
  \end{equation*}
  is well-defined and a completely positive kernel.
\end{lemma}
\begin{proof}
  The fact that the kernel $I_{\mu}$ is well-defined is easy to prove using the corresponding property of
  $I_{\mu_{\kappa}}$, $\kappa\in C'$. To prove complete positivity, we proceed as follows. Let $n\in
  \mathds{N}_{*}$, $\bigl\{\,f_i\bigr\}_{i=1}^n\subset L^2(E,\mu)$ and
  $\bigl\{\tilde{\kappa}_i\bigr\}_{i=1}^n\in (X')_{\mathds{C}}$. Using the decomposition
  $X'_{\mathds{C}}=C'-C'+i(C'-C')$, we see that the map
  \begin{equation*}
    X'_{\mathds{C}}\ni \tilde{\kappa}\mapsto \mu_{\tilde{\kappa}}=\mu_{\kappa^+_{\mathrm{R}}}-\mu_{\kappa^-_{\mathrm{R}}}+i(\mu_{\kappa^+_{\mathrm{I}}}-\mu_{\kappa^-_{\mathrm{I}}})
  \end{equation*}
  defines a linear morphism from $X'_{\mathds{C}}$ to the standard signed measures. Now let
  $\mu_{\tilde{\kappa}}$ be a signed measure, $f$ an everywhere $\mu_{\tilde{\kappa}}$-integrable
  function. Then there is a signed measure $\mu_{\tilde{\kappa}(f)}$ defined by
  $d\mu_{\tilde{\kappa}(f)}(x)=f(x)d\mu_{\tilde{\kappa}}(x)$. If we define in addition
  \begin{equation*}
    \begin{split}
      \mu_{\tilde{\kappa}(f)^{*}}=\mu_{\tilde{\kappa}^{*}(\bar{f})}\; ,\; \mu_{\tilde{\kappa}_1(f_1)+\tilde{\kappa}_2(f_2)}=\mu_{\tilde{\kappa}_1(f_1)}+\mu_{\tilde{\kappa}_2(f_2)}\; ,\; 
      \mu_{\tilde{\kappa}_1(f_1)\tilde{\kappa}_2(f_2)}=\mu_{\tilde{\kappa}_1\tilde{\kappa}_2(f_1f_2)}\; ;
    \end{split}
  \end{equation*}
  then it is easy to see, using property~\eqref{eq:109}, that for any Borel set $b\in \mathscr{B}(E)$,
  $\tilde{\kappa}\in X'_{\mathds{C}}$ and everywhere $\mu_{\tilde{\kappa}}$-integrable $f$:
  \begin{equation*}
    \int_b^{}  \mathrm{d}\mu_{\tilde{\kappa}(f)^{*}\tilde{\kappa}(f)}\geq 0\; .
  \end{equation*}
  Then
  \begin{equation*}
    \begin{split}
      \sum_{i,j=1}^n\tilde{\kappa}_j^{*}\tilde{\kappa}_i\bigl(I_{\mu}(f_i,f_j)\bigr)\mspace{-3mu}=\mspace{-7mu}\sum_{i,j=1}^n\int_E^{}f_i(x)\bar{f}_j(x)  \mathrm{d}\mu_{\tilde{\kappa}_j^{*}\kappa_i}\mspace{-2mu}=\mspace{-7mu}\sum_{i,j=1}^n\int_E^{}  \mathrm{d}\mu_{\tilde{\kappa}_j(f_j)^{*}\tilde{\kappa}_i(f_i)}
      \mspace{-2mu}\\=\mspace{-3mu}\int_E^{}  \mathrm{d}\mu_{\bigl(\sum_{i=1}^n\tilde{\kappa}_i(f_i)\bigr)^{*}\bigl(\sum_{i=1}^n\tilde{\kappa}_i(f_i)\bigr)}\mspace{-3mu}\geq 0\; .
    \end{split}
  \end{equation*}
\end{proof}
\begin{corollary}\label{cor:1}
  Let $(X,X',C)$ be a triple satisfying~\eqref{eq:101}-\eqref{eq:104}, and $\mu$ a finite $C$-valued Borel
  measure on a topological space $E$. If we denote by $L^\infty(E,\mathrm{d}\mu)\subset \mathds{C}^E$ the
  space of $\mu$-a.e.\ bounded functions, \emph{i.e.}\
  \begin{equation*}
    L^\infty(E,\mathrm{d}\mu)=\bigcap_{\kappa\in C'}L^\infty(E,\mathrm{d}\mu_k)\; ;
  \end{equation*}
  then the integral $I_{\mu}:L^{\infty} (E,\mathrm{d}\mu)\to X_{\mathds{C}}$ is a completely positive
  function -- considering $L^{\infty} (E,\mathrm{d}\mu)$ as an abelian multiplicative group.
\end{corollary}

The last ingredient needed to formulate Bochner's theorem is the Fourier transform. The Fourier transform
extends quite naturally to cone-valued measures.
\begin{definition}[Fourier transform of a $C$-valued measure]\label{def:102}
  Let $G$ be a locally compact abelian group, $(X,X',C)$ a triple
  satisfying~\eqref{eq:101}-\eqref{eq:102}, and $\mathcal{M}(\hat{G},C)$ the set of finite $C$-valued
  Borel measures on the character group $\hat{G}$. The \emph{Fourier transform} is a map
  $\hat{\phantom{a}}: \mathcal{M}(\hat{G},C)\to (X_{\mathds{C}})^{G}$, defined by
  \begin{align*}
    \mspace{160mu}&\hat{\mu}(g)=\int_{\hat{G}}^{}\gamma(g)  \mathrm{d}\mu(\gamma)&(\forall g\in G)\; .
  \end{align*}
\end{definition}
\noindent Using the definitions above, Bochner's theorem is written in a rather familiar form.
\begin{thm}[Bochner]\label{thm:102}
  Let $G$ be a locally compact abelian group, $(X,X',C)$ a triple
  satisfying~\eqref{eq:101}-\eqref{eq:104}. The Fourier transform is a bijection between finite $C$-valued
  measures on $\hat{G}$ and completely positive ultraweakly continuous functions from $G$ to
  $X_{\mathds{C}}$.
\end{thm}
\begin{proof}
  Let $\mu$ be a finite $C$-valued measure on $\hat{G}$. Finiteness of the measure implies the
  integrability of $\gamma(g)$, since $(\forall \gamma\in \hat{G})(\forall g\in G)\lvert \gamma(g)
  \rvert_{}^{}=1$. In addition, $\gamma$ is a representation of the abelian group $G$ on the functions
  $L^{\infty}(G,\mu)$. Hence it follows by \cref{cor:1} that $\hat{\mu}(\cdot )$ is completely
  positive. To prove ultraweak continuity, let $\kappa\in C'+i\{0\}$. By \cref{def:104}
  \begin{equation*}
    \kappa\bigl(\hat{\mu}(\cdot )\bigr)=\int_{\hat{G}}^{}\gamma(\cdot )  \mathrm{d}\mu_{\kappa}(\gamma)
  \end{equation*}
  is the Fourier transform of the finite measure $\mu_{\kappa}$, hence continuous. Now by \eqref{eq:105},
  $(X')_{\mathds{C}}=C'-C'+i (C'-C')$ and therefore for any $\tilde{\kappa}\in (X')_{\mathds{C}}$,
  $\tilde{\kappa}\bigl(\hat{\mu}(\cdot )\bigr)\in \mathds{C}^G$ is continuous. By~\eqref{eq:104}, this
  yields the ultraweak continuity of $\hat{\mu}(\cdot )$.

  Now let us consider a completely positive ultraweakly continuous function $f$ from $G$ to
  $X_{\mathds{C}}$. Then for any $\kappa\in C'+i\{0\}$, $\kappa\bigl(f(\cdot )\bigr)$ is a positive
  definite continuous $\mathds{C}$-valued function. Continuity trivially follows from ultraweak continuity
  (since $\kappa\in (X')_{\mathds{C}}$). To prove positive-definiteness, we exploit complete
  positivity. By \cref{def:103}, for any $n\in \mathds{N}_{*}$, $\bigl\{g_i\bigr\}_{i=1}^n\subset G$ and
  $\bigl\{\tilde{\kappa}_i\bigr\}_{i=1}^n\subset (X')_{\mathds{C}}$,
  \begin{equation*}
    \sum_{i,j=1}^n\tilde{\kappa}_{j}^{*}\tilde{\kappa}_i\Bigl(f\bigl(g_ig_j^{-1}\bigr)\Bigr)\geq 0\; .
  \end{equation*}
  Then by property~\eqref{eq:109}, there exists $\tilde{\kappa}_{\kappa}\in (X')_{\mathds{C}}$ such that
  $\kappa=\tilde{\kappa}^{*}_{\kappa}\tilde{\kappa}_{\kappa}$. So we can choose
  $\tilde{\kappa}_{i}=z_i\tilde{\kappa}_{\kappa}$ for any $i\in \bigl\{1,\dotsc,n\bigr\}$, where $z_i\in
  \mathds{C}$. Therefore by linearity we obtain
  \begin{equation*}
    \sum_{i,j=1}^n\bar{z}_jz_i\kappa\Bigl(f\bigl(g_ig_j^{-1}\bigr)\Bigr)\geq 0\; ;
  \end{equation*}
  and hence positive-definiteness of $\kappa\bigl(f(\cdot )\bigr)$.

  The classical Bochner's theorem for locally compact abelian groups \citep[see,
  \emph{e.g.},][]{loomis1953vn
  } implies the existence of a unique positive, finite measure $\mu_{\kappa}$ such that
  $\kappa\bigl(f(\cdot )\bigr)=\hat{\mu}_{\kappa}(\cdot )$. Therefore we have a unique family of positive
  and finite measures $(\mu_{\kappa})_{\kappa\in C'}$. In order for it to define a unique finite
  $C$-valued measure, it is necessary that $\kappa\mapsto \mu_{\kappa}$ is additive. Let
  $\kappa_1,\kappa_2\in C'$.  Then $\kappa_1+\kappa_2\in C'$, and there is a unique measure
  $\mu_{\kappa_1+\kappa_2}$ such that $\hat{\mu}_{\kappa_1+\kappa_2}(\cdot
  )=(\kappa_1+\kappa_2)\bigl(f(\cdot )\bigr)=\kappa_1\bigl(f(\cdot )\bigr)+\kappa_2\bigl(f(\cdot
  )\bigr)=\hat{\mu}_{\kappa_1}(\cdot )+\hat{\mu}_{\kappa_2}(\cdot )$.  However since the Fourier transform
  is a linear bijection, it follows that $\mu_{\kappa_1+\kappa_2}=\mu_{\kappa_1}+\mu_{\kappa_2}$. Hence by
  \cref{thm:101} we have defined a unique $C$-valued measure $\mu$. In addition, by \cref{def:104} for any
  $\kappa\in C'$
  \begin{equation*}
    \kappa\bigl(f(\cdot )\bigr)=\int_{\hat{G}}^{}\gamma(\cdot )  \mathrm{d}\mu_{\kappa}(\gamma)=\kappa\biggl(\int_{\hat{G}}^{}\gamma(\cdot )  \mathrm{d}\mu(\gamma)\biggr)\; .
  \end{equation*}
  Now by \eqref{eq:105}, it follows that for any $\tilde{\kappa}\in (X')_{\mathds{C}}$
  \begin{equation*}
    \tilde{\kappa}\bigl(f(\cdot )\bigr)=\tilde{\kappa}\biggl(\int_{\hat{G}}^{}\gamma(\cdot )  \mathrm{d}\mu(\gamma)\biggr)\; ,
  \end{equation*}
  and therefore by \eqref{eq:104} it follows that
  \begin{equation*}
    f(\cdot )=\int_{\hat{G}}^{}\gamma(\cdot )  \mathrm{d}\mu(\gamma)\; .
  \end{equation*}
\end{proof}

\subsubsection{Topological vector spaces}
\label{sec:locally-conv-spac}

Bochner's Theorem, \cref{thm:102} can be applied to finite dimensional real vector spaces (seen as abelian
groups under addition). In that context, the Fourier transform takes the following form. Let $V$ be a
finite dimensional vector space, $V'$ its continuous dual. Given a $C$-valued measure on $V$, then its
Fourier transform is a function from $V'$ to $X_{\mathds{C}}$ defined by
\begin{equation*}
  \hat{\mu}(\omega)=\int_V^{}e^{2i\omega(v)}  \mathrm{d}\mu(v)\; .
\end{equation*}
Using cylinders, we obtain a variant of Bochner's theorem for cylindrical measures on topological real
vector spaces with arbitrary dimension.
\begin{definition}[$C$-valued cylindrical measure]\label{def:105}
  Let $V$ be a topological real vector space, $F(V)$ the set of its $\sigma(V,V')$-closed subspaces with
  finite codimension (ordered by inclusion); and $(X,X',C)$ a triple satisfying
  \eqref{eq:101}-\eqref{eq:102}. A family of measures $M=(\mu_{\Phi})_{\Phi\in F(V)}$ is a cylindrical
  measure iff it is a projective system of $C$-valued measures on the family $Q(V)$ of finite dimensional
  quotients of $V$.

  In other words, the family $(\mu_{\Phi})_{\Phi\in F(V)}$ satisfies:
  \begin{itemize}
  \item $\bigl(\forall \Phi\in F(V)\bigr)\; \mu_{\Phi}$ is a $C$-valued measure on $V/\Phi$;
  \item Define for any $b\in \mathscr{B}(V/\Phi)$; $p^{-1}_{\Phi\Psi}(b)=\bigl\{\xi\in V/\Psi,
    p_{\Phi\Psi}(\xi)\in b \bigr\}$, and $p_{\Phi\Psi}(\mu_{\Psi})(b)=\mu_{\Psi}\bigl(p^{-1}_{\Phi\Psi}(b)
    \bigr)$. Then
    \begin{equation*}
      \bigl(\forall \, \Phi\supset \Psi\in F(V)\bigr)\; \mu_{\Phi}=p_{\Phi\Psi}(\mu_{\Psi})=p_{\Phi\Psi} \,_{*}\, \mu_{\Psi}  \; .
    \end{equation*}
  \end{itemize}
\end{definition}
\begin{remark}
  The compatibility condition of \cref{def:105} implies that for any $\Phi,\Psi\in F(V)$,
  \begin{equation*}
    \mu_{\Phi}(V/\Phi)=\mu_{\Psi}(V/\Psi)=m\in C_{\infty}\; .
  \end{equation*}
  We call $m$ the \emph{total mass} of the cylindrical measure $M$. A cylindrical measure
  $M=(\mu_{\Phi})_{\Phi\in F(V)}$ is \emph{finite} if for any $\Phi\in F(V)$, the measure $\mu_{\Phi}$ is
  finite.
\end{remark}
We recall that every $C$-valued measure $\mu$ on $V$ induces a cylindrical measure
$M_{\mu}=(\mu_{\Phi})_{\Phi\in F(V)}$. In fact, let $p_{\Phi}:V\to V/\Phi$ be the canonical projection,
then it is sufficient to choose $\mu_{\Phi}=p_{\Phi}\,_{*}\, \mu$, the push-forward (image) of $\mu$ by
means of the canonical projection. On the other hand, for any \emph{finite dimensional} $V$, every
cylindrical measure $M=(\mu_{\Phi})_{\Phi\in F(V)}$ induces a measure $\mu^{(M)}=\mu_{\{0\}}$ (this fails
to be true if $V$ is infinite dimensional). Now we are almost ready to define the Fourier transform of
cylindrical measures. In order to do that, let us recall that the push-forward can be extended to
cylindrical measures. Consider a cylindrical measure $M=(\mu_{\Phi})_{\Phi\in F(V)}$, and a linear
continuous map $u:V\to W$. For any $\Xi\in F(W)$, $u^{-1}(\Xi)\in F(V)$, and $u_{\Xi}:V/u^{-1}(\Xi)\to
W/\Xi$ is the linear application obtained from $u$ passing to the quotients. Then
\begin{equation*}
  u\, _{*}\, M=(u_{\Xi}\, _{*}\, \mu_{u^{-1}(\Xi)})_{\Xi\in F(W)}\; .
\end{equation*}
Finally, we denote by $\Phi^{\circ}\subset V'$ the polar of $\Phi$, since $\Phi$ is a vector space,
$\Phi^{\circ}=\bigl\{\xi\in V', (\forall v\in V)\xi(v)=0\bigr\}$. It is possible to identify $(V/\Phi)'$
and $\Phi^{\circ}$ by means of the adjoint map $^{\mathrm{t}}p_{\Phi}$.
\begin{definition}[Fourier transform of cylindrical measures]\label{def:106}
  Let $V$ be a topological vector space, $(X,X',C)$ a triple satisfying~\eqref{eq:101}-\eqref{eq:102}, and
  let $\mathcal{M}_{\mathrm{cyl}}(V,C)$ be the set of finite $C$-valued cylindrical measures on $V$. The
  \emph{Fourier transform} is a map $\hat{\phantom{a}}: \mathcal{M}_{\mathrm{cyl}}(V,C)\to
  (X_{\mathds{C}})^{V'}$, defined by
  \begin{equation*}
    \hat{M}(\xi)=\int_{\mathds{R}}^{}e^{2i t}  \mathrm{d}(\xi\,_{*}\, M)(t)\; .
  \end{equation*}
  The Fourier transform can be equivalently defined as
  \begin{equation*}
    (\forall \xi\in\Phi^{\circ})\; \hat{\mu}(\xi)=\int_{V/\Phi}^{}e^{2i\xi(v)}  \mathrm{d}\mu_{\Phi}(v)\; .
  \end{equation*}
\end{definition}
\noindent We remark that
\begin{equation}
  \label{eq:1021}
  V'=\bigcup_{\Phi\in F(V)}\Phi^{\circ}\; ,
\end{equation}
and the consistency condition of \cref{def:105} ensure the above definition is consistent. With the aid of
\cref{thm:102} and a projective argument, it is possible to prove the following result. A more direct
proof, for scalar measures, can be found in
\citep{vakhania1987ma
}.
\begin{thm}[Bochner for cylindrical measures]\label{thm:103}
  Let $V$ be a topological vector space, $(X,X',C)$ a triple satisfying~\eqref{eq:101}-\eqref{eq:104}. The
  Fourier transform is a bijection between finite $C$-valued cylindrical measures on $V$ and completely
  positive functions from $V'$ to $X_{\mathds{C}}$ that are ultraweakly continuous when restricted to any
  finite dimensional subspace of $V'$.
\end{thm}

\subsection{Signed and complex vector measures}
\label{sec:signed-measures}

As in the scalar case, it is possible to introduce signed and complex vector measures.
\begin{itemize}
\item Given a convex pointed cone $C$ of a real vector space $X$, we define the relation $\leq_C$ by
  \begin{equation*}
    x\leq_C  y\textrm{ iff } y-x\in C\; .
  \end{equation*}
  If for any $\bigl\{x,y\bigr\}\subset X$, there exist the supremum $x\vee_C y$ with respect to the
  partial order $\leq_C $, then $(X,\leq_C )$ is a Riesz space. In this case, we say that $C$ is a
  \emph{lattice cone} of $X$. Every pointed and generating cone is a lattice cone, if there are two
  elements in $C$ with either infimum or supremum.
\item The extended real line $\mathds{R}\cup \{-\infty,+\infty\}$ is not an additive monoid, since
  $+\infty-\infty$ is not defined. However both $(-\infty,+\infty]$ and $[-\infty,+\infty)$ are additive
  monoids.
\item We define $X_{\infty}=\mathrm{Hom}_{\mathrm{mon}}(C',\mathds{R}\cup \{-\infty,+\infty\})$ as the
  subset of functions $f\in (\mathds{R}\cup \{-\infty,+\infty\})^{C'}$ satisfying the following
  properties:
  \begin{itemize}
  \item If $\pm \infty\in \ran f$, then $\mp\infty\notin \ran f\,$;
  \item $f:C'\to \ran f$ is a monoid homomorphism.
  \end{itemize}
  This definition is justified by the fact that since $+\infty-\infty$ is not defined, signed measures may
  only take either $+\infty$ or $-\infty$ as a value (in order to be additive). This has also to be the
  case for signed vector measures, and therefore they will have $X_{\infty}$ as target space, see
  \cref{def:14} below.
\item $(X_{\mathds{C}})_{\infty}=\mathrm{Hom}_{\mathrm{mon}}(C',\mathds{C}\cup \{\infty\})$.
\end{itemize}
Let us consider the extension to vector measures of the concept of signed measures. This is easily done by
means of $X_{\infty}$ defined above.
\begin{definition}[Signed vector measures]\label{def:14}
  Let $(X,X',C)$ be a triple that satisfies \eqref{eq:101}-\eqref{eq:102}, $E$ a measurable space.  A
  function $\mu\in (X_{\infty})^{\Sigma}$ is a signed vector measure on $E$ iff it is countably additive
  and $\mu(\emptyset)=0$.
\end{definition}
The following useful lemma follows directly from the definition of signed measures.
\begin{lemma}\label{lemma:101}
  Every $C$-valued measure is also a signed measure. Any real linear combination of two $C$-valued
  measures is a signed measure, provided at least one of the two measures is finite.
\end{lemma}
The important \cref{thm:101} can be easily adapted to hold for signed measures as well.
\begin{proposition}\label{thm:8}
  There is a bijection between signed vector measures $\mu$ on $E$ and families of signed measures
  $(\mu_{\kappa})_{\kappa\in C'}$ on $E$ such that for any $b\in \Sigma$, $\bigl(\kappa\mapsto
  \mu_{\kappa}(b) \bigr)\in X_{\infty}$.
\end{proposition}
A signed measure $\mu$ is \emph{finite} iff for any $\kappa\in C'$, $\mu_{\kappa}$ is finite. The idea
behind signed vector measures is that, as in the case of standard measures, they are the sum of two
cone-valued measures. Therefore it is reasonable to define them as a collection indexed only by the dual
cone $C'$, in order to prevent possible ``sign incongruences'' on $\mu_{\kappa}$ due to the action of a
$\kappa\notin C'$.  As a matter of fact, with this definition we can indeed prove the existence of a
unique Jordan decomposition for signed vector measures. The precise statement is contained in the
following result.
\begin{proposition}\label{thm:9}
  Let $(X,X',C)$ be a triple satisfying \eqref{eq:101}-\eqref{eq:102}; and $\mu$ a signed vector measure
  on a measurable space $E$. Then there exist three $C$-valued measures $\mu^+,\mu^-,\lvert \mu
  \rvert_{}^{}$ such that: 
  \begin{itemize}
  \item $\mu=\mu^+-\mu^-$, and the decomposition is unique;
  \item $\lvert \mu \rvert_{}^{}=\mu^++\mu^-$;
  \item At least one between $\mu^+$ and $\mu^-$ is finite;
  \item $\mu$ is finite iff $\lvert \mu \rvert_{}^{}$ is finite.
  \end{itemize}

  If in addition $(X,\leq _C)$ is a Riesz space, then $\mu^+=\mu\vee_C 0$, $\mu^-=\mu\wedge_C0$ and
  $\lvert \mu \rvert_{}^{}=\lvert \mu \rvert_C^{}$.  The operations $+$, $-$, $\vee_C$, $\wedge_C$ and
  $\lvert\, \cdot \,\rvert_C^{}$ on measures are defined pointwise on measurable sets, and $0$ is the
  measure identically zero.
\end{proposition}
\begin{proof}
  Let $\mu$ be a signed vector measure. Then $(\mu_{\kappa})_{\kappa\in C'}$ is the corresponding family
  of signed measures. By Jordan decomposition of signed measures, for any $\kappa\in C'$, there exist a
  unique decomposition $\mu_{\kappa}=\mu_{\kappa}^+-\mu_{\kappa}^-$, with $\mu_{\kappa}^+$ and
  $\mu_{\kappa}^-$ positive measures with at least one of the two finite, and $\mu_{\kappa}$ is finite iff
  $\lvert \mu_{\kappa} \rvert_{}^{}$ is finite. Hence if $(\lvert \mu_k \rvert_{}^{})_{\kappa\in C'}$ is
  the image of a $C$-valued measure $\lvert \mu \rvert_{}^{}$, $\mu$ is finite iff $\lvert \mu
  \rvert_{}^{}$ is finite. In addition, suppose that there exists a $\tilde{\kappa}\in C'$ such that
  $\mu_{\kappa}^+$ is not finite. Then $+\infty\in \ran \mu$, and therefore $-\infty\notin \ran \mu$,
  i.e. for any $\kappa\in C'$, $\mu^-_{\kappa}$ is finite. It follows that if $(\mu_{\kappa}^-)_{\kappa\in
    C'}$ is the image of a $C$-valued measure, such measure is finite. An analogous statement holds with
  plus replaced by minus. By Lemma~\ref{lemma:101}, to prove the first part of the theorem it remains only
  to check that the families $(\mu_{\kappa}^+)_{\kappa\in C'}$ and $(\mu_{\kappa}^-)_{\kappa\in C'}$ are
  $C$-valued measures, i.e. that for any $b\in \Sigma$, the maps $\kappa\mapsto \mu^{\pm}_{\kappa}(b)$ are
  monoid morphisms. On one hand, we have by the fact that $\mu\in X_{\infty}$ and then Jordan
  decomposition that
  \begin{equation*}
    \mu_{\kappa_1+\kappa_2}(b)=\mu_{\kappa_1}(b)+\mu_{\kappa_2}(b)=\mu_{\kappa_1}^+(b)+\mu_{\kappa_2}^+(b)-\bigl(\mu_{\kappa_1}^-(b)+\mu_{\kappa_2}^-(b)\bigr)\; ;
  \end{equation*}
  on the other hand, by Jordan decomposition we have also that
  \begin{equation*}
    \mu_{\kappa_1+\kappa_2}(b)=\mu_{\kappa_1+\kappa_2}^+(b)-\mu_{\kappa_1+\kappa_2}^-(b)\; .
  \end{equation*}
  Now since the decomposition is unique, it follows that
  \begin{equation*}
    \mu_{\kappa_1+\kappa_2}^{\pm}(b)=\mu_{\kappa_1}^{\pm}(b)+\mu_{\kappa_2}^{\pm}(b)\; ,
  \end{equation*}
  i.e. the map is a monoid morphism.

  To prove the last part, let $\mu=\mu^+-\mu^-$ be a vector signed measure with the respective
  decomposition. Then for any $b\in \Sigma$, we have that
  \begin{equation*}
    X\ni \mu(b)=\mu^+(b)-\mu^-(b)\; ;\; \mu^+(b),\mu^-(b)\geq_C 0\; .
  \end{equation*}
  If $(X,\leq_C )$ is a Riesz space, it then follows that $\mu^+=\mu\vee_C 0$, $\mu^-=\mu\wedge_C0$.
\end{proof}

The complex vector measures are defined in an analogous fashion, and they are the sum of four $C$-valued
measures. We quickly mention the basic definitions and results without proof, for they are equivalent to
the ones for signed vector measures.
\begin{definition}[Complex vector measures]\label{def:15}
  Let $(X,X',C)$ be a triple that satisfies \eqref{eq:101}-\eqref{eq:102}, $E$ a measurable space. A
  function $\mu\in \bigl((X_{\mathds{C}})_{\infty}\bigr)^{\Sigma}$ is a complex vector measure on $E$ iff
  it is countably additive and $\mu(\emptyset)=0$.
\end{definition}
\begin{lemma}\label{lemma:8}
  Under the identifications $\mathds{R}\ni \alpha\to \alpha+i0$, $+\infty\to \infty$, $-\infty\to \infty$; every signed vector measure is also
  a complex measure. Any complex linear combination of two signed measures is a complex measure.
\end{lemma}
\begin{proposition}\label{thm:10}
  Let $(X,X',C)$ be a triple satisfying \eqref{eq:101}-\eqref{eq:102}; and $\mu$ a complex vector measure
  on a measurable space $E$. Then there exist five $C$-valued measures
  $\mu^+_{\mathrm{R}},\mu^-_{\mathrm{R}},\mu^+_{\mathrm{I}},\mu^-_{\mathrm{I}},\lvert \mu \rvert_{}^{}$
  such that:
  \begin{itemize}
  \item $\mu=\mu^+_{\mathrm{R}}-\mu^-_{\mathrm{R}}+i(\mu^+_{\mathrm{I}}-\mu^-_{\mathrm{I}})$, and the
    decomposition is unique;
  \item $\lvert \mu
    \rvert_{}^{}=\mu^+_{\mathrm{R}}+\mu^-_{\mathrm{R}}+\mu^+_{\mathrm{I}}+\mu^-_{\mathrm{I}}$;
  \item $\mu$ is finite iff $\lvert \mu \rvert_{}^{}$ is finite, or equivalently if
    $\mu^+_{\mathrm{R}},\mu^-_{\mathrm{R}},\mu^+_{\mathrm{I}},\mu^-_{\mathrm{I}}$ are all finite.
  \end{itemize}
\end{proposition}
\begin{corollary}
  The integral with respect to a \emph{finite} complex vector measure $\mu$ is a map $\int_{(\cdot )} ^{}
  \mathrm{d}\mu:\Sigma\to X_{\mathds{C}}$ defined by
  \begin{equation*}
    \int_b^{}  \mathrm{d}\mu=\int_b^{}  \mathrm{d}\mu^+_{\mathrm{R}}-\int_b^{}  \mathrm{d}\mu^-_{\mathrm{R}}+i\Bigl(\int_b^{}  \mathrm{d}\mu^+_{\mathrm{I}} -\int_b^{}  \mathrm{d}\mu^-_{\mathrm{I}} \Bigr)\; .
  \end{equation*}
\end{corollary}
\subsection{A concrete realization: duals of C*-algebras}
\label{sec:class-concr-exampl}

In this subsection, we discuss explicitly the relevant class of triples satisfying the properties
\eqref{eq:101}-\eqref{eq:104} that was used in this paper.
\begin{itemize}
\item Given a C*-algebra $\mathfrak{A}$, we denote by $\mathfrak{A}_+$ the set of elements with positive
  spectrum, and by $\mathfrak{A}_{*}$ the set of self-adjoint elements.
\item If $\mathfrak{A}'_{*}$ is the continuous dual of the set of self-adjoint elements $\mathfrak{A}_{*}$
  of a C*-algebra $\mathfrak{A}$, we denote by $\mathfrak{A}'_+$ the functionals that are positive when
  acting on $\mathfrak{A}_+$.
\end{itemize}
In order to verify conditions \eqref{eq:101}-\eqref{eq:104}, we make use of the following classical result
\citep[see, \emph{e.g.},][]{takesaki1979I
}.
\begin{proposition}\label{thm:11}
  Let $\mathfrak{A}$ be a C*-algebra. Then:
  \begin{itemize}
  \item $\mathfrak{A}_{*}$ is a real Banach subspace of $\mathfrak{A}$ and
    $\mathfrak{A}=\mathfrak{A}_{*}+i\mathfrak{A}_{*}$.
  \item $\mathfrak{A}_+$ is a closed, pointed and generating convex cone of $\mathfrak{A}_{*}$.
  \item $\mathfrak{A}'_+$ is a pointed and generating convex cone of $\mathfrak{A}'_{*}$; in particular
    for any $\alpha\in \mathfrak{A}'_{*}$ there is a unique decomposition
    \begin{equation*}
      \alpha=\alpha^+-\alpha^-\;,\;\textrm{ with } \alpha^+,\alpha^-\in \mathfrak{A}'_+\; .
    \end{equation*}
  \item $(\mathfrak{A}_*')_{\mathds{C}}=\mathfrak{A}'$.
  \end{itemize}
\end{proposition}

By means of \cref{thm:11}, conditions~\eqref{eq:101}, \eqref{eq:103}-\eqref{eq:104} are immediately
proved. Condition~\eqref{eq:102} is proved using a remarkable result of \citet[][Lemma
I.5]{neeb1998mm
}. In fact, if we call $C'_1$ the set of elements of $C'$ with $\mathfrak{A}_{*}$-norm one, then
$C'_1-C'_1$ is a $0$-neighbourhood of $\mathfrak{A}_{*}$. In \cref{sec:class-char-stat}, the
$\mathfrak{A}_+'$-valued measures played an important role; from the discussion above, it follows that the
usage of all the results in this appendix, and especially Bochner's theorem, is justified.

\section*{Acknowledgments}
\label{sec:acknowledgements}

The author acknowledges the support of MIUR through the FIR grant 2013
\href{http://www.cond-math.it}{Condensed Matter in Mathematical Physics (Cond-Math)} (code
RBFR13WAET). Part of this work has been done at the University of Stuttgart, with the support of the IRTG
\href{http://www.mathematik.uni-stuttgart.de/grk1838/}{GRK1838}.

{\footnotesize
  
\Addresses}
\end{document}